\documentclass[12pt,leqno]{article}
\usepackage{amsmath,amsthm,amssymb}

\def\init{\setcounter{equation}{0}}
\newtheorem{theorem}{Theorem}[section]
\newcommand{\R}{\mathbb{R}}
\newcommand{\Z}{\mathbb{Z}}

\newtheorem{lemma}[theorem]{Lemma}

\newcommand{\e}{{\varepsilon}}

\newcommand{\rw}{\rightarrow}

\usepackage{tikz}
\usetikzlibrary{decorations.pathreplacing}

\usepackage{comment}

\title{Inverse problems for general second order  hyperbolic equations with time-dependent coefficients}

\author{G.Eskin, \ \ \  Department of Mathematics, UCLA,\\ Los Angeles,
CA 90095-1555, USA. \ E-mail: eskin@math.ucla.edu
}

\begin{document}
\maketitle

\begin{abstract}
We study the inverse problems for the second order hyperbolic equations of general form with time-dependent coefficients
assuming that the boundary data are given on a part of the boundary. 
The main  result  of this paper  is the determination  of the time-dependent Lorentzian metric by the boundary measurements.
This is achieved by the adaptation  of a variant of the Boundary Control method  developed by the author in [E1], [E2].  
\end{abstract}

Keywords:  inverse  problems,  hyperbolic equations,  geometric optics.

Mathematics Subject Classification 2010: 35R30,  35L10

\section{Introduction.}
\label{section 1}
\init
Consider a second order 
hyperbolic equation in $\R^{n+1}$ of the form
\begin{multline}                                                                                           \label{eq:1.1}
\sum_{j,k=0}^n\frac{1}{\sqrt{(-1)^ng(x)}}\Big(-i\frac{\partial}{\partial x_j}-A_j(x)\Big)
\sqrt{(-1)^ng(x)}g^{jk}(x)\Big(-i\frac{\partial}{\partial x_k}-A_k(x)\Big)u(x)
\\
=0,
\end{multline}
where
$x=(x_0,x_1,...,x_n)\in \R^{n+1},\ x_0$ is the time variable.
 In (\ref{eq:1.1})  $g(x)=\det[g_{jk}(x)]_{j,k=0}^n$,  where $[g_{jk}(x)]_{j,k=0}^n=\big([g^{jk}]_{j,k=0}^n\big)^{-1}$
 is the metric tensor,  $A(x)=(A_0(x),A_1(x),...,A_n(x))$  is  the vector potential.  We assume  that all coefficients in (\ref{eq:1.1}) 
 belong to 
 $C^\infty(\R^{n+1})$  and that
 \begin{equation}                                                                                           \label{eq:1.2}
 g^{00}(x)\geq c_0>0,\ \ \ \forall x\in \R^{n+1}.
 \end{equation}   
Let $(\xi_0,\xi_1,...,\xi_n)$  be dual 
variables 
to $(x_0,x_1,...,x_n)$.  The strict hyperbolicity  of (\ref{eq:1.1})  with respect to $\xi_0$  means
that the quadratic equation 
\begin{equation}                                             								\label{eq:1.3}
\sum_{j,k=0}^n g^{jk}(x)\xi_j\xi_k=0
\end{equation}
has two real distinct roots $\xi_0^-(\xi_1,...,\xi_n)<\xi_0^+(\xi_1,...,\xi_n)$  for all $(\xi_1,...,\xi_n)\neq(0,...,0)$  and all $x\in \R^{n+1}$.
We have
\begin{multline}           														\label{eq:1.4}
\xi_0^\pm(\xi_1,...,\xi_n)
\\
=\frac{-\sum_{j=1}^ng^{j0}(x)\xi_j
\pm\sqrt{\big(\sum_{j=1}^ng^{j0}(x)\xi_j\big)^2-g^{00}(x)\sum_{j,k=1}^ng^{jk}(x)\xi_j\xi_k}}{g^{00}}.
\end{multline}
The strict hyperbolicity implies that
\begin{equation}             															\label{eq:1.5}
\Big(\sum_{j=1}^n g^{j0}(x)\xi_j\Big)^2-g^{00}(x)\sum_{j,k=1}^ng^{jk}(x)\xi_j\xi_k>0
\end{equation}
for all $(\xi_1,...,\xi_n)\neq 0, \ x\in \R^n$.

In this paper we assume a more restrictive condition that
\begin{equation}															\label{eq:1.6}
\sum_{j,k=1}^ng^{jk}(x)\xi_j\xi_k\leq -c_1\sum_{j=1}^n\xi_j^2,
\end{equation}
i.e.  we assume that the spatial part of the equation (\ref{eq:1.1}) is elliptic for any $x\in\R^{n+1}$.

Note that  the quadratic form  (\ref{eq:1.3}) has the signature  $(+1,-1,...,-1)$.  Therefore  $(-1)^ng(x)>0$.
We assume also that $A_j(x),0\leq j\leq n,$  are real-valued.  Thus the operator  in (\ref{eq:1.1}) is formally self-adjoint.

We consider the following class of domains $D\subset\R^{n+1}$.
Let $D_t=D\cap\{x_0=t\}$  be  the intersection of  $D$  with  the plane $\{t=x_0\}, t\in\R$.  We assume  that $ D_t$  is
a smooth closed bounded domain in $\R^n$  smoothly dependent and uniformly bounded in $t$  and such  that
$D_t$  is diffeomorphic   to  $D_0$  for  all $t\in\R$.   More  precisely   we assume  that there  exists  a diffeomorphism
\begin{equation}																\label{eq:1.7}
y_0=x_0,\ \ \ y_j=\hat y_j(x_0,x_1,...,x_n),\ 1\leq j\leq n,
\end{equation}
that maps  $D_{x_0}$  onto $D_0$ and smoothly depends on $x_0$.   We shall  call  such  domains $D$  admissible.

Let $S(x_0,x_1,...,x_n)=0$  be the equation of $\partial D=\bigcup_{t\in \R}\partial D_t$.  
Sometimes we denote  $\partial D$  by  $S$.
We assume that $S$  is   a time-like 
smooth surface in $\R^{n+1}$,  i.e.  
\begin{equation}             	    													\label{eq:1.8}
\sum_{j,k=0}^ng^{jk}(x)\nu_j\nu_k<0,
\end{equation}
where $x\in S$  and $(\nu_0,\nu_1,...,\nu_n)$  is a normal vector  to   
$S$.
The vector $(\nu_0,...,\nu_n)$ satisfying (\ref{eq:1.8}) is called a space-like vector.  Also,  the surface  $\Sigma$  in  $\R^{n+1}$ 
is called space-like  if 
$\sum_{j,k=0}^ng^{jk}(x)\nu_j(x)\nu_k(x)>0$,  where  $x\in \Sigma$ and $(\nu_0(x),...,\nu_n(x))$  is the normal vector
to $\Sigma$.

Consider the initial-boundary value problem
\begin{align}																\label{eq:1.9}
&Lu = 0 \ \ \ \mbox{in}\ \ D,
\\
&u = 0 \ \ \ \mbox{for}\ \ x_0 \ll 0\ \ \mbox{in}\ \ D,								\label{eq:1.10}
\\
&u\big|_S =f,																\label{eq:1.11}
\end{align}
where
$f$ is a smooth function on $S=0$  with compact support,  $Lu=0$  is the same as in (\ref{eq:1.1}).

It is well known that the initial-value problem (\ref{eq:1.9}), (\ref{eq:1.10}),(\ref{eq:1.11})
is well-posed (cf. [H]),  assuming that (\ref{eq:1.2}),  (\ref{eq:1.5}) and (\ref{eq:1.8})
are satisfied.

Let $S_0\subset S$  be a part of $S$  such that (see Fig. 1.1)  $S_{0t}=\partial D_t\cap S_0$  has a nonempty interior for
all $t\in \R$.  We assume also that for any $x^{(0)}\in \partial S_0$  a vector  $\tau^{(1)}$,  tangent to $S$ and normal to $S_0$,
is not parallel to $(1,0,...,0)$.

We include in the definition of the admissibility  of  $D$  (see  (\ref{eq:1.7}))   that the map  $y=\hat y(x_0,x)$   is such  that 
\begin{equation} 																		\label{eq:1.12}
\hat y(x_0,x)=x  \ \ \ \mbox{on}\ \ \ S_{0x_0},\ \forall x_0\in \R.
\end{equation}
\\
\

\begin{tikzpicture}
\draw (-2,0) arc (180:360: 2 and 1);
\draw[dashed] (2,0) arc (0:180:2 and 1);

\draw (-0.5,3) .. controls (-1,2) and (-2,1) .. (-2,0);
\draw (-2,0).. controls (-2,-1) and (-1,-1.8) .. (-0.6,-3);

\draw (0.5,3) .. controls (1,2) and (2,1) .. (2,0);
\draw (2,0).. controls (2,-1) and (1,-1.8) .. (0.6,-3);

\draw [line width=3pt](-0.3,3) .. controls (-0.8,1.8) and (-1.5,1) .. (-1.1,-0.8);
\draw[line width=3pt] (-1.1,-0.8).. controls (-1,-1) and (-1,-1.8) .. (-0.2,-3);

\draw[line width=3pt] (0.3,3) .. controls (0.8,1.8) and (1.5,1) .. (1.1,-0.8);
\draw [line width=3pt](1.1,-0.8).. controls (1,-1) and (1,-1.8) .. (0.2,-3);

\draw (0.9,0.5) node {$ S_0$};
\draw (-2.5, -0.5)  node {$\partial D_t$};

\draw (1.5,2) node {$S$};

\end{tikzpicture}   
\\
\
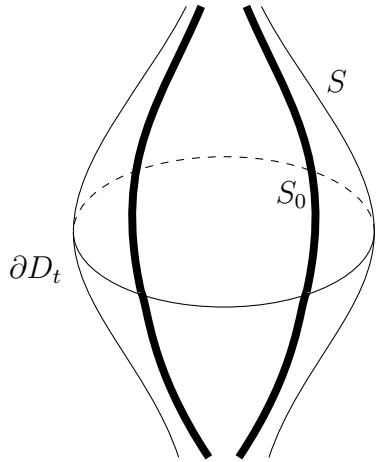
{\bf Fig. 1.1.}  $S_0$ is  part of the boundary $S, \ S_0\cap \partial D_t\neq \emptyset$  for $\forall t$. 
\\
\

Note  that  (\ref{eq:1.7})   maps  the admissible domain $D$   onto  $D=D_0\times\R,\ S_0$  
onto  $\hat S_0=\Gamma_0\times\R$  (cf.  (\ref{eq:1.12})),  where $\Gamma_0=S_0\cap \partial D_0$,  i.e.   $\hat  D,\hat S_0$ 
are cylindrical  domains (see  Fig.  1.2).
\\
\\
\

\begin{tikzpicture}

\draw (-2,0) arc (180: 360 :2 and 1);
\draw[dashed] (2,0) arc(0:180:2 and 1);
\draw (-2,-3) -- (-2,2.7);
\draw (-0.5,-3.2) -- (-0.5,2.2);
\draw (1.5,-3) -- (1.5,2.3);
\draw (2,-2.7) -- (2,2.7);

\draw (-2.5,-0.3)  node {$\hat D_0$};
\draw (1,-1.3) node {$\Gamma_0$};

\draw [ultra thick] (-0.5,-1) .. controls (0,-1.1)  and (1,-1) ..  (1.5,-0.7);

\end{tikzpicture}

\
\\
{\bf Fig. 1.2.}
\ \
$\hat D=\hat D_0\times \R, \ \ \hat S_0=\Gamma_0\times\R$ are cylindrical domains.
\\
\

The Dirichlet-to-Neumann operator $\Lambda$  that maps the Dirichlet data to the Neumann 
data on $\partial D$  is defined as
\begin{equation} 																	\label{eq:1.13}
\Lambda f=\sum_{j,k=0}^n g^{jk}(x)\Big(\frac{\partial u}{\partial x_j}-iA_j(x)u\Big)\nu_k(x)
\Big(\sum_{p,r=0}^n g^{pr}(x)\nu_r\nu_p\Big)^{-\frac{1}{2}}\Big|_S,
\end{equation}
where  $u(x)$  is the solution  of (\ref{eq:1.9}), (\ref{eq:1.10}),(\ref{eq:1.11})
and $(\nu_0(x),...,\nu_n(x))$  is the unit  outward normal to $S$.

Denote by
\begin{equation}															             \label{eq:1.14}
y=y(x)
\end{equation}
any proper diffeomorphism  of $\overline D$  onto   some domain $\overline{\hat D}$ such that 
\begin{equation}                  			  											    \label{eq:1.15}
y=x\ \ \ \mbox{on}\ \ S_0.
\end{equation}
We call a diffeomorphism  of $\overline D$ onto  $\overline {\hat D}$ proper if for  any  $[t_1,t_2]\subset \R$  the image
of  $\overline D\cap \{t_1\leq x_0\leq t_2\}$  is a domain $\overline{\hat D}\cap\{S_-^{(t_1)}(y_1,...,y_n)\leq y_0
\leq S_+^{(t_2)}(y_1,...,y_n)\}$,  where  $y_0=S_+^{(t_1)},y_0=S_-^{(t_0)}$  are  space-like surfaces.

Let
$\hat L\hat u=0$  be the equation  (\ref{eq:1.1})  in $y$-coordinates,  $y\in \hat D$.  We have
\begin{multline}             															\label{eq:1.16}
\hat L\hat u
\equiv \sum_{j,k=0}^n\frac{1}{\sqrt{(-1)^n\hat g(y)}}
\Big(-i\frac{\partial}{\partial y_j}-\hat A_j(y)\Big)\sqrt{(-1)^n\hat g(y)}
\hat g^{jk}(y)
\\
\cdot 
\Big(-i\frac{\partial}{\partial y_k}-\hat A_k(y)\Big)
\hat u=0,
\end{multline}
where
\begin{equation}              														\label{eq:1.17}
\hat g^{jk}(y)=\sum_{p,r=0}^ng^{pr}(x)\frac{\partial y_j}{\partial x_p}\frac{\partial y_k}{\partial x_r},
\end{equation}
\begin{equation}                                 												\label{eq:1.18}
A_j(x)=\sum_{k=0}^n\hat A_k(y)\frac{\partial y_k}{\partial x_j},\ \ \  0\leq j\leq n,
\end{equation}
Here   
\\
$\hat u(y)=u(x)$,  
$y=y(x),\ \hat g(y)=\det [\hat g_{jk}(y)]_{j,k=0}^n,\ [\hat g_{jk}(y)]_{j,k=0}^n=\big([\hat g^{jk}(y)]_{j,k=0}^n\big)^{-1}$.

Note  that (\ref{eq:1.17}), (\ref{eq:1.18}) are equivalent to the equalities
\begin{align}        					  					   	    				\label{eq:1.19}
\sum_{k=0}^nA_k(x)dx_k &=\sum_{k=0}^n\hat A_k(y)dy_k,
\\ 
																			\label{eq:1.20}
\sum_{j,k=0}^n\hat g_{jk}(y)dy_jdy_k &=\sum_{j,k=0}^n g_{jk}(x)dx_jdx_k,
\end{align}
where $y$  and $x$  are related by (\ref{eq:1.14}).
Metric tensors $[g_{jk}(x)]_{j,k=0}^n$  and  $[\hat g_{jk}(x)]_{j,k=0}^n$, related  by (\ref{eq:1.20}),
are called isometric.

We assume that conditions (\ref{eq:1.2}),  (\ref{eq:1.6}) hold also in $y$-coordinates,  i.e.
\begin{equation}                														\label{eq:1.21}
\hat g^{00}(y)\geq C_0,\ \ \ \sum_{j,k=1}^n\hat g^{jk}(y)\xi_j\xi_k\leq -C_1\sum_{j=1}^n\xi_j^2.
\end{equation}       

Let $c(x)\in C^\infty(\overline D)$  
be such that
\begin{equation}              															\label{eq:1.22}
|c(x)|=1,\ \ x\in \overline D,
\ \ \ \ 																
c(x)=1\ \ \mbox{on}\  \ S_0.
\end{equation}
The group $G_0(\overline D)$ of such $c(x)$  is called the gauge group.

If $Lu=0$  then  $u'=c^{-1}(x)u(x)$  satisfies the equation of the form (\ref{eq:1.1})  with $A_j(x)$
replaced   by
\begin{align}																		\label{eq:1.23}
A_j'(x)&=A_j(x)-ic^{-1}(x)\frac{\partial c}{\partial x_j}, \ \  1\leq j\leq n,
\\
\nonumber
A_0'(x)&=A_0(x)+ic^{-1}(x)\frac{\partial c}{\partial x_0}.
\end{align}

We shall call potentials $(A_0',..,A_n'(x))$  and $(A_0(x),..,A_n(x)))$  related by (\ref{eq:1.23}) gauge  equivalent.  Note  that when 
$D$  is simply connected  then  $c(x)=\exp i\varphi$  where $\varphi(x)\in C^\infty(\overline D),\ \varphi(x)$  is
real-valued and  $\varphi(x)=0$  on $S_0$.

Let  $y=y(x)$  be  the change  of variables,  such that  $y(x)=x,\ x\in S_0$,  transforming  the equation $Lu=0$  in $D$ 
 to   the equation  
of the form  (\ref{eq:1.16})  in $\hat D$.  We  consider  the initial-boundary value problem
\begin{align}                                    													\label{eq:1.24}
&\hat L\hat u=0 \ \ \mbox{in}\ \ \hat D,
\\
                                                                                                                                   \label{eq:1.25}
&\hat u = 0 \ \ \mbox{for}\ \ y_0\ll 0,\ y\in \hat D,
\\
																			\label{eq:1.26}
&\hat u\big|_{\hat S}=f.
\end{align}
Note that $\hat S_0=S_0$,  since  $\hat y(x)=x$  on $S_0$.

Since  (\ref{eq:1.21})  holds, 
the initial-boundary value problem (\ref{eq:1.24}), (\ref{eq:1.25}), (\ref{eq:1.26})  is also well-posed.  Let  $\hat c(y)\in G_0(\overline{\hat D})$.
Make  the gauge transformation $u'(y)=\hat c^{-1}(y)\hat u(y)$  and  let $L'$  be such that $L'u'=0$.  We have 
\begin{align}                                    													\label{eq:1.27}
&\hat L'u'=0 \ \ \mbox{in}\ \ \hat D,
\\
\label{eq:1.28}
& u' = 0 \ \ \mbox{for}\ \ y_0\ll 0,\ y\in \hat D,
\\
\label{eq:1.29}
& u'\big|_{\hat S}=f.
\end{align}
Note that  $u'=\hat u$  on $\hat S_0$  since  $\hat c(y)=1$  on $\hat S_0$  and $L'u'$  has  the form
\begin{multline}                                                                                           \label{eq:1.30}
L'u'
\\
=\sum_{j,k=0}^n\frac{1}{\sqrt{(-1)^n\hat g(y)}}\Big(-i\frac{\partial}{\partial y_j}-A_j'(y)\Big)
\sqrt{(-1)^n\hat g(y)}\hat g^{jk}(y)\Big(-i\frac{\partial}{\partial y_k}-A_k'(y)\Big)u'(y)
\\
=0,
\end{multline}
$A_j'(y),\ 0\leq j\leq n,$  are potentials gauge equivalent 
to $\hat A_j(y),\ 0\leq j\leq n$.

Let  $\Lambda'$  be the  DN operator for  
(\ref{eq:1.27}), (\ref{eq:1.28}), (\ref{eq:1.29}) 
\begin{equation} 																	\label{eq:1.31}
\Lambda' f=\sum_{j,k=0}^n \hat g^{jk}(y)\Big(\frac{\partial u'}{\partial y_j}-iA_j'(y)u'\Big)\nu_k(y)
\Big(\sum_{p,r=0}^n \hat g^{pr}(y)\nu_r(y)\nu_p(y)\Big)^{-\frac{1}{2}}\Big|_{\hat S},
\end{equation}
where $f$ is the same  as in (\ref{eq:1.10})  and  (\ref{eq:1.29}).

It can be shown  that 
\begin{equation}  																	\label{eq:1.32}
\Lambda f\big|_{S_0}=\Lambda' f\big|_{S_0},\ \ \forall f\in C_0^\infty(S_0),
\end{equation}
if  the operator $L'$  is obtained from $L$  by the change of variables (\ref{eq:1.14}),  (\ref{eq:1.15})  and the gauge 
 transformation $c(y)$  such that
  (\ref{eq:1.22})  holds.

Therefore the inverse problem
of the determination of the coefficient of
   (\ref{eq:1.1})  can be solved only up to the  changes of variables  
       (\ref{eq:1.14}),  (\ref{eq:1.15})   and the gauge transformations   (\ref{eq:1.22}).

We shall formulate now some conditions which will be required to solve the inverse problem.
\\
\

\underline{1)     Real analyticity  in the time variable}
\\
\\
One of the crucial steps  in solving  the inverse problem  will be the use of the following unique continuation theorem of Tataru  and
 Robbiano and Zuily (cf [T], [RZ])  that requires the analyticity in $x_0$:
 \begin{theorem}																				\label{theo:1.1}
 Let the coefficients of (\ref{eq:1.1}) 
be analytic in $x_0$.   Consider  the equation   $Lu=0$  in a neighborhood $U_0$  of a point $P_0$.
Let  $\Sigma=0$  be  a noncharacteristic surface  with respect  to $L$ passing through 
$P_0$.  If  $u=0$  in  $U_0\cap\{\Sigma <0\}$
then $u=0$  in $U_0\cap \{\Sigma>0\}$  near $\Sigma=0$.
\end{theorem}
We assume  also  that  the gauge  $c(x)$  and  the  map  (\ref{eq:1.7})  are analytic  in $x_0$.

  Let  $y=\varphi(x)$  be a diffeomorphism  of neighborhood  $U_0$  onto the neighborhood $\overline V_0=\varphi(\overline U_0)$.  
  Here $\varphi(x)$  is smooth  but not analytic  in any variable.  It is  clear  that  if the unique continuation property   for 
  the operator $L$   holds in $U_0$  then it holds in $V_0$  for  the operator  $\tilde L=\varphi\circ L$,  though the coefficients
  of  $\tilde L$  are  not analytic. 
  Here $\varphi\circ L$   is the operator  $L$  in  $y$-coordinates  (cf.  (\ref{eq:1.16})).
   Therefore  
  the following 
  more general class  of operators $L$  with non-analytic coefficients has the unique 
  continuation  property:  For each point  $x^{(0)}$  on $D$  there is a neighborhood  $U_0$  and  the diffeomorphism 
  $\psi(x)$  of $U_0$  onto  $V_0=\psi(U_0)$  such that  the coefficients of the operators $\psi\circ L$  in  $V_0$  are analytic  in
  $x_0$.  Thus, the unique  continuation  property holds  for $L$  in $U_0$. 
\\
\

\underline{2)     The Bardos-Lebeau-Rauch condition}
\\
\\
Consider the initial-boundary  value problem
$$
Lu=0,\ u=0  \ \ \mbox{for}\ \ x_0\ll 0,\ u\big|_{\partial D_0\times \R}=f
$$
 in the cylindrical domain $D_0\times\R,\ f$  has a compact  support in $\Gamma_0\times\R,\ 
 \Gamma_0\subset \partial D_0$.  We say that
 BLR  condition holds
 on $[t_0,T_{t_0}]$  if  
the bounded map from 
$f\in H_1(\Gamma_0\times(t_0,T_{t_0}))$ to 
$\big(u\big|_{x_0=T_{t_0}},\frac{\partial u}{\partial x_0}\big|_{x_0=T_{t_0}}\big)\in
H_1(D_0)\times L_2(D_0)$,
is onto  in $H_1(D_0)\times L_2(D_0)$,
where $u=0$  for  $x_0<t_0,\ f=0$  for $x_0< t_0$.

Note that BLR  condition obviously holds on  $[t_0,T]$  for any  $T>T_{t_0}$  if it holds  on $[t_0,T_{t_0}]$.  

Let $\{x=x(s),\xi=\xi(s)\}\in T_0^*(\overline D_0\times[t_0,T_{t_0}])$,  where 
\begin{align}                   														\label{eq:1.33}
&\frac{dx_j}{ds}=\frac {\partial  L_0(x(s),\xi(s))}{\partial \xi_j},\ \ x_j(0)=y_j, \ \ 0\leq j\leq n,
\\
\nonumber
&\frac{d\xi_j}{ds}=-\frac {\partial  L_0(x(s),\xi(s))}{\partial x_j},\ \ \xi_j(0)=\eta_j, \ \ 0\leq j\leq n,
\end{align}
be the equations of null-bicharacteristics.
Here $L_0(x,\xi)=\sum_{j,k=0}^ng^{jk}(x)\xi_j\xi_k,\\ L_0(y,\eta)=0$.

We assume  that  for any  $t_0$  there exists 
$T_{t_0}$   depending  continuously  on  $t_0$  such that  the  BLR  condition  holds  on $[t_0,T_{t_0}]$.                                                                                                                                                                                                                                                                                                                                                                                                     
    It follows  from  [BLR]     that BLR  
  condition  holds   if   any  null bicharacteristic in
$T_0^*(\overline D_0\times[t_0,T])$  intersects $T_0^*(\overline \Gamma_0\times[t_0,T])$  when $T\geq T_{t_0}$. 
\\
\

\underline{3)     Domains of dependence}
\\
\\
Let $G(x,\xi)=\sum_{j,k=0}^ng_{jk}(x)\xi_j\xi_k,\ [g_{jk}]_{j,k=0}^n=\big([g^{jk}]_{j,k=0}^n\big)^{-1}$.
We say  that  $x=x(\tau)$ is a forward time-like ray in $D_0\times\R$ if $x=x(\tau)$  is piece-wise smooth,  
$G(x(\tau),\frac{d x(\tau)}{d\tau})>0$
and $\frac{d x_0}{d\tau}>0,\ 0\leq \tau$.  If
$G(x(\tau),\frac{d x(\tau)}{d\tau})>0$ and $\frac{d x_0}{d\tau}<0$  the ray $x=x(\tau)$  is  called  the backward  time-like ray.

One can show  (cf [CH])  that the forward domain  of influence  $D_+(F)$  of a closed set $F\subset  D_0\times\R$  is 
the closure  of the union of all piece-wise smooth  forward time-like rays in $D_0\times \R$  starting  on $F$.

Analogously,   the backward  domain  of influence  $D_-(F)$  of the closed set  $F\subset D_0\times\R$  is the closure of the union of
all backward  time-like piece-wise smooth rays in $D_0\times\R$ starting at $F$.
The domain of dependence of $F$ is the intersection  $D_+(F)\cap D_-(F)$.

Let  $\Gamma\subset\partial D_0$  and  let  $Lu=0$  in  $D_0\times \R$.
A consequence  of  the unique  continuation  property  is that $u\big|_{\Gamma\times(t_1,t_2)}
=\frac{\partial u}{\partial \nu}\Big|_{\Gamma\times(t_1,t_2)}=0$  implies $u=0$  in the domain  of dependence  of
$\Gamma\times[t_1,t_2]$.  Here  $\frac{\partial}{\partial \nu}$  is the normal derivative to $\Gamma$.
This fact follows from  [KKL1]   in the case  of time-independent  coefficients.  The proof  in the time-dependent   case  is similar.

The following fact  follows from the BLR condition:

Consider  $\Gamma\times[t_1,t_2],\Gamma\subset \partial D_0$.  Suppose  $[t_1,t_2]$  is arbitrary large.  Then the domain 
 of dependence of $\overline\Gamma\times[t_1,t_2]$  contains  $\overline D_0\times[t_1+\delta,t_2-\delta]$  for some  $\delta>0$
 dependent  of the metric  and the domain.
 
 In this paper  we will not attempt to estimate  $\delta>0$  since 
 $[t_0+\delta,t_2-\delta]$  is also  arbitrary large if  $[t_1,t_2]$  is arbitrary large.
 \qed

Now we shall  state  the main result of this paper.

Consider  an admissible domain  $D$  in $\R^{n+1}$  and an initial-boundary value  problem in $D$.  

Using  the map of the form (\ref{eq:1.7})  defining the admissibility
  of the domain  $D$  we get   a cylindrical domain  $D_0\times \R$   with  $S_0=\Gamma_0\times\R$   
  (cf.  Fig. 1.2)  and the initial-boundary  value  problem
\begin{align}																\label{eq:1.34}
&Lu = 0 \ \ \ \mbox{in}\ \ D_0\times\R,
\\
&u = 0 \ \ \ \mbox{when}\ \ x_0 \ll 0,    										\label{eq:1.35}
\\
&u\big|_{\partial D_0\times \R} =f,																\label{eq:1.36}
\end{align}
where $L$ has the form  (\ref{eq:1.1}) and $f$  has a compact support in $\overline \Gamma_0\times\R$.  
Consider 
another  admissible  domain  $\hat D$. Making  again 
the change of variables (\ref{eq:1.7})  we get a cylindrical  domain  $\hat D_0\times\R$   and
another initial-boundary value problem
\begin{align}																\label{eq:1.37}
&L'u' = 0 \ \ \ \mbox{in}\ \ \hat D_0\times\R,
\\
&u' = 0 \ \ \ \mbox{when}\ \ y_0 \ll 0,											\label{eq:1.38}
\\
&u'\big|_{\partial\hat D_0\times\R} =f',																\label{eq:1.39}
\end{align}
where $L'u'$  has  the form  (\ref{eq:1.30}),  $f'$  has a compact support in $\overline\Gamma_0\times\R$.
Therefore  the inverse  problems for the  admissible  domains  are reduced  to the inverse  problems  in cylindrical  domains.

We shall prove  the following  theorem: 
\begin{theorem}      													      \label{theo:1.2}
Consider  two initial-boundary value problems (\ref{eq:1.34}), (\ref{eq:1.35}), (\ref{eq:1.36})  and 
(\ref{eq:1.37}), (\ref{eq:1.38}), (\ref{eq:1.39}) 
in domains $D_0\times\R$  and $\hat D_0\times\R$,  respectively.  
Suppose  $A_j(x), A_j'(y),0\leq j\leq n,$  are real-valued.  Assume  that  $\Gamma_0\subset\partial D_0\cap \partial\hat D_0$  is
nonempty  and open.  Let  $\Lambda$  and $\Lambda'$  be the  corresponding  DN operators for  $L$  and  $L'$.
Assume that $\Lambda f\big|_{\Gamma_0\times\R}=\Lambda' f\big|_{\Gamma_0\times\R}$  for all smooth  $f$  with compact support
in $\overline \Gamma_0\times \R$.  Suppose the conditions  (\ref{eq:1.2}),  (\ref{eq:1.6})  hold for $L$  and  $L'$.  Assume
that the coefficients of $L$  and $L'$  are analytic   in $x_0$  and $y_0$,  respectively.
Suppose also  that BLR  condition holds  for (\ref{eq:1.34}), (\ref{eq:1.35}), (\ref{eq:1.36}) on $[t_0,T_{t_0}]$  for each  $t_0\in \R$.
Then there exists a proper
map  $y=y(x)$  of $\overline D_0\times\R$  onto $\overline{\hat D_0}\times\R,\ y=x$  on  $\Gamma_0\times\R$, 
 and there 
exists a gauge transformation  with  the gauge $c'(y)\in G_0(\overline{\hat D_0}\times\R),\ c'(y)=1$  
on $\overline \Gamma_0\times\R$ such  that 
$L'= c'\circ y^* L$.
Here  $y^*\circ L$  is the operator  with $[\hat g^{jk}(y)]_{j,k=0}^n$  and  $\hat A_k(y),\ 0\leq k\leq n$ 
 as  in (\ref{eq:1.17}), (\ref{eq:1.18}),
$c'\circ y^*\circ L$  is the operator  with potentials  
$A_j'(y),\   0  \leq j\leq n$,  gauge  equivalent to $\hat A_k(y),\ 0\leq k\leq n$.	
\end{theorem}

We end the introduction  with the outline of the previous work  and a short description  of the content of the paper.

The first result on inverse  hyperbolic problems with the data on the part  of the boundary  was obtained by Isakov in [I1].  The powerful Boundary
Control (BC) method  was discovered by Belishev [B1] and was further developed by Belishev [B2], [B3], [B4],  
 Belishev and Kurylev  [BK],
Kurylev and Lassas [KL1],  [KL2] and others (see [KKL1],  [KKL2]).
In [E1], [E2]  the author proposed  a new approach  to hyperbolic inverse problems  that uses substantially the idea of BC method.  This
approach was extended in [E3]  to some class of time-independent metrics with time-dependent  vector
potentials  and in [E4] to the case of hyperbolic equations of
general form with time-independent coefficients  without vector potentials. 
The generalization to the case of Yang-Mills potentials was considered in [E7].
The inverse problems for the D'Alambert equation with the time-dependent scalar potentials were considered 
earlier by Stefanov [S]  and Ramm
and Sjostrand [RS] (see also Isakov [I2]).  The case of the D'Alambert  equation with time-dependent vector potentials was studied by
Salazar  [S1], [S2]:

As it was  mentioned in [BLR]  the study  of hyperbolic  equations with  time-dependent coefficients  is very important  
because  the linearization  of basic  nonlinear  hyperbolic  equations of mathematical physics leads  to time-dependent  linear  
hyperbolic  equations.

The main result  of  the present paper is the determination  of the  time-dependent  Lorentzian  metric by the boundary 
measurements given  on the part of the boundary.  We  consider  the second order hyperbolic equations of general form (\ref{eq:1.1})
with time-dependent  coefficients  and vector potentials.  The method  is the extension  of the approach  in [E1],  [E2],  [E4]   to the 
case of time-dependent metrics.  We adapt  some lemmas  of  [E1], [E2], [E3], [E4]  to the 
time-dependent situation and simultaneously give sharper and simpler proofs.

The main step
in the proof
 is the local step of solving the inverse problem in a small neighborhood  near the boundary.
This is done in \S\S 2-6. 

In \S 2 we make a change of variables  in a neighborhood $U_0\subset \R^{n+1}$  of a point  $x^{(0)}\in S_0$.
We called the new coordinates  the Goursat coordinates since  they are similar  to coordinates arising  in a solution 
of the Goursat  problem  in the case   of hyperbolic equations with  one space  variable.  The Goursat  coordinates  allow
to simplify  the equation  (\ref{eq:1.1}).  Another  advantage  of the Goursat coordinates  is that  
the characteristic  surface  is a plane in these coordinates.  Also it is proven in  \S 2 that the original DN  operator  determines
the new  DN  operator  corresponding to the equation  in Goursat  coordinates.

In  \S  3  we  derive  the  Green  formula  in Goursat  coordinates   and prove  the crucial  density  lemma   (Lemma  \ref{lma:3.1}).

In \S 4  we establish  the main formula  that states that some integrals of solutions  of the initial-boundary  value problems are determined
by the DN  operator  (see Theorem  \ref{theo:4.5},  formula (\ref{eq:4.29})).   To establish  this formula  one needs to compare 
Sobolev norms  on the characteristic surfaces  corresponding   to different  operators having  the same DN  data.  
This is done  by  using  the BLR  condition.  Note that in the case  of time-independent  coefficients   there is an additional  energy 
identity  that allows to avoid 
the use  of the BLR  condition  (see Remark  4.1).  
  Also note  that  proofs  in \S 4   require that the hyperbolic operators  are formally self-adjoint.  Thus the vector  potentials   $A=(A_0,A_1,...,A_n)$  are  required    to be  real-valued.
  
  In \S5  we construct  geometric  optics  type  solutions  depending  on a large parameter $k$.    
Substituting  the geometric  optics solutions  in the main formula ,    we prove  
in \S5 and  \S6 the local version  
(Theorem \ref{theo:6.2})   of the main theorem  (Theorem  \ref{theo:1.2}).  In the  last  \S  7  we consider  the global  case.
At first  we study  the case  of a finite  time  interval  (Theorem  \ref{theo:7.7}))
and then finally  prove  Theorem  \ref{theo:1.2}.

\section{The Goursat coordinates}
\label{section 2}
\init
We shall prove  first the Theorem \ref{theo:1.2}
in the small neighborhood of the boundary $\partial D$.

Let  $x^{(0)}\in S_0$  and let 
$U_0\subset\R^{n+1}$  be a small neighborhood  of  $x^{(0)}$.  

Suppose that we already did the change of variables (\ref{eq:1.7})  to make $\partial D$  and  $S_0$  cylindrical,  i.e. 
$\partial D=\partial D_0\times \R$  and $S_0=\Gamma_0\times \R$.  
   We assume that we have chosen  the coordinates  $(x_0,x',x_n),\ x'=(x_1,...,x_{n-1})$  in $U_0$   such  that  $x_n=0$
    is the equation  of $U_0\cap\partial D$  and $U_0\cap D$  is contained in
the half-space $x_n>0$.   Let  $(x_0^{(0)},x_1^{(0)},...,x_{n-1}^{(0)},0)$ be  the coordinates   of the point $x^{(0)}$.  Let
$T_1<x_0^{(0)}<T_2, \ T_2-T_1$  is small.

Consider the initial-boundary value  problem in $U_0\cap D$:
\begin{align}                                                              									\label{eq:2.1}
&Lu=0, \ \ \ x_n>0, \ \ T_1<x_0<T_2,
\\
&u\big|_{x_0=T_1}=0,\ \ \ \frac{\partial u}{\partial x_0}\Big|_{x_0=T_1}=0,				\label{eq:2.2}
\\
&u\big|_{x_n=0}=g(x_0,x').															\label{eq:2.3}
\end{align}
We assume that  $L$  has the form (\ref{eq:1.1}).  For the simplicity,  we shall not change the notations when choosing the local coordinates  such that
 the equation  of $U_0\cap S_0$  is  $x_n=0$.  Assume that $\mbox{supp}\,g \subset U_0\cap(\Gamma_0\times[T_1,T_2]),\ 
g=0$   for  $x_0<T_1$.  Note   
that $\mbox{supp}\,u(x_0,x',x_n)\cap[T_1,T_2]\subset U_0\cap[T_1,T_2]$  for $x_n>0$  if  $T_2-T_1$  is small.

We introduce new coordinates  to simplify  the operator $L$ (cf. [E4],  pages 327-329)  that we called the Goursat coordinates.

Denote by $\psi^\pm(x),x=(x_0,x',x_n)$  the solutions of the eikonal  equations
\begin{equation}              																\label{eq:2.4}
 \sum_{j,k=0}^ng^{jk}(x_0,x',x_n)\psi_{x_j}^\pm(x)\psi_{x_k}^\pm(x)
=0,\ \ \ x_n>0,
\end{equation}
with 
initial conditions
\begin{equation}																		\label{eq:2.5}
\psi^+\big|_{x_n=0}=x_0-T_1,\ \  \psi^-\big|_{x_n=0}=T_2-x_0.
\end{equation}
It is well known  that the solution  $\psi^\pm(x)$  of  (\ref{eq:2.4}),  (\ref{eq:2.5})   exists  for  $0\leq x_n\leq \e$,
where $\e>0$  is  small  (see,  for example,  [E5],  \S 64).  

Since  (\ref{eq:2.4}) is a quadratic  equation  in $\psi_{x_n}^\pm$  one has to specify  the sign  of the square root.  
We have 
$$
g^{nn}(\psi_{x_n}^\pm)^2+2\sum_{j=0}^{n-1}g^{nj}\psi_{x_j}^\pm\psi_{x_n}^\pm+\sum_{j,k=0}^{n-1} g^{jk}\psi_{x_j}^\pm\psi_{x_k}^\pm
=0.
$$
We will need
below that
 $\psi_{x_n}^++\psi_{x_n}^-<0$  for $x_n>0$.  So we choose  the plus sign  
of the square root:
\begin{equation} 																		\label{eq:2.6}
\psi_{x_n}^\pm=
\frac{-\sum_{j=0}^{n-1}g^{nj}\psi_{x_j}^\pm
+\sqrt{\big(\sum_{j=0}^{n-1}g^{nj}\psi_{x_j}^\pm\big)^2-g^{nn}\big(\sum_{j,k=0}^{n-1}g^{jk}\psi_{x_j}^\pm\psi_{x_k}^\pm\big)}}
{g^{nn}(x)}
\end{equation}
Note that $g^{nn}(x)<0,\ \psi_{x_0}^\pm\big|_{x_n=0}=\pm 1$.   Therefore $\psi_{x_n}^\pm\big|_{x_n=0}=
\frac{\mp g^{n0}+\sqrt{(g^{n0})^2-g^{nn}g^{00}}}{g^{nn}}$.
The solutions $\psi^\pm(x)$  exists for  $0<x_n<\delta,\ \delta$  is small. 
For  given  $T_1, T_2$ we assume that  $\delta$  is  such that surfaces 
$\psi^+=0$  and  $\psi^-=0$  intersect  when $x_n<\delta$  and  are inside $U_0$  when $x_n<\delta$.

Let  $\varphi_j(x_0,x',x_n),1\leq j\leq n-1,$  be solutions of the linear  equation
\begin{equation}																		\label{eq:2.7}
\sum_{p,k=0}^ng^{pk}(x_0,x',x_n)\psi_{x_p}^-\varphi_{jx_k}=0,\ \ x_n>0,
\end{equation}
with initial condition
\begin{equation}  																		\label{eq:2.8}
\varphi_j(x_0,x',0)=x_j,\ \ \ 1\leq j\leq n-1.
\end{equation}
Make  the following  change of variables  in $U_0\cap[T_1,T_2]$:
\begin{align}             																	\label{eq:2.9}
&s=\psi^+(x_0,x',x_n),
\\
\nonumber
&\tau=\psi^-(x_0,x',x_n),
\\
\nonumber
&y_j=\varphi_j(x_0,x',x_n),\ 1\leq j\leq n-1.
\end{align}
Equation  (\ref{eq:1.1})  has  the following  form in $(s,\tau,y')$  coordinates where $y'=(y_1,...,y_{n-1})$
\begin{align}																			\label{eq:2.10}
\hat L\hat u \overset{def}{=}
&
\frac{2}{\sqrt{|\hat g|}}
\Big(\frac{\partial}{\partial s}+i\hat A_+(s,\tau,y')\Big) \sqrt{|\hat g|}\,\hat g^{+,-}(s,\tau,y')\Big(\frac{\partial}{\partial\tau}+i\hat A_-\Big)\hat u
\\
\nonumber
&
+\frac{2}{\sqrt{|\hat g|}}
\Big(\frac{\partial}{\partial \tau}+i\hat A_-(s,\tau,y')\Big) \sqrt{|\hat g|}\,\hat g^{+,-}(s,\tau,y')\Big(\frac{\partial}{\partial s}+i\hat A_+\Big)\hat u
\\
\nonumber
&
-\sum_{j=1}^{n-1}\frac{2}{\sqrt{|\hat g|}}
\Big(\frac{\partial}{\partial y_j}-i\hat A_j(s,\tau,y')\Big) \sqrt{|\hat g|}\,\hat g^{+,j}(s,\tau,y')\Big(\frac{\partial}{\partial s}+i\hat A_+\Big)\hat u
\\
\nonumber
&
-\sum_{j=1}^{n-1}\frac{2}{\sqrt{|\hat g|}}
\Big(\frac{\partial}{\partial s}+i\hat A_+(s,\tau,y')\Big) \sqrt{|\hat g|}\,\hat g^{+,j}(s,\tau,y')\Big(\frac{\partial}{\partial y_j}-i\hat A_j\Big)\hat u
\\
\nonumber
&
-\sum_{j,k=1}^{n-1}
\frac{1}{\sqrt{|\hat g|}}
\Big(\frac{\partial}{\partial y_j}-i\hat A_j(s,\tau,y')\Big) \sqrt{|\hat g|}\,\hat g^{jk}(s,\tau,y')\Big(\frac{\partial}{\partial y_k}-i\hat A_k\Big)\hat u,
\end{align} 
where
\begin{equation}   																			\label{eq:2.11}
\hat g=-(2\hat g^{+,-})^{-2}\big(\det[\hat g^{jk}]_{j,k=1}^{n-1}\big)^{-1}.
\end{equation}
Note that terms containing $\frac{\partial^2}{\partial s^2},\frac{\partial^2}{\partial \tau^2},\frac{\partial^2}{\partial y_j\partial\tau}$
vanished because of (\ref{eq:2.4}), (\ref{eq:2.7}),  and
\begin{align}																				\label{eq:2.12}
& 2\hat g^{+,-}=-\sum_{j,k=0}^ng^{jk}\psi_{x_J}^+\psi_{x_k}^-,
\\
\nonumber
&2\hat g^{+,j}=\sum_{p,r=0}^n g^{pr}\psi_{x_p}^+\varphi_{jx_r},\ \ 1\leq j\leq n-1,
\\
\nonumber
&\hat g^{jk}=\sum_{p,r=0}^n g^{pr}\varphi_{j x_p}\varphi_{kx_r},\ \ 1\leq j,k\leq n-1,
\end{align}
It follows  from  (\ref{eq:2.6})    for  $x_n=0$  that  $g^{+,-}>0$.

In  (\ref{eq:2.10})  $\hat u(s,\tau,y')=u(x_0,x',x_{n}),$
\begin{equation}                           															\label{eq:2.13}
A_k(x)=\sum_{j=1}^{n-1}\hat A_j(s,\tau,y')\varphi_{jx_k}-\hat A_+\psi_{x_k}^+ -\hat A_-\psi_{x_k}^-,\ \ 0\leq k\leq n.
\end{equation}
Now we shall introduce a new system of coordinates (cf. [E4])
\begin{align}                 																	\label{eq:2.14}
&y_0=\frac{s-\tau+T_2+T_1}{2}
=\frac{\psi^+-\psi^-+T_2+T_1}{2},
\\
\nonumber
&y_j=\varphi_j(x),\ \ \ 1\leq j\leq  n-1,
\\
\nonumber
&y_n=\frac{T_2-T_1-s-\tau}{2}=\frac{T_2-T_1-\psi^+(x)-\psi^-(x)}{2},
\end{align}
where $\psi^+,\psi^-,\varphi_j,\ 1\leq j\leq n-1$,  are  the same as in (\ref{eq:2.4}),  (\ref{eq:2.7}).

Note that
\begin{align}                 																	\label{eq:2.15}
&y_0\big|_{x_n=0}=\frac{x_0-T_1-(T_2-x_0)+T_2+T_1}{2}=x_0,
\\
\nonumber
&y_j\big|_{x_n=0}=x_j,\ \ \ 1\leq j\leq  n-1,
\\
\nonumber
&y_n\big|_{x_n=0}=\frac{T_2-T_1-s-\tau}{2}=\frac{T_2-T_1-\psi^+(x)-\psi^-(x)}{2}=0,
\end{align}
Therefore $y=\varphi(x)=(\varphi_0(x_1),\varphi_1(x),...,\varphi_n(x))$  is the identity on $x_n=0$: 
\begin{equation}																		\label{eq:2.16}
\varphi(x)=I\ \ \ \mbox{when}\ \ x_n=0.
\end{equation}
Here 
$$
\varphi_0=\frac{\psi^+(x)-\psi^-(x)+T_2+T_1}{2},\ \ \varphi_n=\frac{T_2-T_1-\psi^+-\psi^-}{2}.
$$
Note that $y_n=\varphi_n(x)>0$  when  $x_n>0$  since $\psi_{x_n}^++\psi_{x_n}^-<0$  (cf. (\ref{eq:2.6})),
\begin{equation} 				 														\label{eq:2.17}
u_s=\frac{1}{2}(u_{y_0}-u_{y_n}),\ \ \ u_\tau=-\frac{1}{2}(u_{y_0}+u_{y_n}).
\end{equation}
Thus one can easily  rewrite (\ref{eq:2.10})  in $(y_0,y',y_n)$  coordinates .

We shall further simplify  (\ref{eq:2.10})  by making a gauge transformation
\begin{equation}   																		\label{eq:2.18}
u'=e^{-id(s,\tau,y')}\hat u.
\end{equation}
Then $u'$  satisfies the equation
\begin{equation}           																	\label{eq:2.19}
L'u'=0,
\end{equation}
where $L'$  is the same as $\hat L$  with $\hat A_j, \hat A_+, \hat A_-$  replaced  by  $A_j', A_+', A_-',\
1\leq j\leq n-1$,  where
\begin{align} 																		      \label{eq:2.20}
&A_j'=\hat A_j-\frac{\partial d}{\partial y_j},\ \ 1\leq j\leq n-1,
\\
\nonumber
&A_+'=\hat A_+-\frac{\partial d}{\partial s},\ \ \ A_-'=\hat A_--\frac{\partial d}{\partial\tau}.
\end{align}
We choose $d(s,\tau,y')$  such that
\begin{align}																			\label{eq:2.21}
&A_+'=-\frac{\partial d}{\partial s}+\hat A_+=0 \ \ \mbox{for}\ \ y_n>0,
\\
\nonumber
&d\big|_{y_n=0}=0.
\end{align}
Let
\begin{equation}																		\label{eq:2.22}
g_1=\big|\det[\hat g^{jk}]_{j,k=1}^{n-1}\big|^{-1},\ \ A=\ln(g_1)^{\frac{1}{4}}.
\end{equation}
Note that 
\begin{align}																			\label{eq:2.23}
&\frac{\partial A}{\partial y_j}=\frac{g_{1y_j}}{4g_1}
=\frac{1}{2}\frac{1}{\sqrt{g_1}}\frac{\partial}{\partial y_1}\sqrt{g_1},\ \ \ 1\leq j\leq n-1,
\\
\nonumber
&\frac{\partial A}{\partial s}=\frac{g_{1s}}{4g_1}=\frac{1}{2}\frac{1}{\sqrt{g_1}}\frac{\partial}{\partial s}\sqrt{g_1},
\\
\nonumber
&\frac{\partial A}{\partial \tau}=\frac{g_{1\tau}}{4g_1}=\frac{1}{2}\frac{1}{\sqrt{g_1}}\frac{\partial}{\partial \tau}\sqrt{g_1}.
\end{align}
Since  
$\sqrt{|\hat g|}=\frac{\sqrt{g_1}}{2\hat g^{+,-}}$  (cf (\ref{eq:2.11}))  we can  rewrite  $L'u'=0$  in the form (cf.  [E1]):

\begin{align}																			\label{eq:2.24}
L'u'=&2\hat g^{+,-}\Big(\frac{\partial}{\partial s}+\frac{\partial A}{\partial s}
\Big)  \Big(\frac{\partial}{\partial \tau}+iA_-'+\frac{\partial A}{\partial \tau}\Big)
\\
\nonumber
&+2\hat g^{+,-}\Big(\frac{\partial}{\partial\tau} +iA_-'+\frac{\partial A}{\partial \tau}\Big)\,\Big(\frac{\partial}{\partial s}+
\frac{\partial A}{\partial s}\Big)u'
\\
\nonumber
&-2\hat g^{+,-}\sum_{k=1}^{n-1}\Big(\frac{\partial}{\partial s}+\frac{\partial A}{\partial s}\Big)
\frac{\hat g^{+,k}}{\hat g^{+,-}}\Big(\frac{\partial}{\partial y_k}-iA_k'+\frac{\partial A}{\partial y_k}\Big)u'
\\
\nonumber
&-2\hat g^{+,-}\sum_{k=1}^{n-1}\Big(\frac{\partial}{\partial y_k}-iA_k'+\frac{\partial A}{\partial y_k}\Big)
\frac{\hat g^{+,k}}{\hat g^{+,-}}
\Big(\frac{\partial}{\partial s}+\frac{\partial A}{\partial s}\Big)
u'
\\
\nonumber
&-\sum_{j,k=1}^{n-1}\hat g^{+,-}\Big(\frac{\partial}{\partial y_j}-iA_j'+\frac{\partial A}{\partial y_j}\Big)
\frac{\hat g^{jk}}{\hat g^{+,-}}\Big(\frac{\partial}{\partial y_k}-iA_k'+\frac{\partial A}{\partial y_k}\Big)u'
\\
\nonumber
&+\hat g^{+,-}V_1u'=0,
\end{align}
where
\begin{align}																				\label{eq:2.25}
V_1=
-&\sum_{j,k=1}^{n-1}\Bigg(\frac{\hat g^{jk}}{\hat g^{+,-}}\frac{\partial A}{\partial y_j}\frac{\partial A}{\partial y_k}
+\frac{\partial}{\partial y_k}\Big(\frac{\hat g^{jk}}{\hat g^{+,-}}\frac{\partial A}{\partial y_j}\Big)\Bigg)
\\
\nonumber
+&4\frac{\partial^2 A}{\partial s \partial\tau} +4\frac{\partial A}{\partial s}\frac{\partial A}{\partial \tau}-
4\sum_{j=1}^{n-1}\frac{\hat g^{+,j}}{\hat g^{+,-}}\frac{\partial A}{\partial s}\frac{\partial A}{\partial y_j}
\\
\nonumber
-&2\sum_{j=1}^{n-1}\Bigg(\frac{\partial}{\partial s}\Big(\frac{\hat g^{+,j}}{\hat g^{+,-}}\frac{\partial A}{\partial y_j}\Big) 
+\frac{\partial }{\partial y_j}\Big(\frac{\hat g^{+,j}}{\hat g^{+,-}}\frac{\partial A}{\partial s}\Big)\Bigg).
\end{align}
Make the change of unknown function
\begin{equation}																			\label{eq:2.26}
u_1= g_1^{\frac{1}{4}} u',
\end{equation}
where  $ g_1=|\det[\hat g^{jk}]_{j,k=1}^{n-1}|^{-1}$  (cf.  (\ref{eq:2.11})).
Then dividing $L'u'=0$  by  $\hat g^{+,-}$  we get (cf.  [E4])
$$
L_1u_1=0,
$$
where   $L_1u_1=0$  has the form  (cf. (\ref{eq:2.24}))
\begin{align}           																			\label{eq:2.27}
L_1u_1=&2\frac{\partial}{\partial s}\Big(\frac{\partial}{\partial\tau}+iA_-'\Big)u_1+2\Big(\frac{\partial}{\partial\tau}+iA_-'\Big)
\frac{\partial}{\partial s}u_1
\\
\nonumber
&-2\sum_{j=1}^{n-1}\frac{\partial}{\partial s}\Big(g_0^{+,j}\Big(\frac{\partial}{\partial y_j}-iA_j'\Big)u_1\Big)
\\
\nonumber
&-\sum_{j=1}^{n-1}2\Big(\frac{\partial}{\partial y_j}-iA_j'\Big)g_0^{+,j}\frac{\partial u_1}{\partial s}
\\
\nonumber
&-\sum_{j,k=1}^{n-1}\Big(\frac{\partial}{\partial y_j}-iA_j'\Big)g_0^{jk}\Big(\frac{\partial}{\partial y_k}-iA_k'\Big)u_1+V_1u_1=0,
\end{align}
where
$$
g_0^{jk}=\frac{\hat g^{jk}}{\hat g^{+,-}},\ \ \  g_0^{+,j}=\frac{\hat g^{+,j}}{\hat g^{+,-}},
$$
and $V_1$  is the same as in (\ref{eq:2.25}).
Using that $\frac{\partial}{\partial \tau}+iA_-'=\frac{1}{2}\big(-\frac{\partial}{\partial y_0}-\frac{\partial}{\partial y_n}\big)+iA_-'
= -\frac{1}{2}\big[\big(\frac{\partial}{\partial y_0}-iA_-'\big)+\big(\frac{\partial}{\partial y_n}-iA_-'\big)\big]$
and $\frac{\partial}{\partial s}=\frac{1}{2}\big(\frac{\partial}{\partial y_0}-\frac{\partial}{\partial y_n}\big)
= \frac{1}{2}\big[\big(\frac{\partial}{\partial y_0}-iA_-'\big)-\big(\frac{\partial}{\partial y_n}-iA_-'\big)\big]$
we can rewrite $L_1u_1$  in $(y_0,y',y_n)$  coordinates:
\begin{align}																			\label{eq:2.28}
&L_1u_1=-\Big(\frac{\partial}{\partial y_0}-iA_-'\Big)^2u_1+ \Big(\frac{\partial}{\partial y_n}-iA_-'\Big)^2u_1
\\
\nonumber
&-\sum_{j=1}^{n-1}\Big(\frac{\partial}{\partial y_0}-iA_-'\Big) g_0^{+,j}\Big(\frac{\partial}{\partial y_j}-iA_j'\Big)u_1
-\sum_{j=1}^{n-1}\Big(\frac{\partial}{\partial y_j}-iA_j'\Big) g_0^{+,j}\Big(\frac{\partial}{\partial y_0}-iA_-'\Big)u_1
\\
\nonumber
&+\sum_{j=1}^{n-1}\Big(\frac{\partial}{\partial y_n}-iA_-'\Big)g_0^{+,j}\Big(\frac{\partial}{\partial y_j}-iA_j'\Big)u_1
+\sum_{j=1}^{n-1}\Big(\frac{\partial}{\partial y_j}-iA_j'\Big)g_0^{+,j}\Big(\frac{\partial}{\partial y_n}-iA_-'\Big)u_1
\\
\nonumber
&-\sum_{j,k=1}^n\Big(\frac{\partial}{\partial y_j}-iA_j'\Big)g_0^{jk}\Big(\frac{\partial}{\partial y_k}-iA_k'\Big)u_1
+V_1u_1=0.
\end{align}
Note that  we transformed  the equation  $Lu=0$  to the equation  $L_1u_1=0$  in two steps.  First,  we transformed  $Lu=0$  to
$L'u'=0$  by making  the change  of variables $y=\varphi(x)$  of  the form  (\ref{eq:2.15})  and gauge 
transformation  with the gauge  $e^{-id(s,\tau,y')}$   belonging to the
 group  $G_0$  (cf.  (\ref{eq:1.22}).  Then we transform  $L'u'=0$  to  $L_1u_1=0$   by using  the change of variables  (\ref{eq:2.26}),  i.e.     
 by   using  a gauge  $e^A$,  where  $A=\ln g_1^{\frac{1}{4}}$  and  then dividing  $L'u'=0$    by  $\hat  g^{+,-}$.
 \qed

The DN operator  for $L$  has the form
$$
\Lambda g=-\sum_{j=0}^n g^{jn}(x)\Big(\frac{\partial u}{\partial x_j}-iA_j(x)u\Big)(-g^{nn}(x))^{-\frac{1}{2}}\Big|_{x_n=0}
$$
since the outward normal to $x_n=0$  is $(0,0....,-1)$.

Rewrite $L'u'=0$  in $(y_0,y',y_n)$  coordinates using  (\ref{eq:2.17}). 

Denote  $-\hat g^{00}=\hat g^{nn}=\hat g^{+,-},\ \hat g^{nj}=\hat g^{jn}=\hat g^{+,j},\ 1\leq j\leq n-1$.  Note that $\hat g^{+,-}>0.$

The DN  operator  for  $L'u'=0$  has the following  form  in $(y_0,y',y_n)$  coordinates:
\begin{equation}                         															\label{eq:2.29}
\Lambda'g=(\hat g^{+,-})^{\frac{1}{2}}\Big[\Big(\frac{\partial u'}{\partial y_n}-iA_-'u'\Big)
+\sum_{j=1}^{n-1}\frac{\hat g^{nj}}{\hat g^{+,-}}\Big(\frac{\partial u'}{\partial y_j}-iA_j'u'\Big)\Big]\Big|_{y_n=0},
\end{equation}
where (cf.  (\ref{eq:2.18}))
$$
u'(y)=e^{-id(y)}u(\varphi^{-1}(y)),
$$
$y=y(x)$  is  the same  as in  (\ref{eq:2.14}).

Since $L'$ is obtained from (\ref{eq:1.1})  by the change of variables
(\ref{eq:2.14})and the gauge transformation  (\ref{eq:2.18}) 
and 
since
(\ref{eq:2.15}),  (\ref{eq:2.21}) hold,  we have 
$\Lambda g=\Lambda' g $ on $\{y_n=0\}\cap U_0$  for  all $g$  with $\mbox{supp}\,g$  in  $(\Gamma_0\times[T_1,T_2])\cap U_0$.
Using  the expression  of $L_1u_1=0$  in $(y_0,y',y_n)$  coordinates  (see (\ref{eq:2.28}))  we get  that DN operator $\Lambda_1g$ 
 has the form
\begin{equation}																		\label{eq:2.30}
\Lambda_1g=\Bigg(\frac{\partial u_1}{\partial y_n}-iA_-'u_1
+\sum_{j=0}^{n-1}g_0^{+,j}\Big(\frac{\partial u_1}{\partial y_j}-iA_j'u_1\Big)\Bigg)\Bigg|_{y_n=0},
\end{equation}
where $g_0^{+,j}=\frac{\hat g^{nj}}{\hat g^{+,-}}$.

We shall show  that the DN  operators $\Lambda'$  determines  the DN operator $\Lambda_1$   in $U_0\cap \Gamma_0$.

The following lemma  is well known,  especially in the elliptic case  (cf.  [LU], [E5], \S57).  
For the hyperbolic case see [E1], Remark 2.2.
\begin{lemma}																			\label{lma:2.1}
The DN operator $\Lambda'$ determines 
\begin{equation}																		\label{eq:2.31}
\hat g^{+,-}\big|_{y_n=0},\ \ \frac{\hat g^{nj}}{\hat g^{+,-}}\Big|_{y_n=0},\ \ \frac{\hat g^{jk}}{\hat g^{+,-}}\Big|_{y_n=0},\ \ 
1\leq j\leq n-1,\ 1\leq k\leq n-1,
\end{equation}
and the derivatives of  (\ref{eq:2.31}) in $y_n$  at $y_n=0$.
\end{lemma}
{\bf Proof:}.  The principal 
symbol of operator $L'$  has the form  $\hat g^{+,-}p(y,\eta)$,  where  (cf.  (\ref{eq:2.10}) in $y$-coordinates) 
\begin{equation}																		\label{eq:2.32}
p(y,\eta)=\eta_0^2-\eta_n^2+2\sum_{j=1}^{n-1}g_0^{+,j}(\eta_0-\eta_n)\eta_j+\sum_{j,k=1}^{n-1} g_0^{jk}\eta_j\eta_k,
\end{equation}
where
\begin{equation}																		\label{eq:2.33}
g_0^{+.j}=\frac{\hat g^{+,j}}{\hat g^{+,-}},\ \  g_0^{jk}=\frac{\hat g^{jk}}{\hat g^{+,-}}.
\end{equation}
Since   (\ref{eq:1.6})  holds   the quadratic  form  $\sum_{j,k=1}^{n-1}g_0^{jk}\eta_j\eta_k$  is negative  definite.  Therefore  for  $\e>0$
 in the region $\Sigma=\{\eta_0^2
+\big(\sum_{j=1}^{n-1}g_0^{+,j}\eta_j\big)^2
-\e\sum_{j=1}^{n-1}\eta_j^2<0\}$  of the cotangent  space  
$T^*=U_0\times(\R^{n+1}\setminus\{0\})$   the  operator  $p(y,\eta)$  is elliptic.  We shall call $\Sigma$  the elliptic  region.

There is a parametrix  of the Dirichlet  problem  in the elliptic  region  and DN  operator  microlocally  in $\Sigma$ is 
a pseudodifferential  operator  on $y_n=0$.   We shall  find the principal symbol of  this operator  in $\Sigma$.  Let  $\lambda_\pm$
be the roots in $\eta_n$  of  $p(y,\eta_0,\eta',\eta_n)=0$:
\begin{multline}																			\label{eq:2.34}
\lambda_\pm=-\sum_{j=1}^{n-1}g_0^{+,j}\eta_j\pm\sqrt{\Big(\sum_{j=1}^{n-1}g_0^{+,j}\eta_j\Big)^2+
\Big(\eta_0^2 +2\sum_{j=1}^{n-1} g_0^{+,j}\eta_j\eta_0+\sum_{j,k=1}^{n-1} g_0^{jk}\eta_j\eta_k\Big)} 
\\
\stackrel{def}{=}-\sum_{j=1}^{n-1} g_0^{+,j}\eta_j\pm\sqrt Q,
\end{multline}
where $\Im \lambda_+>0$ in $\Sigma$.  Therefore the symbol  of DN  in $\Sigma$ is (cf. [E5],  \S57):
\begin{equation}     																		\label{eq:2.35}
\big(\hat g^{+,-}\big)^{\frac{1}{2}}\sqrt Q.
\end{equation}
Knowing   $\Lambda'$  we know  the symbol  (\ref{eq:2.35}) for all $\eta_0,\eta'$.  In particular,  we can  find (\ref{eq:2.31}).  
Computing  the next term  of the parametrix  (cf.  [E5],  \S57)  we 
can find the normal  derivatives of (\ref{eq:2.31}).
\qed

We have 
\begin{equation}																		\label{eq:2.36}
\Lambda_1f'=\Big(\frac{\partial u_1}{\partial y_n}-iA_-'u_1
+\sum_{j=1}^{n-1} g_0^{+,j}\Big(\frac{\partial u_1}{\partial y_j}-iA_j'u_1\Big)\Big)\Big|_{y_n=0},
\end{equation}
where
$u_1\big|_{y_n=0}=f',\ u_1= g_1^{\frac{1}{4}}u',\ f'=g_1^{\frac{1}{4}}\big|_{y_n=0}f$.  Note that  $\frac{\partial}{\partial y_k}u_1
= g_1^{\frac{1}{4}}\frac{\partial}{\partial y_k}u'+\big(\frac{\partial}{\partial y_k} g_1^{\frac{1}{4}}\big)u'$.   Therefore
\begin{equation}																		\label{eq:2.37}
\Lambda_1f'= g_1^{\frac{1}{4}} (g^{+,-}\big)^{-\frac{1}{2}}\Lambda'  f +
\Big(\frac{\partial  g_1^{\frac{1}{4}}}{\partial y_n}
+\sum_{j=1}^{n-1}g_0^{+,j}\frac{\partial  g_1^{\frac{1}{4}}}{\partial y_j}\Big)
g_1^{-\frac{1}{4}}f\Big|_{y_n=0}.
\end{equation}	
It follows from  the Lemma \ref{lma:2.1}  that $ g_1, \frac{\partial g_1}{\partial y_n}, \hat g^{+,-},  g_0^{+,j}$  
are known  on $y_n=0$  if $\Lambda'$  is known.
Therefore  knowing  $\Lambda'f$  we can  determine  $\Lambda_1f'$.  Note that 
\begin{equation}																		\label{eq:2.38}
u_1=g_1^{\frac{1}{4}}e^{-id(y)}u(\varphi^{-1}(y)),
\end{equation}
where   $y=\varphi(x)$  is  given  by  (\ref{eq:2.14}),  (\ref{eq:2.15}).

\section{The Green's formula}
\label{section 3}
\init
First,  we introduce some notations.  

Let $\Gamma_1\subset  U_0\cap\Gamma_0$.
Denote  by  $D_{1T_1}$ 
 the forward domain of influence of $\overline \Gamma_1\times [T_1,T_2]$  in the  half-space $y_n\geq 0$.
 We shall  define  $\Gamma_2$  as the intersection  $D_{1T_1}\cap\{y_n=0\}\cap\{y_0=T_2\}$.
 Analogously,  let $\Gamma_3=D_{2T_1}\cap\{y_n=0\}\cap\{y_0=T_2\}$,
 where  $D_{2T_1}$  is the forward  domain  of influence  of $\overline\Gamma_2\times[T_1,T_2]$.  
 We assume  that  $\Gamma_1\subset \Gamma_2\subset\Gamma_3\subset(\Gamma_0\cap U_0)$.
 \\
 \\
 \

\begin{tikzpicture}[scale=1.2]

\draw [dashed](0,0) -- (0,3);
\draw (0,3) -- (2,5);
\draw[line width = 3pt] (1,2) --(3,4);
\draw[line width = 3pt]  (-1,2) -- (3,6);
\draw (0,0) -- (2,2);
\draw (0,0)-- (1,2);
\draw (0,0) -- (2,2);
\draw [dashed](2,2) -- (2,5);
\draw(2,2)--(3,4);

\draw (-2,1)--(3.5,6.5);
\draw[line width = 3pt] (-1,2) .. controls (-0.5,1.5) and (0.5,1.8)..(1,2);
\draw(-1,2) .. controls (-0.8, 0.8) and (-0.5,0.2) .. (0,0);


\draw[line width = 3 pt](3,6) .. controls (3.5,5.7) and (3.2,4.3) .. (3,4);
\draw[dashed](3,6) .. controls (3.2,5.4)  and (2.85,4) .. (2,2);
\draw (1,3) node {$Y_{jT_1}$};
\draw (1.5,5.2) node {$\Gamma_{j+1}$};
\draw(1,0.5) node {$\Gamma_j$};
\end{tikzpicture}
\\
\
\\
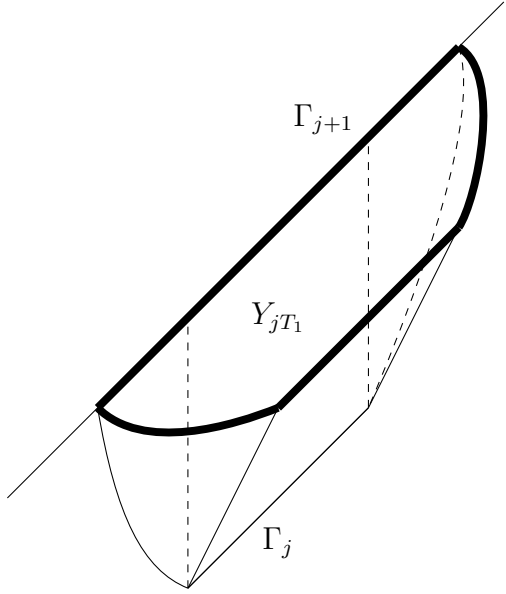
{\bf Fig. 3.1.}  $Y_{jT_1}$  is 
the intersection of the plane $\tau=0$  with  $D_{jT_1}$, 
\\
 $\Gamma_{j+1}$  is the intersection  of $Y_{jT_1}$  with the plane
$y_n=0$.
\\
\\
\
Let  $D_{js_0}$  be the forward domain  of influence  of $\overline\Gamma_j\times[s_0,T_2],1\leq j\leq 3$,  where 
$T_1\leq s_0\leq T_2$.
Denote by $Y_{js_0}$  
the  intersection  of $D_{js_0}$  with  the plane $\tau=T_2-y_n-y_0=0$  (cf. Fig. 3.1).  Let  $X_{js_0}$  be the part of
$D_{js_0}$  below  $Y_{js_0}$  
and  let $Z_{js_0}=\partial X_{js_0}\setminus(Y_{js_0}\cup\{y_n=0\})$.

We assume  also that  $D_{3,T_1}\cap\{y_n=0\}\subset \Gamma_0\cap U_0$   and that  $D_{3T_1}$  does not  intersect
$\Gamma_0\times[T_1,T_2]$  outside  of  $y_n=0$.

Consider the following initial-boundary value problem:
\begin{align}                                                                                  \label{eq:3.1}
&L_1u^f=0
\\
\nonumber
&u^f=u_{y_0}^f=0\ \ \mbox{for}\ \ y_0=T_1,\ y_n>0,
\\
\nonumber
&u^f\big|_{y_n=0}=f,
\end{align}
where $\mbox{supp}\,f\subset \overline\Gamma_3\times[T_1,T_2]$.  Also let $v^g$  be such  that
\begin{align}                                                                                  \label{eq:3.2}
&L_1v^g=0 \ \ \mbox{for}\ \  y_n>0,
\\
\nonumber
&v^g=v_{y_0}^g=0\ \ \mbox{for}\ \ y_0=T_1,\ y_n>0,
\\
\nonumber
&v^g\big|_{y_n=0}=g,
\end{align}
where  $\mbox{supp}\,g\subset\overline\Gamma_3\times[T_1,T_2]$.

Note that  $L_1^*=L_1,$  i.e. $L_1$  is formally  self-adjoint.

Let $(u,v)$  be  the $L_2$  inner  product  $\int_{X_{3T_1}}u(y)\overline v(y)dy$.
We have 
\begin{equation}                            							\label{eq:3.3}
(L_1u^f,v^g)-(u^f,L_1v^g)=0,
\end{equation}
since $L_1u^f=0,\ L_1v^g=0$.
The Jacobian $\frac{\partial(y_n,y_0)}{\partial(s,\tau)}$  is equal 
to $\frac{1}{2}$.   Thus $dy_0dy_n=\frac{1}{2}dsd\tau$.  Integrating  by parts in $s$  we get
\begin{align}   												\label{eq:3.4}
&-\int\limits_{X_{3T_1}}\frac{\partial}{\partial s}\Big(\frac{\partial}{\partial\tau}+iA_-'\Big)u^f\overline{v^g}dsd\tau
\\
\nonumber
&=\int\limits_{X_{3T_1}}\Big(\frac{\partial}{\partial\tau}+iA_-'\Big)u^f\frac{\overline{\partial v^g}}{\partial s}dsd\tau
-\int\limits_{y_n=0}\Big(\frac{\partial}{\partial\tau}+iA_-'\Big)u^f\overline{v^g}d\tau.
\end{align}
Integrating by parts in $\tau$  we get
\begin{align}   												\label{eq:3.5}
&-\int\limits_{X_{3T_1}}\frac{\partial^2}{\partial\tau\partial s}u^f\overline{v^g}dsd\tau
\\
\nonumber
&=\int\limits_{X_{3T_1}}\frac{\partial u^f}{\partial s} \frac{\overline{\partial v^g}}{\partial \tau}dsd\tau
-\int\limits_{y_n=0}\frac{\partial u^f}{\partial s}\overline{v^g} ds
+\int\limits_{\tau=0}\frac{\partial u^f}{\partial s}\overline{v^g}ds.
\end{align}
We used in (\ref{eq:3.4}),  (\ref{eq:3.5}) that $u^f, v^g$  are equal  to zero on $Z_{3T_1}$.   Note that  
$s=y_0-T_1,\tau=T_2-y_0$  on $y_n=0$,  and $\frac{\partial}{\partial s}=\frac{1}{2}\big(\frac{\partial}{\partial y_0}
-\frac{\partial}{\partial y_n}\big),  
\frac{\partial}{\partial \tau}=-\frac{1}{2}\big(\frac{\partial}{\partial y_0}+\frac{\partial}{\partial y_n}\big)$.
Therefore,  making changes of variable  $\tau=T_2-y_0$  in the first integral  and $s=y_0-T_1$  in the second,  we get
\begin{equation}    													\label{eq:3.6}
-\int\limits_{y_n=0}\Big(\frac{\partial}{\partial\tau}+iA_-'\Big)u^f\overline{v^g}d\tau
-\int\limits_{y_n=0}\frac{\partial u^f}{\partial s}\overline{v^g}ds=
\int_{y_n=0}\Big(\frac{\partial}{\partial y_n}-iA_-'\Big)u^f\overline{v^g}dy_0
\end{equation}
Analogously,  integrating  by parts other terms of  
$\int_{X_{3T_1}}(L_1u^f)\overline{v^g}dsd\tau$  we get  (cf.  [E3],  p.316)
\begin{align}														\label{eq:3.7}
0&=(L_1u^f,v^g)-(u^f,L_1v^g)
\\
\nonumber
&=\int\limits_{Y_{3T_1}}\Big(\frac{\partial u^f}{\partial s}\overline{v^g}-u^f\frac{\overline{\partial v^g}}{\partial s}\Big)dy'ds
+\int\limits_{\Gamma_3\times[T_1,T_2]}(\Lambda_1f\overline g-f\overline{\Lambda_{1}g})dy'dy_0,
\end{align}
where $\Lambda_1$ has the  form  (\ref{eq:2.30})
\begin{equation} 																	\label{eq:3.8}
\Lambda_1f=\Big(\frac{\partial u^f}{\partial y_n}-iA_-'u^f\Big)
-\sum_{j=1}^{n-1}g_0^{+,j}\Big(\frac{\partial u^f}{\partial y_j}-iA_j'u^f\Big)\Big|_{y_n=0}.
\end{equation}
It follows  from  (\ref{eq:3.7})  that 
\begin{equation}   																\label{eq:3.9}
\int\limits_{Y_{3T_1}}(u_s^f\overline{v^g}-u^f\overline{v_s^g})dsdy'
\end{equation}
is determined by the boundary data,  i.e.   by the DN  operator  on  $\Gamma^{(3)}\times(T_1,T_2)$.

We shall denote  the $L_2$  inner product  in $Y_{3T_1}$  by  $(u,v)_{Y_{3T_1}}$,  or simply  $(u,v)$
when  it is clear  what  is the domain of integration.

Let $D_j^-$  be the backward domain of influence of $\overline\Gamma_j\times[T_1,T_2]$.  
Thus   $ D_{jT1}\cap D_j^-$  is the  domain  of dependence  of  $\overline\Gamma_j\times[T_1,T_2]$.
Denote by  $Q_j$  the intersection  of $D_j^-$  
with $\tau=0$.  
 Let $R_{js_0}=Y_{js_0}\cap Q_j$  be the rectangle
$\{s_0-T_1\leq s\leq T_2-T_1,\tau=0,y'\in\overline\Gamma_j\}$.
Note that  $R_{js_0}$  belongs  to the domain  of dependence of $\overline\Gamma_j\times[s_0,T_2]$. 
Let $H_0^1(R_{js_0})$ be the subspace of the Sobolev space $H^1(R_{js_0})$  
  consisting of 
$w\in H^1(R_{js_0})$  such that $w=0$  on $\partial R_{j0}\setminus\{y_n=0\}$.  Analogously,  let  $H_0^1(Y_{js_0})$ 
 be the subspace of $H^1(Y_{js_0})$ 
 consisting  of $v\in H^1(Y_{js_0})$  such that $v=0$  on $\partial Y_{js_0}\setminus\{y_n=0\}$.   Note
that $R_{js_0}\subset Y_{js_0}\subset R_{j+1,s_0}$  (cf.  Fig. 3.2).
\\
\

\begin{tikzpicture}
\draw(0,0)-- (0,4);
\draw(0,0)-- (7,0);
\draw(0,4)-- (7,4);
\draw(7,0)-- (7,4);
\draw(2.5,4)--(2.5,0);
\draw(4.5,4)--(4.5,0);

\draw[ultra thick] (1,4) --(6,4);
\draw[ultra thick] (2.5,0)--(4.5,0);

\draw[ultra thick] (1,4) ..controls (0.8,2.5) and (1,1.5) .. (2.5,0);
\draw[ultra thick] (6,4) ..controls (6.2,2.5) and (6,1.5) .. (4.5,0);
\draw(3.5,-0.3) node {$\Gamma_j$};
\draw(4,4.3)  node {$\Gamma_{j+1}$};
\draw(3.5,2) node {$R_{js_0}$};
\draw(2,2.5) node {$Y_{js_0}$};
\draw(5,2.5) node {$Y_{js_0}$};

\draw(-1.2,4) node {$s=T_2-T_1$};

\draw(-1.2,0) node {$s=s_0-T_1$};

\end{tikzpicture}
\\
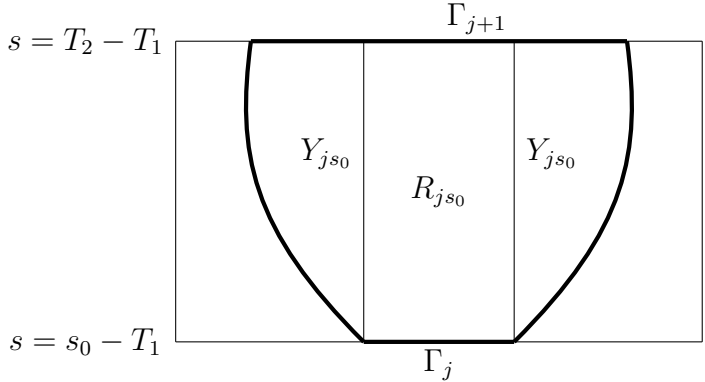
{\bf Fig.3.2.}   The rectangle 
$R_{js_0} =\{s_0-T_1\leq s\leq T_2-T_1,\tau=0,y'\in \Gamma_j\}$ , \linebreak 
$Y_{js_0}$ is the intersection  of the domain of
influence of  $[s_0,T_2]\times\overline\Gamma_j$ \\ with the plane $\tau=0$.  Note that 
$R_{js_0}\subset Y_{js_0}\subset R_{j+1.s_0}$.
\\
\

Note that $H_0^1(R_{js_0})$  is a subspace  of $H_0^1(Y_{js_0})$.
\begin{lemma}[Density lemma]         														\label{lma:3.1}
For any $w\in H_0^1(R_{js_0})$   there exists a sequence $\{u^{f_n}\}$  where $u^{f_n}$  are solutions of the initial-boundary 
value problem  (\ref{eq:3.1}),  $f_n(y_0,y')\in  H_0^1(\Gamma_j\times[s_0,T_2])$,   such that  $\|w-u^{f_n}\|_{1,Y_{js_0}}\rw 0$
when  $n\rw\infty, \ j=1,2,3.$
\end{lemma}
Here $\|w\|_{1,Y_{js_0}}$  is the norm  in $H_0^1(Y_{js_0})$    and
$f\in H_0^1(\Gamma_j\times[s_0,T_2])$,  i.e.  $f=0$  on $\partial(\Gamma_j\times[s_0,T_2])\setminus(\Gamma_j\times\{y_0=T_2\})$.
\\
\

{\bf Proof:}  The proof of Lemma \ref{lma:3.1} is a simplification of the proof 
of Lemmas 2.2  and 3.2  in  [E3].  We shall prove  Lemma \ref{lma:3.1}  for the case $s_0=T_1$.
The proof for the case $T_1<s_0<T_2$  is identical.
\\
\

\begin{tikzpicture}[scale=1.]

\draw (0,-3.5)--(0,0);
\draw(-0.5,0)--(3.9,0);
\draw(0,0)--(3,-3);
\draw(1,-2) node {$\Delta_2'$};
\draw (1.5,-0.3) node {$\tau<0$};
\draw (2.5,-1.4) node {$\Gamma_2'$};
\draw (-0.5,-1.5) node {$\Gamma_1'$};

\draw(-0.1,0.3) node {$y_0=T_2$};

\end{tikzpicture}
\
\\
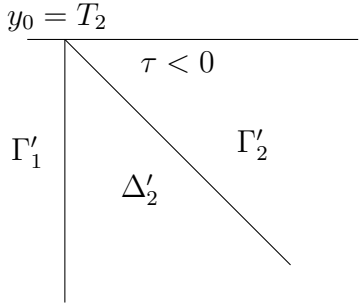
{\bf Fig. 3.3.}  The   domain $\Delta_2'$  is bounded  by 
 $\Gamma_1'$  and  $\Gamma_2'$.
\\
\

Denote by  $\Delta_2'$  the domain bounded  by the half-plaines
$\Gamma_1'=\{y_n=0,y_0<T_2,y'\in\R^{n-1}\}$  and
$\Gamma_2'=\{\tau=T_2-y_n-y_0=0,s<T_2,y'\in\R^{n-1}\}$  (cf. Fig. 3.3).   Let 
$\Gamma_\infty'$  be  the plane $\tau=0$.    Denote  by  $H_0^{-1}(\Gamma_2')$  the Sobolev  space 
 of $h\in H^{-1}(\Gamma_\infty')$
such  that  
$\mbox{supp}\,h\subset\overline\Gamma_2'$,  i.e. $h(s,y')=0$   when  $s>T_2$.  
Note that  $H_0^{-1}(\Gamma_2')$  is dual  to $H^1(\Gamma_2')$  with  respect  to the  extension of the  $L_2$  inner  product
on $\Gamma_\infty'$  (cf.  [E5]).
\begin{lemma}																	\label{lma:3.2}
For any  $h(s,y')\in  H_0^{-1}(\Gamma_2')$  there exists  a distribution  $u(s,\tau,y')$   such that
\begin{align}																	\label{eq:3.10}
&L_1u=0\ \ \ \mbox{in}\ \ \Delta_2',
\\
&\frac{\partial u}{\partial s}\Big|_{\Gamma_2'}=h,									\label{eq:3.11}		
\\
&u\big|_{y_n=0}=0.																\label{eq:3.12}
\end{align}
\end{lemma}
{\bf Proof:}  
Since  $h(s,y')=0$  for  $s>T_2$,  there exists  $v(s,y')=0$  for  $s>T_2,  v(s,y')$  belongs  to $L_2$  in $s$  and  to $H^{-1}$ 
in $y'$  and such  that  $\frac{\partial v}{\partial s}=h$  in  $\Gamma_\infty'$.    We can  define  $v(s,y')$   by the formula 
$$
v(s,y')=\lim_{\e\rw 0}e^{\e(s-T_2)}F^{-1}\frac{\tilde h(z_0,\xi_1,...,\xi_n)}{z_0+i\e},
$$
where  $\tilde h(z_0,\xi_1,...,\xi_{n-1})$  is the Fourier 
transform  of $h(s,y')$  and
 $F^{-1}$  is the inverse  Fourier transform,  $z_0$  is the dual  variable  to $s$.
         
The distribution  $\theta(-\tau)u$  satisfies  the equation
\begin{equation}																\label{eq:3.13}
L_1(\theta(-\tau)u)=4h\delta(-\tau)
\end{equation}
in the half-space  $y_n>0$  with  the boundary  condition
\begin{equation} 																\label{eq:3.14}
\theta(-\tau)u\big|_{y_n=0}=0,
\end{equation}
where $\theta(s)=1$  for  $s>0$   and  $\theta(s)=0$  for  $s<0$.

We  look  for  $\theta(-\tau)u$  in the  form
\begin{equation}                                                                							\label{eq:3.15}
\theta(-\tau)u=\theta(-\tau)v+w,
\end{equation}
where $w$   satisfies 
\begin{align}																	\label{eq:3.16}
&L_1w=\varphi,
\\
&w\big|_{y_n=0}=-\theta(-\tau)v\big|_{y_n=0},									\label{eq:3.17}		
\\
&\varphi=L_{11}(\theta(-\tau)v),										\label{eq:3.18}
\end{align}
where  $L_{11}=L_1+\frac{4\partial^2}{\partial s\partial \tau}$.  
 Note  that  $L_{11}$  is a differential  operator  in $\frac{\partial}{\partial s},\frac{\partial}{\partial y_k},1\leq k\leq n-1$.

We  impose  the zero  initial  conditions  on $w$  requiring that
\begin{equation}															\label{eq:3.19}
w=0\ \ \ \mbox{for}\ \ \ y_0>T_2.
\end{equation}
Therefore  $w$  is the solution  of  the hyperbolic equation  $L_1w=\varphi$ in the  half-space   $y_n>0$  with  the boundary  
condition  (\ref{eq:3.17})  and the zero  initial conditions (\ref{eq:3.19}).
  It follows  from  ([H]  and [E6])  that  initial-boundary  value problem  has a unique solution  
in appropriate  Sobolev  space  of negative order.

Since  $\varphi $  belongs  to  $L_2$  
in $\tau$  and to  Sobolev spaces of negative order in  $s$  and  $y'$,   we get  that  $w$  belongs 
to  $H^1$  in  $\tau$.
Therefore  $w\big|_{\tau=\tau_0}$  is continuous  function   of  $\tau_0$  with the values  in Sobolev's spaces  of negative  order  in 
$(s,y')$.
Since  $\varphi=0$  for  $\tau<0$  we have  that  $w=0$  for  $\tau<0$  by the domain of influence argument.
Therefore by the  continuity 
$w\big|_{\tau=0}=0$  and  $\frac{\partial w}{\partial s}\big|_{\tau=0}=0$.

Therefore  $u=v(s)+w(s,\tau,y')$  is  the distribution  solution  of (\ref{eq:3.10}),  (\ref{eq:3.11}),  (\ref{eq:3.12})   in $\Delta_2'$.

Note  that the restrictions of  any distribution  solution  of $L_1u=0$  to $y_n=0$  exists since  $y_n=0$  is not  a characteristic  surface  for  $L_1$.   This property  is called the partial  hypoellipticity  (cf.,  for example,  [E5]).
\qed

Now using  Lemma  \ref{lma:3.2}   we can  prove  Lemma \ref{lma:3.1}.  
If  $\{u^f,f\in H_0^1(\Gamma_{jT_1}),
\linebreak
\Gamma_{jT_1}=\Gamma_j\times[T_1,T_2]\}$  is not dense  in
$H_0^1(R_{jT_1})$  then there exists  nonzero $h\in  H_0^{-1}(R_{jT_1})$  such 
 that $(u^f,h)=0,\forall u^f,  f\in H_0^1(\Gamma_{jT_1})$.
Let  $v$  be  such  that  $\frac{\partial v}{\partial s}=h$.   Then
$\big(u^f,\frac{\partial v}{\partial s}\big)_{Y_{jT_1}}=0,
\forall f\in H_0^1(\Gamma_{jT_1})$.
Let $u$  be the same  
as in Lemma \ref{lma:3.2},  i.e. $L_1u=0$  in $\Delta_2'$,  $u\big|_{Y_{T_1}}=v,\ u\big|_{y_n=0}=0$.  Then
\begin{equation}    																				\label{eq:3.20}
-2\Big(u^f,\frac{\partial v}{\partial s}\Big)_{Y_{jT_1}}=\int\limits_{\Gamma_{jT_1}}f\frac{\overline{\partial u}}{\partial y_n}dy'dy_0,\ \ 
\forall f\in H_0^1(\Gamma_{jT_1}).
\end{equation}
Note that  $(u^f,h)_{Y_{jT_1}}$  is understood as the extension  of $L_2$  inner  product  in  
$(u_1,h)$  in  $\Gamma_\infty'$,  where $u_1$  is an  arbitrary   extension  of $u^f$  for $s>T_2$.   Analogously
for the right hand  side of (\ref{eq:3.20}).  (Note that  $u=0$  for  $y_0>T_2)$.

To justify  (\ref{eq:3.20})  we take a sequence $h_j\in C_0^\infty(\Gamma_2'), \ h_j\rw h$  in  $H_0^{-1}(\Gamma_2')$.  
By Lemma \ref{lma:3.2}  there exists  smooth $v_j$  such that  $L_1v_j=0$  in $\Delta_2',v_j\big|_{y_n=0}=0,\linebreak
\frac{\partial v_j}{\partial s}\big|_{\tau=0}=h_j$.
Applying  the Green's formula (\ref{eq:3.7})  to $u^f$  and  $v_j$  we get
$$
\int\limits_{Y_{3T_1}}\Big(\frac{\partial u^f}{\partial s}\overline v_j-u^f\frac{\overline{\partial v_j}}{\partial s}\Big)dy'ds=
\int\limits_{\Gamma_{3T_1}}f\frac{\overline{\partial v_j}}{\partial y_n}dy'dy_0,
$$
since  $v_j\big|_{y_n=0}=0.$  Integrating  by parts  we get
$$
\int\limits_{Y_{3T_1}}\frac{\partial u^f}{\partial s}\overline v_jdsdy'
=-\int\limits_{Y_{3T_1}}u^f\frac{\overline{\partial v_j}}{\partial s}dsdy' +u^f\overline v_j\big|_{y_n=0,y_0=T_2}.
$$
Since $v_j=0$  for  $s>T_2$ we have  
 that $v_j\big|_{y_n=0,y_0=T_2}=0$.  Therefore,  taking the
limit when $j\rw \infty$  we get (\ref{eq:3.20}).  

   Since  $f$ is arbitrary  and 
 $\big(u^f,\frac{\partial v}{\partial s}\big)=0$  we get  that  $\frac{\partial u}{\partial y_n}=0$  on  $\Gamma_{jT_1}$.
Therefore  $L_1u=0$  in $\Delta_2'$  and  $u$  has zero  Cauchy data  on $\Gamma_{2T_1}$.  Then  $u=0$  in the domain of
dependence
 $D_j^-\cap D_{jT_1}$,  in particular,  $u=0$  on  $R_{jT_1}$.   Thus  $v=0$  on   $R_{jT_1}$  and this  contradicts 
 the assumption  that  $h\neq 0$.
 \qed
 
 We shall prove  two  more  theorems  in this  section  that will be  used  in \S 4.     

We shall need  some  known  results on the initial-boundary  hyperbolic  problem.  The following theorem holds:
\begin{lemma}     																	\label{lma:3.3}
Let  $L_1u=F$  in  $\R_+^n\times(-\infty,T_2)$  where $F\in  H_+^s(\R^n\times(-\infty,T_2),\R_+^n=\{y_n>0,y'\in\R^{n-1}\}$.
Let  $u\big|_{y_n=0}=f$,  where $f\in H_+^{s+1}(\R^{n-1}\times(-\infty,T_2))$. 
Then  for any  $s\geq 0$  and  any  $f\in  H_+^{s+1}(\R^{n-1}\times(-\infty,T_2))$  and
$F\in  H_+^s(\R_+^n\times(-\infty,T_2))$   there exists  a unique  $u\in H_+^{s+1}(\R_+^n\times(-\infty,T_2)) $
 such that 
\begin{equation}																	\label{eq:3.21}
\|u\|_{s+1}\leq C([f]_{s+1}+\|F\|_s).
\end{equation}
Moreover,
\begin{equation}																		\label{eq:3.22}
\Big[\frac{\partial u(y_0,y',0)}{\partial y_n}\Big]_s \leq C(\|F\|_s+[f]_{s+1}).
\end{equation}
Here  $H_+^s(\R_+^n\times (-\infty,T_2))$  is the  Sobolev's  space  $H^s(\R_+^n\times(-\infty,T_2))$  with norm  $\|u\|_s$
consisting  of $u(y)$   with   the support  in  $y_0\geq  T_1,\ [f]_s$  is the norm  in $H_+^s(\R^{n-1}\times(-\infty,T_2))$.

We assume  that  $F(y)$  and  $f(y_0,y')$  have compact  supports  in $y'$.   Note that $f=0,\ F=0,\ u=0$   for  $y_0<T_1$.

Then  $u(y_0,y',y_n)$  has  also  a compact support in $y'$.
\end{lemma}
The proof of 
Lemma \ref{lma:3.3}  in  the case  of  time-dependent  coefficients   is given  in [H]  and [E6].

Note  that Lemma \ref{lma:3.3}  holds  also  in the case  when  $\R_+^n$  is replaced by an arbitrary  smooth domain  
$\Omega\subset\R^n$.

The following lemma  follows   from Lemma \ref{lma:3.3}.
\begin{lemma}																			\label{lma:3.4}
Let,  for  the simplicity,  $F=0$.  For any  $f\in  H_+^1(\R^{n-1}\times(-\infty,T_2]),f=0$  for  $y_0\leq T_1$,  there exists  
$u\in  C(H^1(\R_+^n),[T_1,T_2])\cap C^1(L_2(\R_+^n), [T_1,T_2])$  such that  
$L_1u=0$  in $\R_+^n\times(-\infty,T_2], u=0$  for  $y_0<T_1$,
\begin{equation}																		\label{eq:3.23}
\max_{T_1\leq y_0\leq T_2}\|u(y_0,y',y_n)\|_1^2
+\max_{T_1\leq y_0\leq T_2}\|\frac{\partial u}{\partial  y_0}(y_0,y',y_n)\|_0^2\leq  C[f]_1^2.
\end{equation}
Here  $C(H^1(\R_+^n),[T_1,T_2])\cap C^1(L_2(\R_+^n),(T_1,T_2))$
  means  that  $u(y),\frac{\partial u(y)}{\partial y_0}$  are continuous  functions of  $y_0$ 
with values in  $H^1(\R_+^n), L_2(\R_+^n)$,  respectfively.
\end{lemma}
{\bf Proof:}   Take  $s>\frac{3}{2}$.   Consider  the equation 
$L_1u=0$
in  $\R_+^n\times(-\infty,T_2), u=0$  for   $y_0<T_1$,   using  $(y_0,y',y_n)$   coordinates  (cf.  (\ref{eq:2.28}).   We have 
\begin{equation}															\label{eq:3.24}
g_0^{00}=-g_0^{nn}=1,\ \ \ \ \  g_0^{0j}=g_0^{nj}, \ \ \ \ \ A_0'=A_n'=A_-'.
\end{equation}
Let  $(u,v)_{T'}$  be  the  $L_2$-inner  product  in  $\R_+^{n+1}\times(T_1,T'),T'\leq  T_2$.
Integrating  by  parts   the identity
$$
0=(L_1u,u_{y_0})_{T'}+(u_{y_0}L_1u)_{T'},
$$
we get  (cf.  [E4])
\begin{equation} 															\label{eq:3.25}
E_{T'}(u,u)+\Lambda_0(f,f)+I_1=0,
\end{equation}
where
\begin{multline}												 			\label{eq:3.26}
E_{T'}(u,u)
\\
=\int\limits_{\R_+^n}\Big(|u_{y_0}-iA_0'u|^2
-\sum_{j,k=1}^ng_0^{jk}(u_{y_j}-iA_j'u)\overline{(u_{y_k}-iA_ku')}+V_1|u|^2\Big)
dy'dy_n\big|_{y_0=T'},
\end{multline}
\begin{equation}															\label{eq:3.27}
\Lambda_0(f,f)=\int\limits_{T_1}^{T_2}\int\limits_{\R^{n-1}}
\big[(\Lambda_1f)\overline{f_{y_0}} +f_{y_0}\overline{\Lambda_1f}\big]dy'dy_0,
\end{equation}
$\Lambda_1f$  is  the same  as  in  (\ref{eq:2.36}),
\begin{equation}																\label{eq:3.28}
|I_1|\leq  C\int\limits_{\R_+^n\times[T_1,T']} \sum_{k=0}^n|u_{y_k}|^2dy_0dy'dy_n.
\end{equation}
Note   that  $I_1=0$   when  the coefficients  of  $L_1$  do not depend  on $y_0$.

Let  $\|u\|_{s,T'}$  be the  norm  in $H^s(\R_+^n)$  when  $y_0=T'$.  We have
\begin{equation}            															\label{eq:3.29}
|I_1|\leq C\int\limits_{T_1}^{T'}(\|u\|_{1,t}^2+\|u_{y_0}\|_{0,t}^2)dt\leq C\int\limits_{T_1}^{T'}|[u]|_t^2dt
\end{equation}
where
\begin{equation}																\label{eq:3.30}
|[u]|_t^2=\|u\|_{1,t}^2+\|u_{y_0}\|_{0,t}^2.
\end{equation}
We have 
\begin{equation} 																\label{eq:3.31}
E_{T'}(u,u)\geq C|[u]|_{1,T'}^2
\end{equation}
if $T_2-T_1$  is  small.

Since  $T_2-T_1$  is small,  (\ref{eq:3.25})  implies
\begin{equation}																\label{eq:3.32}
\max_{T_1\leq T'\leq T_2}|[u]|_{1,T'}^2\leq C(T_2-T_1)(\max_{T_1\leq T'\leq T_2} |[u]|_{1,T'}^2+|\Lambda_0(f,f)|)
\end{equation}
Note that 
\begin{equation} 																\label{eq:3.33}
|\Lambda_0(f,f)|\leq C\Big([f]_1^2+\Big[\frac{\partial u(y_0,y',0)}{\partial y_n}\Big]_0\Big).
\end{equation}
Therefore (\ref{eq:3.22}), (\ref{eq:3.25}),  (\ref{eq:3.31}),  (\ref{eq:3.32}),  (\ref{eq:3.33})  imply  (\ref{eq:3.23}).

Since  $H_+^s$  is  dense in  $H_+^1$  when $s>1$  we can approximate  $f\in  H_+^1(\R^{n-1}\times (-\infty,T_2))$  
by functions from  $H_+^s(\R^{n-1}\times(-\infty,T_2)),s>\frac{3}{2}$   and  therefore   the inequality  (\ref{eq:3.23})   holds
for  $f\in H_+^1$. 
\qed

We shall study the Goursat  problem  (see Fig. 3.4):
\\
\

\begin{tikzpicture}[scale=1.5]

\draw[->] (0,-3.5)--(0,0.9);
\draw[->] (-0.5,0)--(3.9,0);

\draw(1.5,-0.5) node {$\Delta_1$};

\draw(0,-3)--(3,0);
\draw(0,0)--(1.5,-1.5);
\draw (2.5,-1.6) node {$(\frac{T_1+T_2}{2},\frac{T_2-T_1}{2})$};

\draw (-0.3,0.8) node {$y_0$};
\draw (4,-0.3) node {$y_n$};

\draw(-0.3,-3.1)  node {$T_1$};
\draw(-0.3,-0.25) node {$T_2$};

\draw(1.0,0.2) node {$\Gamma_4$}; 
\draw (3.1,0.2) node {$T_2-T_1$};

\draw (1.0,-1.3) node {$\Gamma_2$};
\draw (2.5,-0.8) node {$\Gamma_3$};
\filldraw[gray] (0,0) circle (1.5pt)
                      (1.5,-1.5) circle (1.5pt)
                       (3,0) circle (1.5pt)
                       (0,-3) circle (1.5pt);

\end{tikzpicture}
\
\\
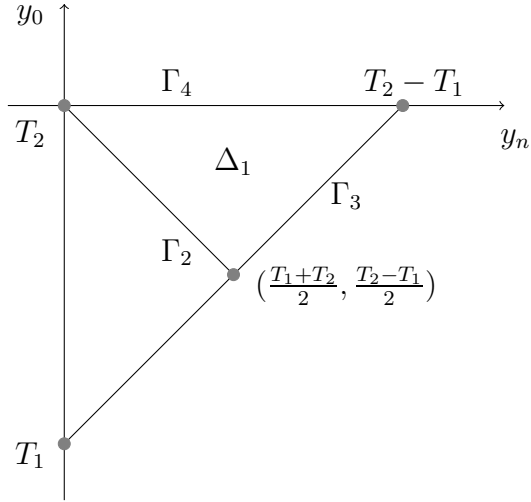
{\bf Fig. 3.4.}  Domain $\Delta_1$  is bounded  by the planes $\Gamma_2,\Gamma_3,\Gamma_4$.
\\
\

We use the same notations that we used in the proof of the Lemma 3.1 in [E3]:
Let  $\Gamma_2,\Gamma_3$  and  $\Gamma_4$  be the following  planes:

 $\Gamma_2=\{\tau=T_2-y_0-y_n=0, 0\leq y_n\leq \frac{T_2-T_1}{2},y'\in\R^{n-1}\},$
 
 $\Gamma_3=\{s=y_0-y_n-T_1=0, \frac{T_2+T_1}{2}\leq y_0\leq T_2,y'\in\R^{n-1}\},$
 
 $\Gamma_4=\{y_0=T_2, 0\leq y_n\leq T_2-T_1,y'\in\R^{n-1}\}$.
\\
Let $\Delta_1$  be the domain  bounded  by $\Gamma_2,\Gamma_3,\Gamma_4$.  The following lemma is similar  to  
Lemma 3.1  in [E3]:
\begin{lemma}															\label{lma:3.5}
For any  $v_0\in H^1(\Gamma_4), v_1\in L_2(\Gamma_4)$  there exists $u\in H^1(\Delta_1)$  such that 
$L_1u=0$  in $\Delta_1,\ u\big|_{\Gamma_4}=v_0,u_{y_0}\big|_{\Gamma_4}=v_1$.   Moreover,
the traces  $\varphi=u\big|_{\Gamma_2},\psi=u\big|_{\Gamma_3}$  exists and belongs  to $H^1(\Gamma_2),H^1(\Gamma_3)$,
respectively.  The following estimate holds:
\begin{equation}														\label{eq:3.34}
\big\|u\big|_{\Gamma_2}\big\|_{1,\Gamma_2}^2+\big\|u\big|_{\Gamma_3}\big\|_{1,\Gamma_3}^2
\leq C\Big(\big\|u\big|_{\Gamma_4}\big\|_{1,\Gamma_4}^2+\big\|u_{y_0}\big|_{\Gamma_4}\big\|_{0,\Gamma_4}^2\Big).
\end{equation}
Vice versa,  for any  $\varphi\in H^1(\Gamma_2),\psi\in H^1(\Gamma_3),\varphi=\psi$
  at $y_0=\frac{T_2+T_1}{2}$   there exists  $u\in  H^1(\Delta_1),\ L_1u=0$  in  $\Delta_1$  such  that 
   $u\big|_{\Gamma_2}=\varphi,\ u\big|_{\Gamma_3}=\psi$  and  the following estimate holds:
  \begin{equation}														\label{eq:3.35}
\big\|u\big|_{\Gamma_4}\big\|_{1,\Gamma_4}^2+\big\|u_{y_0}\big|_{\Gamma_4}\big\|_{0,\Gamma_4}^2
\leq C\Big(\big\|u\big|_{\Gamma_2}\big\|_{1,\Gamma_2}^2+\big\|u\big|_{\Gamma_3}\big\|_{1,\Gamma_3}^2\Big).
\end{equation}
\end{lemma}
{\bf Proof:}
Let $\Delta_{1,T'}$  be the domain  bounded  by  $\Gamma_2,\Gamma_3$  and  $\Gamma_{4,T'}$,   where 
$\Gamma_{4,T'}$  is the plane  $y_0=T',\frac{T_1+T_2}{2}\leq T'\leq T_2$.  Denote  by  $(u,v)_{\Delta_{1,T'}}$
the $L_2$-inner product
in  $\Delta_{1,T'}$.   Integrating  by  parts  the identity
$$
(L_1u,u_{x_0})_{\Delta_{1T'}}+(u_{x_0},L_1u)_{\Delta_{1T'}}=0
$$
we get,  as  in [E4]:
\begin{equation}														\label{eq:3.36}
E_{T'}(u,u)+Q_{T'}(u,u)+Q_{T'}^{(1)}(u,u)=I_2,
\end{equation}
where  $E_{T'}(u,u)$  is  the same as in  (\ref{eq:3.26}),
\begin{align}															\label{eq:3.37}
&Q_{T'}(u,u)=\frac{1}{2}\int\limits_{\Gamma_{2T'}}\Big[4|u_s|^2-
\sum_{j,k=1}^{n-1}g_0^{jk}\Big(\frac{\partial u}{\partial y_j}-iA_j'u\Big)\overline{\Big(\frac{\partial u}{\partial y_k}-iA_k'u\Big)}
\\
\nonumber
&-2\sum_{j=1}^{n-1}\Bigg(g_0^{0j}\Big(\frac{\partial u}{\partial s}+iA_-'u\Big)
\overline{\Big(\frac{\partial  u}{\partial y_j}-iA_j'u\Big)}
+g_0^{0j}\Big(\frac{\partial u}{\partial y_j}-iA_j'u\Big)\overline{\Big(\frac{\partial u}{\partial s}+iA_-'u\Big)}\Bigg)
\\
\nonumber
&+V_1|u|^2\Big]dy'ds,
\end{align}
(cf.  (3.22)  in  [E4]),
\begin{multline}																		\label{eq:3.38}
Q_{T'}^{(1)}(u,u)
\\
=\frac{1}{2}\int\limits_{\Gamma_{3T'}}\Bigg(|u_\tau+iA_-'u|^2
-\sum_{j,k=1}^{n-1}g_0^{jk}(u_{y_j}-iA_j'u)\overline{(u_{y_k}-iA_k'u)}+V_1|u|^2\Bigg)dy'd\tau,
\end{multline}
\begin{equation}																		\label{eq:3.39}
|I_2|\leq C\int\limits_{\Delta_{1T'}}\sum_{j=0}^{n}\Big|\frac{\partial u}{\partial y_j}\Big|^2dy_0dy'dy_n.
\end{equation}
Here 
$\Gamma_{2T'}, \Gamma_{3T'}$  are  parts of $\Gamma_2,\Gamma_3$  for   $\frac{T_1+T_2}{2}\leq T'$.
When  $T_2-T_1$  is small,  $Q_{T'}(u,u)$  is positive definite (cf.  [E4],  (3.23)).  Therefore
\begin{equation}																		\label{eq:3.40}
C_1\|u\|_{1,\Gamma_{2T'}}^2\leq Q_{T'}(u,u)\leq C_2\|u\|_{1,\Gamma_{2T'}}^2.
\end{equation}
Analogously,
\begin{equation}																		\label{eq:3.41}
C_1'\|u\|_{1,\Gamma_{3T'}}^2\leq Q_{T'}^{(1)}(u,u)\leq C_2'\|u\|_{1,\Gamma_{3T'}}^2.
\end{equation}    
Having  (\ref{eq:3.31}),  (\ref{eq:3.39}),  (\ref{eq:3.40}),  (\ref{eq:3.41})  we can complete the proof   of Lemma \ref{lma:3.5}  exactly as
the proof of Lemma 3.1  in [E3].
\qed

Combining Lemmas \ref{lma:3.4}  and \ref{lma:3.5}  we can  prove the following  lemma:
\begin{lemma}																			\label{lma:3.6}
The map  $f\rw u^f$  is a bounded operator  from  $H_0^1(\Gamma_j\times[s_0,T_2])$ to  $H_0^1(Y_{js_0})$:
\begin{equation}																		\label{eq:3.42}
\|u^f\|_{1,Y_{js_0}}\leq C[f]_1.
\end{equation}
\end{lemma}
{\bf Proof:}
It follows  from Lemma \ref{lma:3.4} that
\begin{equation}																		\label{eq:3.43}
\|u^f\|_{1,\Gamma_4}^2+\|u_{y_0}^f\|_{0,\Gamma_4}^2\leq C[f]_{1}^2 .
\end{equation}
Then  (\ref{eq:3.34})  gives
\begin{equation}     																		\label{eq:3.44}
\|u^f\|_{1,\Gamma_2}^2\leq C\big(\|u^f\|_{1,\Gamma_4}+\|u_{y_0}^f\|_{0,\Gamma_4}^2\big).
\end{equation}
Combining  (\ref{eq:3.43})  and (\ref{eq:3.44})
  and taking  into account  that $\mbox{supp}\, u^f\big|_{\Gamma_2}=Y_{js_0}$,  we get (\ref{eq:3.42}).

\section{The Main formula}
\label{section 4}
\init
Let  $L_1^{(i)},i=1,2,$  be two operators of the form  (\ref{eq:2.27})  such that the corresponding  DN  operators  $\Lambda_1^{(1)}$  
and  $\Lambda_1^{(2)}$  
are equal on  $U_0\cap \{y_n=0\}$.  We choose  $\Gamma_1^{(1)}=\Gamma_1^{(2)}=\Gamma_1$  in 
a neighborhood of $x^{(0)}$  in  $U_0\cap \{y_n=0\}$.   Let  $\Gamma_j^{(i)}, j=2,3,i=1,2,$  be defined as before  (see Fig. 3.1)
for  $i=1,2,$  respectively.
\begin{lemma}																			\label{lma:4.1}
We have  $\Gamma_j^{(1)}=\Gamma_j^{(2)},j=2,3$  (cf. [E1]).
\end{lemma}
{\bf Proof:}
Let  $\Delta_{2T_1}^{(i)}$  be  the intersection  of the domain of influence $D_{2T_1}^{(i)},i=1,2,$
with the plane  $y_n=0$.  Note that  $\Delta_{2T_1}^{(i)}$ 
is the  intersection of $y_n=0$  with the closure  of the union  $\bigcup\mbox{supp}\, u_i^f$  where
the union is taken  over all  $f\in  H_0^1(\Gamma_1\times[T_1,T_2]), L_1^{(i)}u_i^f=0$.

Let  $\tilde\Delta_{2T_1}^{(i)}$  be the closure of  the union  $\bigcup \mbox{supp}\,\Lambda_1^{(i)}f$,  where  the union  is
taken  also over  all $f\in  H_0^1(\Gamma_1\times [T_1,T_2])$.   We shall show that 
$\tilde\Delta_{2T_1}^{(i)}=\Delta_{2T_1}^{(i)}$.

If $x^{(0)}\not\in \Delta_{2T_1}^{(i)}$
then  $u_i^f=0,\ \forall f$,  in some neighborhood of  $x^{(0)}$  in  $U_0$.   Then  $\Lambda_1^{(i)} f=0$  in a neighborhood  
of  $x^{(0)},\ \forall  f$.  Thus  $x^{(0)}\not\in\tilde\Delta_{2T_1}^{(i)}$,
i.e.
$\tilde\Delta_{2T_1}^{(i)}\subset \Delta_{2T_1}^{(i)}$.  Let now  $x_0'\not\in \tilde\Delta_{2T_1}^{(i)}$.
Then  $\Lambda_1^{(i)}f=0$  in a neighborhood  of $x_0'$
for any  $f\in  H_0^1(\Gamma_1\times[T_1,T_2])$  and also  $f=0$  in a neighborhood 
of  $x_0'$.  Then  by the uniqueness  of the Cauchy problem  
(see  [T], [RZ]) we have that all 
$u^f=0$  in a neighborhood of $x_0'$  in $\R^{n+1}$.   Therefore  $x_0'\not\in \Delta_{2T_1}^{(i)}$.  Thus   
$\Delta_{2T_1}^{(1)}=\tilde\Delta_{2T_1}^{(1)}$.  Since $\Lambda_1^{(1)}=\Lambda_1^{(2)}$,  we have
$\tilde\Delta_{2T_1}^{(1)}=\tilde\Delta_{2T_1}^{(2)}$.  Therefore $\Delta_{2T_1}^{(1)}=\Delta_{2T_1}^{(2)}$,
i.e. $\Gamma_2^{(1)}=\Gamma_2^{(2)}$.  Analogously one shows that  $\Gamma_3^{(1)}=\Gamma_3^{(2)}$.
\qed

Since  $\Gamma_j^{(1)}=\Gamma_j^{(2)}$  we shall write  $\Gamma_j, 1\leq j\leq 3,$  instead  of $\Gamma_j^{(i)}$.  
It follows  from (\ref{eq:3.7})  that  (\ref{eq:3.9})  is  determined  by the boundary  data.  Integrating  by part  we have
\begin{align}																	\label{eq:4.1}
&\int\limits_{Y_{3T_1}}(u_s^f\overline{v^g}-u^f\overline{v_s^g})dsdy'
\\
\nonumber
=2&\int\limits_{Y_{3T_1}}u_s^f\overline{v^g}dsdy'
-\int\limits_{\partial Y_{3T_1}\cap\{y_n=0\}} u^f(T_2-T_1,0,y')\overline{v^g}(T_2-T_1,0,y')dy'.
\end{align}
Since $u^f(T_2-T_1,0,y')=f(T_2,y'),\ v^g(T_2-T_1,0,y')=g(T_2,y')$,  we  have that  
\begin{equation}																\label{eq:4.2}
(u_s^f,v^g)=\int\limits_{Y_{3T_1}}u_s^f\overline{v^g}dsdy'
\end{equation}
is also  determined  by the  boundary data.

\begin{lemma}																	\label{lma:4.2}
Let  $f\in H_0^1(\Gamma_1\times[T_1,T_2])$.  For any  $s_0\in [T_1,T_2)$ there  exists  $u_0\in H_0^1(R_{2s_0})$
such  that
\begin{equation}							 									\label{eq:4.3}
(u_s^f,v')=(u_{0s},v')
\end{equation}
for any  $v'\in  H_0^1(Y_{3s_0})$.
Note  that  $R_{2s_0}=\{\tau=0,s_0-T_1\leq  s\leq  T_2-T_1,y'\in \Gamma_2\}$   (cf.  Fig. 4.1).
\end{lemma}
{\bf  Proof:}
Note  that
$Y_{1T_1}\cap\{s_0-T_1\leq  s\leq T_2-T_1\}\subset R_{2s_0}$.
Let  $w_1$  be such  that  $w_{1s}=0$  in  $R_{2s_0},\ w_1=u^f$  when  $s=s_0-T_1,y'\in \Gamma_2$.  Then
$u_0=u^f-w_1$    for 
$s\geq  s_0-T_1,\ u_0=0$  for  $s\leq  s_0-T_1$,  belongs  to  $H_0^1(R_{2s_0})$  and solves  (\ref{eq:4.3}).
\qed
\\
\begin{tikzpicture}
\draw[->](0,-0.5)-- (0,5);
\draw[->](-0.5,0)-- (8,0);
\draw[ultra thick] (1,4) --(6,4);
\draw[ultra thick] (2.5,0)--(4.5,0);
\draw (1,2)--(1,4);
\draw(1,2)--(6,2);
\draw(6,2)--(6,4);

\draw[ultra thick](1,4)..controls (1.5,3) and (1,2) .. (2.5,0);
\draw[ultra thick](6,4)..controls (5.5,3) and (6,2) .. (4.5,0);

\draw(3.5,-0.3) node {$\Gamma_j$};
\draw(3.5,4.3)  node {$\Gamma_{j+1}$};
\draw(4,1.5) node {$Y_{jT_1}$};
\draw((6.7,3.1) node {$R_{j+1,s_0}$};
\draw(-1.2,4) node {$s=T_2-T_1$};
\draw(0,5.5) node {$s$};
\draw(-1.2,2) node {$s=s_0-T_1$};
\draw(-1.2,0.1) node {$s=0$};

\filldraw[gray](0,0) circle (2pt)
                      (0,2) circle (2pt)
                      (0,4) circle (2pt);

\end{tikzpicture}
\\
\
{\bf Fig. 4.1.} 
\\  
$R_{j+1,s_0}$ is  
the rectangle $\{s_0-T_1\leq s\leq T_2-T_1,\tau=0,y'\in\Gamma_{j+1}\},\\ 
  Y_{jT_1}\cap\{s_0-T_1\leq s\leq T_2-T_1\}\subset R_{j+1,s_0}$.
\\
\\
\
If $v'=v^{g'}$  where  $g'\in  H_0^1(\Gamma_2\times[s_0,T_2])$  then  $(u_0,v^{g'})=(u^f,v^{g'})$  is determined  by the DN  operator.
Let  $g\in H_0^1(\Gamma_1\times[T_1,T_2])$.  We shall  show  that  still  $(u_0,v^g)$  is determined  by the 
DN operator.  The following theorem  holds.
\begin{theorem}																	\label{theo:4.3}
Let  $L_1^{(i)},i=1,2$,  be  two operators of the form  (\ref{eq:2.27}).  Let $f$  be in $H_0^1(\Gamma_1\times[T_1,T_2])$  
and let $u_0^{(i)}$
be the same  as in (\ref{eq:4.3})  
for  $i=1,2$.  Then
\begin{equation} 																	\label{eq:4.4}
(u_{0s}^{(1)},v_1^g)\big|_{Y_{2s_0}^{(1)}}=(u_{0s}                                                                                                                                                                                                                                                                                                                                                                                                                                                                                                                                                                                                                                                                                                                                                                                                                                                                                                                                                                                                                                                                                                                                                                                                                                                                                                                                                                                                                                                                                                                                                                                                                                                                                                                                                                                                                                                                                                                                                                                                                                                                                                                                                                                                                                                                                                                                                                                                                                                                                                                                                                                                                                                                                                                                                                                                                                                                                                                                                                                                                                                                                                                                                                                                                                                                                                                                                                                                                                                                                                                                                                                                                                                                                                                                                                                                                                                                                                                                                                                                                                                                                                                                                                                                                                                                                                                                                                                                                                                                                                                                                                                                                                                                                                                                                                                                                                                                                                                                                      ^{(2)},v_2^g)\big|_{Y_{2s_0}^{(2)}}
\end{equation}
for all $g\in H_0^1(\Gamma_2\times[T_1,T_2]).$

Here  $u_i^f,v_i^g$  are  the same  as in  (\ref{eq:3.1}), (\ref{eq:3.2})  for  $i=1,2$,  respectively.  Operators $L^{(i)}$  and,
consequently,  $L_1^{(i)}$  are formally  self-adjoint,  $(u_0^{(i)},v_i^g)_{Y_{2s_0}^{(i)}}$  is the $L_2$-inner product
over $Y_{2s_0}^{(i)}, i=1,2.$
\end{theorem}
To prove Theorem \ref{theo:4.3}  we will need  the Density Lemma \ref{lma:3.1}
and  the following lemma that uses the BLR condition:
\begin{lemma}																	\label{lma:4.4}
Let $L^{(1)}$  and  $L^{(2)}$  be  two operators in  $D\cap[t_0,T_2]$  having the same  DN  operator on 
$\Gamma_0\times[t_0,T_2]$.  Suppose  $L^{(1)}$  satisfies  the BLR  condition  on $[t_0,T_2]$.

Let  $L_1^{(i)},u_i^f,X_{2s_0}$  be the same as in (\ref{eq:3.1}),  $i=1,2,f\in H_0^1(\Gamma_{2s_0})$,
where $\Gamma_{2s_0}=\Gamma_2\times[s_0,T_2]$.
Then
\begin{equation} 																\label{eq:4.5}
\big\|u_2^f\big\|_{1,Y_{2s_0}^{(2)}}\leq C_2\big\|u_1^f\big\|_{1,Y_{2s_0}^{(1)}}.
\end{equation}
\end{lemma}
{\bf Proof of Lemma \ref{lma:4.4}} (cf.  Lemma 2.3  in [E3]):

Suppose  that  BLR  condition  (see [BLR])  is satisfied for $L^{(1)}$  on  $[t_0,T_{t_0}]$  and  
$t_0<T_1,\ T_2\geq T_{t_0}$.
The BLR condition implies  that  the  map  
$f\rw\big(u_1^f(x)\big|_{D_{T_2}},\frac{\partial u_1^f(x)}{\partial x_0}\big|_{D_{T_2}}\big)$  of
$H_+^1(\Gamma_0\times(t_0,T_2))$  to  $H^1(D_{T_2})\times L_2(D_{T_2})$  is onto,
where $D_{T_2}=D_0^{(1)}\times\{x_0=T_2\}$.  It follows  
from [H] (cf. also Lemma \ref{lma:3.6})  that
\begin{equation}															\label{eq:4.6}
\big\|u_1^f(x)\big\|_{1,D_{T_2}}^2+\big\|\frac{\partial u_1^f}{\partial x_0}\big\|_{0,D_{T_2}}
\leq C_0\big [f\big]_{1,\Gamma_0\times (t_0,T_2)}.
\end{equation}
By the closed  graph  theorem  we have
\begin{equation}															\label{eq:4.7}
\inf_{\mathcal F} \big[f'\big]_{1,\Gamma_0\times[t_0,T_2]}\leq  C_1\Big(\big\|u_1^f\big\|_{1,D_{T_2}}^2+
\big\|\frac{\partial u_1^f}{\partial x_0}\big\|_{0,D_{T_2}}^2\Big),
\end{equation}
where
$\mathcal F\subset  H_+^1(\Gamma_0\times(t_0,T_2))$  is  the set  of $f'$ such that 
\begin{equation}																\label{eq:4.8}
u_1^{f'}(x)\big|_{D_{T_2}}=u_1^f(x)\big|_{D_{T_2}},\ \ \ \frac{\partial u_1^{f'}(x)}{\partial x_0}\Big|_{D_{T_2}}
=\frac{\partial u_1^{f}(x)}{\partial x_0}\Big|_{D_{T_2}}.
\end{equation}
It follows  from  (\ref{eq:4.7})
that there exists  $f_0\in  \mathcal F,f_0=0$  for $x_0<t_0$  such  that
\begin{equation}     																\label{eq:4.9}
\big[f_0\big]_{1,\Gamma_0\times[t_0,T_2]}\leq  C\Big(\big\|u_1^f\big\|_{1,D_{T_2}}^2+
\big\|\frac{\partial u_1^f}{\partial x_0}\big\|_{0,D_{T_2}}^2\Big),
\end{equation}
Note that 
\begin{equation} 																\label{eq:4.10}
u_1^{f_0}(x)\big|_{D_{T_2}}=u_1^f(x)\big|_{D_{T_2}},\ \ \ \frac{\partial u_1^{f_0}(x)}{\partial x_0}\Big|_{D_{T_2}}
=\frac{\partial u_1^{f}(x)}{\partial x_0}\Big|_{D_{T_2}}.
\end{equation}
Let  $\Lambda^{(i)}$  be  the DN operator  corresponding to  $L^{(i)},i=1,2,$
and
let $u_{i0}^{f_0}$  be the solution  of  $L^{(i)}u_{i0}^{f_0}=0$  in   $D\times[t_0,T_2),\ u_{i0}^{f_0}=f_0$  on  
$\Gamma_0\times[t_0,T],\ i=1,2$.   Let  $L_1^{(i)}u_i^{f_0}=0$  be  the solutions in a neighborhood  $U_0$  obtained  from  
$L^{(i)}u_{i0}^{f_0}=0$  as in  \S 2  (cf.  (\ref{eq:2.38}))  for  $i=1,2.$

Consider the identity
\begin{equation}																\label{eq:4.11}
(L_1^{(i)}u_i^{f_0},v_i^g)\big|_{X_{js_0}^{(i)}}-(u_i^{f_0},L_1^{(i)}v_i^g)\big|_{X_{js_0}^{(i)}}=0,
\end{equation}
where $v^g$  is the same as in (\ref{eq:3.2}).
Since $\mbox{supp}\,v^g \subset  D(\Gamma_{js_0})$,  where $D(\Gamma_{js_0})$  is the  domain  of 
 influence   of $\Gamma_{js_0}$
for  $y_0\leq T_2$, 
we have  that  $v_i^g=0$  on $Z_{js_0}$.  Therefore  integrating   by  parts in (\ref{eq:4.11})  we get as in  (\ref{eq:3.7}):
$$
(u_{1s}^{f_0},v_1^g)\big|_{Y_{js_0}^{(1)}}-(u_1^{f_0},v_{1s}^g)\big|_{Y_{js_0}^{(1)}}=
-(\Lambda_1^{(1)}f_0,g)\big|_{\Gamma_{js_0}}+(f_0,\Lambda_1^{(1)}g)\big|_{\Gamma_{js_0}}.
$$
Analogously,  we have  for $L^{(2)}u_2^{f_0}=0, \ L^{(2)}v_2^g=0$:
$$
(u_{2s}^{f_0},v_2^g)\big|_{Y_{js_0}^{(2)}}-(u_2^{f_0},v_{2s}^g)\big|_{Y_{js_0}^{(2)}}=
-(\Lambda_1^{(2)}f_0,g)\big|_{\Gamma_{js_0}}+(f_0,\Lambda_1^{(2)}g)\big|_{\Gamma_{js_0}}.
$$
We have that  $\Lambda^{(1)}f_0=\Lambda^{(2)}f_0$   on
$\Gamma_0\times[t_0,T_2]$.   
Therefore  (\ref{eq:2.38})  implies  that  $\Lambda_1^{(1)}f_0=\Lambda_1^{(2)}f_0$  in  $\Gamma_{js_0}$.   Also  
$\Lambda_1^{(1)}g=\Lambda_1^{(2)}g$  in  $\Gamma_{js_0}$.
Integrating  by parts  we get
$$
-(u_i^{f_0},v_{is}^g)=(u_{is}^{f_0},v_i^g)-\int\limits_{R^{n-1}}\big(u_i^{f_0}\overline{v_i^g}\big|_{s=T_2-T_1}
-u_i^{f_0}v_i^g\big|_{s=0}\big)dy'.
$$
Note that $v_i^g\big|_{s=0}=0$  and  $u_i^{f_0}\overline{v_1^g}\big|_{s=T_2-T_1}=f_0(T_2,y')\overline{g(T_2,y')}$.
Therefore
\begin{equation}																\label{eq:4.12}
(u_{1s}^{f_0},v_1^g)=(u_{2s}^{f_0},v_2^g)
\end{equation}
for all  $g\in  H_0^1(\Gamma_{3s_0})$.

Let $\Gamma_2,\Gamma_3,\Gamma_4$  be  the same  as  in Lemma \ref{lma:3.5}.   It was  proven  there that
\begin{equation}														\label{eq:4.13}
\big\|u_1^f\big\|_{1,\Gamma_4}^2+\big\|u_{1y_0}^f\big\|_{0,\Gamma_4}^2
\leq C\Big(\big\|u_1^f\big\|_{1,\Gamma_2}^2+\big\|u_1^f\big\|_{1,\Gamma_3}^2\Big),
\end{equation}
\begin{equation}														\label{eq:4.14}
\big\|u_1^f\big\|_{1,\Gamma_2}^2+\big\|u_{1}^f\big\|_{1,\Gamma_3}^2
\leq C\Big(\big\|u_1^f\big\|_{1,\Gamma_4}^2+\big\|u_{1,y_0}^f\big\|_{0,\Gamma_4}^2\Big).
\end{equation}
It follows from $u_1^f\big|_{\Gamma_4}=u_1^{f_0}\big|_{\Gamma_4},\  u_{1y_0}^f\big|_{\Gamma_4}=
u_{1y_0}^{f_0}\big|_{\Gamma_4}$  that
\begin{equation}														\label{eq:4.15}
u_1^f\big|_{\Gamma_2}=u_1^{f_0}\big|_{\Gamma_2}
\end{equation}
by the domain of dependence argument.  Comparing (\ref{eq:4.12})   with 
$(u_{1s}^{f},v_1^g)=(u_{2s}^{f},v_2^g)$  and  taking into account 
(\ref{eq:4.15})  we get
\begin{equation}														\label{eq:4.16}
(u_{2s}^{f_0},v_2^g)\big|_{Y_{2s_0}^{(2)}}=(u_{2s}^{f},v_2^g)\big|_{Y_{2s_0}^{(2)}},\ \ \ 
\forall g\in H_0^1(\Gamma_3\times(s_0,T_2)).
\end{equation}
By Lemma \ref{lma:3.1}
$\{v_2^g\}$  are dense in $H_0^1(R_{3s_0}^{(2)})$.  Since  
$Y_{2s_0}^{(2)}\subset R_{3s_0}^{(2)}$  we  get  that $\{v_2^g\}$  are  dense in $H_0^1(Y_{2s_0}^{(2)})$
and  therefore $u_{2s}^{f_0}=u_{2s}^f$  in  $Y_{2s_0}^{(2)}$.   Since  $u_2^f\big|_{s=T_2-T_1}=f(T_2,y')=u_1^f(T_2,y',0),
\ u_2^{f_0}\big|_{s=T_2-T_1}=f_0(T_2,y')=u_1^{f_0}(T_2,y',0)$  and since
$u_1^{f_0}(T_2,y',0)=u_1^{f}(T_2,y',0)$  we  get  that $u_2^f\big|_{s=T_2-T_1}=u_2^{f_0}\big|_{s=T_2-T_1}.$  Thus
\begin{equation}  															\label{eq:4.17}
u_2^{f_0}=u_2^f\ \ \mbox{on}\ \ \ Y_{2s_0}^{(2)}.
\end{equation}
It follows from  (\ref{eq:4.13})  that
\begin{equation}																\label{eq:4.18}
\big\|u_1^f\big\|_{1,\Gamma_4}^2+\big\|u_{1y_0}^f\big\|_{0,\Gamma_4}^2\leq  C\big\|u_1^f\big\|_{1,\Gamma_2}^2,
\end{equation}
since that $u_1^f=0$  on $\Gamma_3$  by the  domain  of dependence argument.  

Since  $Y_{2s_0}^{(2)}$  belongs  to the domain  of dependence  of  $D_{T_2}$  we get,  similarly to
  (\ref{eq:4.14}),   that
\begin{equation} 																	\label{eq:4.19}
\big\|u_2^{f_0}\big\|_{1,Y_{2s_0}^{(2)}}^2
\leq C_1\big(\big\|u_2^{f_0}\big\|_{1,D_{T_2}^{(2)}}^2+\big\|u_{2y_0}^{f_0}\big\|_{0,D_{T_2}^{(2)}}^2\big),
\end{equation}
where $D_{T_2}^{(2)}=D^{(2)}\cap\{y_0=T_2\}$.

We also have  (cf.  Lemma \ref{lma:3.6})
\begin{equation}															\label{eq:4.20}
\big\|u_2^{f_0}\big\|_{1,D_{T_2}^{(2)}}^2 +  \Big\|\frac{\partial u_2^{f_0}}{\partial x_0}\Big\|_{0,D_{T_2}^{(2)}}^2
\leq C\big[f_0\big]_{1,\Gamma_0\times[t_0,T_2]}^2.
\end{equation}
Combining (\ref{eq:4.18}), (\ref{eq:4.9})   with (\ref{eq:4.19}), (\ref{eq:4.20})  and taking  into account   
(\ref{eq:4.17}),  we get
\begin{equation} 															\label{eq:4.21}
\big\|u_2^{f}\big\|_{1,Y_{2s_0}^{(2)}}\leq  C\big\|u_1^{f}\big\|_{1,Y_{2s_0}^{(1)}}.
\end{equation}
Now we shall  prove  Theorem \ref{theo:4.3}.

{\bf Proof of Theorem \ref{theo:4.3}}
Since  $u_0^{(1)}\in H_0^1(R_{2s_0}^{(1)})$  we get,  using  the Density Lemma  \ref{lma:3.1},  that there exists  
$u_1^{f_n},f_n\in  H_0^1(\Gamma_2\times[s_0,T_2])$  such that $\big\|u_0^{(1)}-u_1^{f_n}\big\|_{1,Y_{2s_0}^{(1)}}\rw 0$.
By Lemma \ref{lma:4.4}
$\{u_2^{f_n}\}$  also  converges  in $H_0^1(Y_{2s_0}^{(2)})$  to some function  $w\in  H_0^1(Y_{2s_0}^{(2)})$.  
Passing to the limit  in 
\begin{equation}															\label{eq:4.22}
\big(u_{1s}^{f_n},v_1^g\big)=\big(u_{2s}^{f_n},v_2^g\big),
\end{equation}
we get
\begin{equation}															\label{eq:4.23}
\big(u_{0s}^{(1)},v_1^g\big)=\big(w_{s},v_2^g\big)\ \ \ \mbox{for any}\ \ \ g\in H_0^1(\Gamma_{3T_1}),
\end{equation}
where  $\Gamma_{3T_1}=\Gamma_3\times[T_1,T_2]$
Note that  (\ref{eq:4.22})  and  therefore  (\ref{eq:4.23})  hold also  for any $g'\in  H_0^1(\Gamma_{3s_0}),$
i.e.
\begin{equation}															\label{eq:4.24}
\big(u_{0s}^{(1)},v_1^{g'}\big)=\big(w_{s},v_2^{g'}\big).
\end{equation}
For such  $g'$  the equality (\ref{eq:4.3}) holds,  i.e. 
\begin{equation}															\label{eq:4.25}
\big(u_{0s}^{(1)},v_1^{g'}\big)=\big(u_{0s}^{(2)},v_2^{g'}\big).
\end{equation}
Comparing (\ref{eq:4.24})  and  (\ref{eq:4.25})  we get   
\begin{equation}															\label{eq:4.26}
\big(u_{0s}^{(2)},v_2^{g'}\big)=\big(w_{s},v_2^{g'}\big),
\end{equation}
Since
$v_2^{g'}\in  H_0^1(Y_{3s_0}^{(2)})$  are dense  in  $H_0^1(R_{3s_0}^{(2)})$  and  $w\in  H_0^1(Y_{2s_0}^{(2)})
\subset   H_0^1(R_{3s_0}^{(2)})$,  we have  that  $u_{0s}^{(2)}=w_s$.   Since  $u_0^{(2)}$     and 
$w$   are zero on  $\partial Y_{3s_0}^{(2)}\setminus\{y_n=0\}$  we get that 
\begin{equation}															\label{eq:4.27}
u_0^{(2)}=w\ \ \mbox{in}\ \ \  Y_{2s_0}^{(2)}.
\end{equation} 
Therefore (\ref{eq:4.23}) and  (\ref{eq:4.27})  gives
\begin{equation}															\label{eq:4.28}
\big(u_{0s}^{(1)},v_1^g)=(u_{0s}^{(2)},v_2^g\big)
\end{equation}
for all  $g\in  H_0^1(\Gamma_{3T_1})$,  i.e.      (\ref{eq:4.4})  holds.
\qed
\\

The following formula will be the main tool in solving the  inverse problem.
\begin{theorem}																\label{theo:4.5}
For any  $T_1\leq s_0\leq T_2$  the integral  
\begin{equation}															\label{eq:4.29}
\int\limits_{Y_{jT_1}\cap\{0\leq s\leq s_0-T_1\}}\frac{\partial u^f}{\partial s}\overline{v^g}dsdy',\ \ 
\forall f\in H_0^1(\Gamma_{jT_1}),\ \forall  g\in  H_0^1(\Gamma_{jT_1}),\  j=1,2,
\end{equation}
is determined  by the DN operator  on $\Gamma_{jT_1}=\Gamma_j\times[T_1,T_2]$.
\end{theorem}
{\bf Proof:}
Since
$u_0^{(i)} =\frac{\partial  u^f}{\partial s}$  for  $s\geq  s_0-T_1, \ u_{0s}=0$  for  $s\leq s_0-T_1$,  formula 
(\ref{eq:4.28})  gives that 
$\int_{Y_{jT_1}\cap\{s> s_0-T_1\}}\frac{\partial u^f}{\partial s}\overline{v^g}dsdy'$ is determined by  the DN
operator  on  $\Gamma_{jT_1}$.  The integral  (\ref{eq:4.29})  is the difference  $(u_s^f,v^g)-(u_{0s},v^g)$
thus  (\ref{eq:4.29}) is determined  by DN operator.

{\bf Remark 4.1}.
 When the coefficients  of  $L_1^{(i)}, i=1,2$,  do not depend on $y_0$,  we can  obtain the estimate  (\ref{eq:4.5})
without assuming the  BLR  condition.  In this case  we can derive,  in addition  to (\ref{eq:3.3}),  another Green's formula  
(cf. (3.18)  in  [E3]):

Consider the identity
\begin{equation}															\label{eq:4.30}
0=(L_1u,v_{y_0})+(u_{y_0},L_1 v).
\end{equation}
Integrating  by parts  as in [E3]  and using  that  $\frac{\partial}{\partial  y_0}$  and $L_1$ are commute,  we get  (cf.  (3.20)  in
[E3])
\begin{equation}                               											\label{eq:4.31}
\tilde Q(u,v)=-\tilde\Lambda_0(f,g),
\end{equation}
where
\begin{align}															\label{eq:4.32}
&\tilde Q(u,v)
=\frac{1}{2}\int\limits_{Y_{20}}
\Bigg[ 
2(u_s+iA_-'u)\overline{v_s}+2u_s\overline{(v_s+iA_-'v)}
\\
\nonumber
&-2\sum_{j=1}^{n-1}\Bigg(
g_0^{0j}\Big(\frac{\partial u}{\partial s}+iA_-'u\Big)
\overline{\Big(\frac{\partial v}{\partial y_j} -iA_j'v\Big)}
+g_0^{0j}\Big(\frac{\partial u}{\partial y_j}-iA_j'u\Big)\overline{\Big( \frac{\partial v}{\partial s}+iA_-'v\Big)}\Bigg)
\\
\nonumber
&-\sum_{j,k=1}^{n-1}g_0^{jk}\Big(\frac{\partial}{\partial y_j}-iA_j'\Big)u \overline{ 
\Big(\frac{\partial}{\partial y_k}-iA_k'
\Big)v
}
+V_1u\overline v
\Bigg]dsdy'
\end{align}
and
\begin{equation}																\label{eq:4.33}
\Lambda_0(f,g)=\int\limits_{\Gamma^{(2)}\times[0,T]}
\big(\Lambda_1 f \ \overline{g_{y_0}}+f_{y_0}\overline{\Lambda_1g}\big)dy'dy_0.
\end{equation}
As in [E3]  (cf. (3.23)  in  [E3]) we have  that  $\tilde Q(u,u)$  is a positive  definite form  when  $T_2-T_1$  
is small  and
\begin{equation}																\label{eq:4.34}
C_2\big\|u\big\|_{1,Y_{2s_0}}^2\leq \tilde Q(u,u)\leq C_1\big\|u\big\|_{1,Y_{2s_0}}^2.
\end{equation}
Let $u_i^f, i=1,2,$  be such that  
$L_1^{(i)}u_i^f=0$  in $X_{2s_0}^{(i)},\ u_i^f\big|_{y_n=0}=f,\ u_i^f=0$  for $y_0<T_1,\  \mbox{supp}\,f$  is contained in  
$\Gamma_2\times(T_1,T_2]$.  We assume that  $\Lambda_1^{(1)}=\Lambda_1^{(2)}$  on  
$\Gamma_{2T_1}=\Gamma_2\times(T_1,T_2)$.  
It follows from (\ref{eq:4.13}),  (\ref{eq:4.31}),  (\ref{eq:4.33})   that 
\begin{equation}																\label{eq:4.35}
\tilde Q_1(u_1^f,u_1^f)=\tilde Q_2(u_2^f,u_2^f),
\end{equation}
where $\tilde Q_i$  corresponds  to $L_1^{(i)},i=1,2$.  Thus,  (\ref{eq:4.34})  implies  that
\begin{equation}       																\label{eq:4.36}
C_1\big\|u_1^f\big\|_{1,Y_{2s0}^{(1)}}\leq \big\|u_2^f\big\|_{1,Y_{2s0}^{(2)}}\leq C_2\big\|u_1^f\big\|_{1,Y_{2s0}^{(1)}},
\end{equation}
i.e.  the estimate (\ref{eq:4.5})  is proven.

\section{The geometric optics construction}
\label{section 5}
\init
It follows from Theorem \ref{theo:4.5} that the DN operator allows to determine\linebreak
$
\int_{Y_{2s_0}\cap\{s\leq s_0 - T_1\}}u_{s}^f\overline{v^g}dsdy'$
for all  $f\in H_o^1(\Gamma_{jT_1}),g\in H_0^1(\Gamma_{jT_1}),j=1,2,$
  i.e.  if $u_i^f,v_i^g$  satisfy (\ref{eq:3.1}),  (\ref{eq:3.2}),  $i=1,2,$
then 
$$
\int\limits_{Y_{2s_0}^{(1)}\cap\{s\leq s_0 - T_1\}}u_{1s}^f\overline{v_1^g}dsdy'=
\int\limits_{Y_{2s_0}^{(2)}\cap\{s\leq s_0 - T_1\}}u_{2s}^f\overline{v_2^g}dsdy'
$$
Let $u_i$  be  the solution of $L^{(i)}u_i=0$  such that
\begin{equation}																\label{eq:5.1}
u_i=u_N^{(i)}+u_i^{(N+1)},\ \ \ u_N^{(i)}=\sum_{p=0}^N\frac{a_p^{(i)}(s,\tau,y')}{(ik)^p}e^{ik(s-s_0')},\ \ s_0'=s_0-T_1,
\end{equation}
$u_i^{(N+1)}$  will be chosen below,  $k$  is a large parameter.
We have the following equations  for $a_p^{(i)},0\leq p\leq N$  (see  [E4]  and [E5],  \S64,  for more  details on  
the construction  of geometric optics type solutions):
 \begin{multline}																\label{eq:5.2}
 -4i\Big(\frac{\partial}{\partial \tau}+iA_-^{(i)}\Big)a_0^{(i)}
 +2i\sum_{j=1}^{n-1}g_{i0}^{oj}\Big(\frac{\partial}{\partial y_j}-iA_j^{(i)}\Big)a_0^{(i)}
 \\
 +2i\sum_{j=1}^{n-1}\Big(\frac{\partial}{\partial y_j}-iA_j^{(i)}\big)\big(g_{i0}^{0j}a_0^{(i)}\big)
 =0,
 \end{multline}
\begin{multline}																\label{eq:5.3}
 -4i\Big(\frac{\partial}{\partial \tau}+iA_-^{(i)}\Big)a_p^{(i)}
 +2i\sum_{j=1}^{n-1}g_{i0}^{oj}\Big(\frac{\partial}{\partial y_j}-iA_j^{(i)}\Big)a_p^{(i)}
 \\
 +2i\sum_{j=1}^{n-1}\Big(\frac{\partial}{\partial y_j}-iA_j^{(i)}\big)\big(g_{i0}^{0j}a_0^{(i)}\big)
 =
 -L_1^{(i)}a_{p-1}^{(i)},\ p\geq 1,
 \end{multline}
with the initial conditions
\begin{equation}																\label{eq:5.4}
a_0^{(i)}(s,\tau,y')\big|_{\tau=\tau_0}=\chi_1(s)\chi_2(y'),\ \tau_0=T_2-T_1-s,
\end{equation}
\begin{equation}																\label{eq:5.5}
a_p^{(i)}(s,\tau,y')\big|_{\tau=\tau_0}=0,\ \ p\geq 1,
\end{equation}
where 
$ \chi_1(s)=0$ for 
$|s-s_0'|>2\delta,\  \chi_1(s)=1$ for $|s-s_0'|\leq \delta$,  $\chi_2(y')\in C_0^\infty(\Gamma_2),\ 
\chi_2(y')\neq 0$  when  $|y'-y_0'|<\delta,\ y_0'\in \Gamma_2$  is  arbitrary,
$g_{i0}^{j0}$ corresponds  to  $L_1^{(i)},i=1,2$.  Note  that  $y_n=\frac{T_2-T_1-s-\tau}{2}=0$  when  $\tau=\tau_0$.

Let  $u_i^{(N+1)}$  be such  that 
\begin{equation}																\label{eq:5.6}
L_1^{(i)}u_i^{(N+1)}=-\frac{1}{(ik)^N}\big(L_1^{(i)}a_N^{(i)}\big)e^{ik(s-s_0')},\ \ y_n>0,\ y_0<T_2,
\end{equation}
$u_i^{(N+1)}=u_{iy_0}^{(N+1)}=0$  when  $y_0=T_1, y_n>0, i=1,2, \ u_i^{(N+1)}\big|_{y_n=0}=0,\ y_n\leq T_2$.
Such  $u_i^{(N+1)}$  exists  (cf. [H])  and $L_1^{(i)}(u_N^{(i)}+u_i^{(N+1})=0$.

Since $\mbox{supp}\,u_N^{(i)}$ is contained  in a small  neighborhood of the line 
$\{s=s_0-T_1,y'=y_0'\}$,  we have that 
$\mbox{supp}\,(u_N^{(i)}+u_i^{N+1})\subset D_+(\Gamma_2\times[T_1,T_2])$  when  $s_0-T_1>0$.
Here ,  as in \S 1,  $D_+(\Gamma_2\times [T_1,T_2])$  is  the forward  domain  of influence  of $\Gamma_2\times [T_1,T_2]).$

Let $\beta^{(i)}(s,\tau,\hat y')=(\beta_1^{(i)},\beta_2^{(i)},..,\beta_{n-1}^{(i)})$  be  the solution  of  the  system (cf. [E4]) 
\begin{equation}																\label{eq:5.7}
\frac{\partial \beta_j^{(i)}(s,\tau,\hat y')}{\partial\tau}=-g_{i0}^{0j}(s,\tau,\beta^{(i)}(s,\tau,\hat y')),\ 1\leq j\leq n-1,\ y_n>0,
\end{equation}
\begin{equation}																\label{eq:5.8}
\beta^{(i)}(s,\tau,\hat y')\big|_{\tau=\tau_0}=\hat y_i',\ \ i=1,2,\ \tau_0=T_2-T_1-s,
\end{equation}
where $\hat y'=(\hat y_1,...,\hat y_{n-1})\in \Gamma_2$,  $s$  is a parameter in  (\ref{eq:5.7}).

Let
\begin{equation}																\label{eq:5.9}
\hat s=s,\hat\tau=\tau,\ \hat y'=\alpha^{(i)}(s,\tau,y'),\ \alpha^{(i)}=(\alpha_1^{(i)},...,\alpha_{n-1}^{(i)})
\end{equation}
be  the inverse  to the map
\begin{equation}																\label{eq:5.10}
s=\hat s,\ \tau=\hat\tau, \ y'=\beta^{(i)}(\hat s,\hat \tau,\hat y'),
\end{equation}
i.e. 
\begin{equation}																\label{eq:5.11}
\alpha_j^{(i)}(s,\tau,\beta^{(i)}(s,\tau,\hat y'))=\hat y_j,\ \ 1\leq j\leq n-1.
\end{equation}
Note that  $\alpha_j^{(i)}(s,\tau,y'),1\leq j\leq n-1,$  satisfy  the equation
\begin{equation}																\label{eq:5.12}
\frac{\partial \alpha_j^{(i)}(s,\tau,y')}{\partial\tau}-\sum_{k=0}^{n-1}g_{i0}^{k0}(s,\tau,y')\frac{\partial\alpha_j^{(i)}}{\partial y_k}=0,
\ \ \alpha_j^{(i)}\big|_{\tau=\tau_0}=y_j,\ 1\leq j\leq n-1.
\end{equation}
Let $\hat a_0^{(i)}(s,\tau,\hat y')=a_0^{(i)}(s,\tau,y')$,  where
$y'=\beta^{(i)}(s,\tau,\hat y')$.  Then using  (\ref{eq:5.7})  and  (\ref{eq:5.2})  we get
\begin{multline}                                                                                                               \label{eq:5.13}
\frac{\partial\hat a_0^{(i)}}{\partial\tau}=\frac{\partial a_0^{(i)}}{\partial\tau}+
\sum_{j=1}^{n-1}\frac{\partial a_0^{(i)}}{\partial y_j}\frac{\partial\beta_j^{(i)}}{\partial\tau}
\\
=\frac{\partial a_0^{(i)}}{\partial t}-\sum_{j=1}^n g_{i0}^{0j}\frac{\partial a_0^{(i)}}{\partial y_j}
=\hat B^{(i)}(s,\tau,\hat y')\hat a_0^{(i)}(s,\tau,\hat y'),
\end{multline}
where
$B^{(i)}(s,\tau,y')=-iA_-'-i\sum_{j=1}^{n-1}g_{i0}^{0j} A_j'+\frac{1}{2}\sum_{j=1}^{n-1}\frac{\partial g_{0i}^{0j}}{\partial y_j},\ 
\hat B^{(i)}(s,\tau,\hat y')=B^{(i)}(s,\tau,\beta^{(i)}(s,\tau,\hat y')),\ \hat a_0^{(i)}(s,\tau,\hat y')\big|_{\tau=\tau_0}
=\chi_1(s)\chi_2(\hat y')$.

Therefore
\begin{equation}																	\label{eq:5.14}
a_0^{(i)}(s,\tau,y')=\chi_1(s)\chi_2(\alpha^{(i)}(s,\tau,y'))e^{b^{(i)}(s,\tau,\alpha^{(i)})}
\end{equation} 
where $ b^{(i)}(s,\tau,\hat y')=\int_{\tau_0}^\tau\hat B^{(i)}(s,\hat\tau,\hat y')d\hat\tau$.  Substituting  $u=u_N^{(i)}+u_i^{(N+1)}$  into 
(\ref{eq:4.29})  instead  of  $u_i^f$,  integrating  by parts  in $s$  and taking  the  limit  when  $k\rw \infty$,  we get 
\begin{multline}																	      \label{eq:5.15}
\int\limits_{\R^{n-1}}e^{b^{(1)}(s_0',0,\alpha^{(1)})}  \chi_2(\alpha^{(1)}(s_0',0,y'))\overline{v_1^g(s_0',0,y')}dy'
\\
=
 \int\limits_{\R^{n-1}}e^{b^{(2)}(s_0',0,\alpha^{(2)})}  \chi_2(\alpha^{(2)}(s_0',0,y'))\overline{v_2^g(s_0',0,y')}dy'.
 \end{multline} 
Note that  $\tau=0$  on  $Y_{2T_1}^{(i)},i=1,2.$  In (\ref{eq:5.15}) $s_0\in  (T_1,T_2]$  is arbitrary,  $s_0'=s_0-T_1$.

Denote by  $Y_{2T_1}^{(i)}(\tau')$  the intersection  of the plane  $\tau=\tau'$  with  $X_{2T_1}^{(i)}$.  Let  
$R_{2T_1}^{(i)}(\tau')\subset Y_{2T_1}^{(i)}(\tau')$  be  the rectangle  $\{\tau=\tau',0\leq s\leq T_2-T_1-\tau',y'\in\Gamma_2\}$.
Note  that $Y_{2T_1}^{(i)}(0)=Y_{2T_1}^{(i)}$  and  $R_{2T_1}^{(i)}(0)=R_{2T_1}^{(i)}$.

Repeating  the proof  of Theorem \ref{theo:4.5} with  $Y_{2T_1}^{(2)},R_{2T_1}^{(2)}$  replaced  by 
$Y_{2T_1}^{(2)}(\tau'),\linebreak R_{2T_1}^{(2)}(\tau'),\ 0\leq \tau'\leq T_2-T_1$,  we get  again,  using
 the geometric optics construction  
(\ref{eq:5.1}),  that  (\ref{eq:5.15})   holds   for any  $(s,\tau)\in  \Sigma$,  where 
$\Sigma=\{(s,\tau),s\geq 0,\tau\geq 0, s+\tau\leq T_2-T_1\}$.  Thus, we have
 \begin{multline}																	\label{eq:5.16}
\int\limits_{\R^{n-1}}e^{b^{(1)}}  \chi_2(\alpha^{(1)}(s,\tau,y'))\overline{v_1^g(s,\tau,y')}dy'
\\
=
 \int\limits_{\R^{n-1}}e^{b^{(2)}}  \chi_2(\alpha^{(2)}(s,\tau,y'))\overline{v_2^g(s,\tau,y')}dy'.
 \end{multline} 
for  $(s,\tau,y')\in X_{2T_1}^{(i)}$.  

Let  $\beta^{(i)}(\Sigma\times\overline\Gamma_2)$  be  the image  of
$\Sigma\times\overline\Gamma_2$
under  the map  (\ref{eq:5.10}).   Note that the support  of  geometric optics solution  $u_N^{(i)}+u_i^{(N+1)}$  is
contained  in $D(\Gamma_2\times[T_1,T_2])$.  Also  we have that 
the curve  $y'=\beta^{(i)}(s,\hat\tau,\hat y')$  for  $\tau_0\leq\hat\tau\leq\tau$,  is contained in  $X_{2T_1}^{(i)}$.  
Therefore $\beta^{(i)}(\Sigma\times\overline\Gamma_2)\subset X_{2T_1}^{(i)}$.  Denote  by
$X_{\Gamma_2}^{(i)}$  the intersection  of  $\beta^{(i)}(\Sigma\times\overline\Gamma_2)$   with  
$\Sigma\times\overline\Gamma_2$.  Note that  $\Sigma\times\overline\Gamma_2=
\bigcup_{0\leq \tau'\leq T_2-T_1} R_{2T_1}^{(i)}(\tau')$.

Finally,  denote by  $\tilde X_{\Gamma_2}^{(i)}$  the image  of  
$X_{\Gamma_2}^{(i)}$ under the inverse  map (\ref{eq:5.9}).  
 Note that  $\tilde X_{\Gamma_2}^{(i)}\subset\Sigma\times\overline\Gamma_2$.

Making the change  of variables  (\ref{eq:5.10})  in (\ref{eq:5.16})  we get
\begin{multline}														\label{eq:5.17}
\int\limits_{\Gamma_2}e^{b^{(1)}(s,\tau,\hat y')}  \chi_1(\hat y')\overline{v_1^g(s,\tau,\beta^{(1)}(s,\tau,\hat y'))}
J_1(s,\tau,\hat y')d\hat y'
\\
=
\int\limits_{\Gamma_2}e^{b^{(2)}(s,\tau,\hat y')}  \chi_2(\hat y')\overline{v_2^g(s,\tau,\beta^{(2)}(s,\tau,\hat y'))}
J_2(s,\tau,\hat y')d\hat y',
 \end{multline} 
where  $J_i$  is the Jacobian of  the map  (\ref{eq:5.10}),  $(s,\tau,\hat y')\in  \Sigma\times \Gamma_2$.

Let  $b^{(i)}=b_1^{(i)}+ib_2^{(i)}$,  where  $b_1^{(i)}, b_2^{(i)}$  are real.

Since $\chi_2(y')\in C_0^\infty(\Gamma_2)$  is  arbitrary,  we have 
\begin{equation}																\label{eq:5.18}
e^{b_1^{(1)}-ib_2^{(1)} }v_1^g(s,\tau,\beta^{(1)}) J_1=e^{b_1^{(2)}-ib_2^{(2)}}v_2^g(s,\tau,\beta^{(2)}) J_2.
\end{equation}
Let 
\begin{align}																	\label{eq:5.19}
&w_i^g(s,\tau,\hat y')=v_i^g(s,\tau,\beta^{(i)}(s,\tau,\hat y')),  \ \ \ \hat y'\in \Gamma_2,
\\
\nonumber																			
 &\tilde w_i^g(s,\tau,\hat y')=w_i^g(s,\tau,\hat y')e^{-b^{(i)}(s,\tau,\hat y')}.
 \end{align}
Our  strategy  will  be to show  that  $w_1^g(s,\tau,\hat y')=w_2^g(s,\tau,\hat y')$  in $\tilde X_{\Gamma_2}^{(1)}$  and then
to show  that the equations 
$\tilde L_1^{(1)}w_1^g=0$  and  $\tilde L_1^{(2)}w_1^g=0$  have  the same coefficients  in 
$\tilde X_{\Gamma_2}^{(1)}$.  Here
$\tilde L_1^{(i)}$  is obtained  from  $L_1^{(i)}$  by the change of variables(\ref{eq:5.10}),  $i=1,2$.

We shall show first  that  $e^{2b_1^{(1)}}J_1(s,\tau,\hat y')=e^{2b_1^{(2)}}J_2(s,\tau,\hat y')$.
Consider the geometric optics   solutions $v_{i,k}^g$  of  the form  (\ref{eq:5.1}),  where 
$g=\chi_1(s)\chi_3(y'),\chi_3(y')\in C_0^\infty(\Gamma_2)$  is  arbitrary.  Substituting
$v_{i,k}^g$  into  (\ref{eq:5.16}),  integrating by parts  and passing  to the limit when  $k\rw\infty$,
we get 
\begin{multline}																\label{eq:5.20}
\int\limits_{\R^{n-1}}e^{2b_1^{(1)}}\chi_2(\alpha^{(1)}(s_0',\tau,y'))\overline{\chi_3(\alpha^{(1)}(s_0',\tau,y'))}dy'
\\
=
\int\limits_{\R^{n-1}}e^{2b_1^{(2)}}\chi_2((\alpha^{(2)}(s_0',\tau,y'))\overline{\chi_3(\alpha^{(2)}(s_0',\tau,y'))}dy',
\end{multline}
where  $s_0'=s_0-T_1$.

Note  that $e^{b^{(i)}}e^{\overline{b^{(i)}}}=e^{2b_1^{(i)}}$. 

Making the change of variables $y'=\beta^{(i)}(s_0,\tau,\hat y')$  and using  that $\chi_2$  and $\chi_3$  are arbitrary we get
\begin{equation}   															\label{eq:5.21}
e^{2b_1^{(1)}}J_1(s_0',\tau,\hat y')=e^{2b_1^{(2)}}J_2(s_0',\tau,\hat y').
\end{equation}
Therefore,  (\ref{eq:5.18})  and  (\ref{eq:5.21})
imply
\begin{align}																	\label{eq:5.22}
&e^{-b^{(1)}(s,\tau,\hat y')}v_1^g(s,\tau,\beta^{(1)}(s,\tau,\hat y'))
=
e^{-b^{(2)}(s,\tau,\hat y')}v_2^g(s,\tau,\beta^{(2)}(s,\tau,\hat y'))
\ \ \ \mbox{in}\ \ \Sigma\times\Gamma_2,
\\
\nonumber
&\mbox{i.e.}\ \ \ 
 \tilde w_1^g(s,\tau,\hat y')=\tilde w_2^g(s,\tau,\hat y').
\ \ \ \ \ \ \ \ \ \ \ \ \ \ \ \ \ \ \ \ \ \ \ \ \ \ \ \ \ \ \ \ \ \ \ \ \ \ \ \ \ \ \ \ \ \ \ \ \ \ \ \ \ \ \ \ \ 
\qed
\end{align}

As in  (\ref{eq:4.12})  the integration  by parts  gives 
$$
\int\limits_{Y_{3T_1}}(u_s^f\overline{v^g}-u^f\overline{v_s^g})dsdy'
=-2\int\limits_{Y_{3T_1}}u^f\overline{v_s^g}dsdy'
+\int\limits_{\partial Y_{3T_1}\cap\{y_n=0\}}u^f\big|_{y_n=0}\overline{v^g}\big|_{y_n=0}dy'.
$$
Therefore  $\int\limits_{Y_{3T_1}}u^f\overline{v_s^g}dsdy'$  is determined  by the boundary data
since  $u^f\big|_{y_n=0}=f(T_2,y'),\overline{v^g}\big|_{y_n=0}=   {\overline g(T_2,y')}$,  i.e.  the roles  of  $u^f$  and  $v^g$
are reversed  in  comparison  with  (\ref{eq:4.12}).  Therefore  we get,  as in (\ref{eq:4.28}),
\begin{equation}																\label{eq:5.23}
\int\limits_{Y_{2s_0}^{(1)}\cap\{s\leq s_0'\}}u_1^f\overline{v_{1s}^g}dsdy'  
=\int\limits_{Y_{2s_0}^{(2)}\cap\{s\leq s_0'\}}u_2^f\overline{v_{2s}^g}dsdy'.
\end{equation}																			 
Substituting  in  (\ref{eq:5.23})  the geometric  optics solution  (\ref{eq:5.1}),  integrating  by parts  in $s$,
multiplying by  $ik$  and,  finally,  taking the limit  when  $k\rw\infty$,    we get  (\ref{eq:5.16})  with  $v_i^g$  replaced  by
$v_{is}^g$.  Note that we assumed that  $v_i^g\in  H_0^2(Y_{2T_1})$  when  integrating  by parts  in  (\ref{eq:5.23}).
This  can be  achieved  by requiring that $g\in  H_0^2(\Gamma_{2T_1})$  
and using  the regularity  results  for hyperbolic  initial-boundary value problems
(cf.  [H], [E6]).
Therefore  we  get  (\ref{eq:5.18}),   with  $v_i^g$  replaced
 by $v_{is}^g$:
 \begin{equation}																
 \nonumber
 e^{b_1^{(1)}  -ib_2^{(1)}      } v_{1s}^g(s,\tau,\beta^{(1)}(s,\tau,\hat y'))J_1=
e^{b_1^{(2)} -ib_2^{(2)}      } v_{2s}^g(s,\tau,\beta^{(2)}(s,\tau,\hat y'))J_2.
\end{equation}
Using  (\ref{eq:5.21})  we get
\begin{equation}																\label{eq:5.24}
e^{-b^{(1)}}v_{1s}^g(s,\tau,\beta^{(1)})=e^{-b^{(2)}}v_{2s}^g(s,\tau,\beta^{(2)}).
\end{equation}
We shall  need  the following  lemma:
\begin{lemma}																	\label{lma:5.1}
The equalities  
\begin{equation}																\label{eq:5.25}
\alpha_{js}^{(1)}(s,\tau,\beta^{(1)}(s,\tau,\hat y'))=\alpha_{js}^{(2)}(s,\tau,\beta^{(2)}(s,\tau,\hat y')), \ \ 1\leq j\leq n-1,
\end{equation}
\begin{equation}																\label{eq:5.26}
b^{(1)}(s,\tau,\hat y')= b^{(2)}(s,\tau,\hat y')
\end{equation}
hold  on $\tilde X_{\Gamma_2}^{(1)}$.
 \end{lemma}
{\bf Proof:}  Making the change
of variables $\hat y'=\alpha^{(i)}(s,\tau,y')$  in  (\ref{eq:5.19}),  we  get
\begin{equation}																\label{eq:5.27}
e^{-b^{(i)}(s,\tau,\alpha^{(i)}(s,\tau,y'))}v_i^g(s,\tau,y')=\tilde w_i^g(s,\tau,\alpha^{(i)}(s,\tau,y')).
\end{equation}
Differentiating in $s$  we have
\begin{align}																		\label{eq:5.28}
&\Big(-\frac{d}{ds}b^{(i)}(s,\tau,\alpha^{(i)}(s,\tau,\hat y'))\Big)e^{-b^{(i)}}v_i^g(s,\tau,y')+
e^{-b^{(i)}}v_{is}^g(s,\tau,y')
\\
\nonumber
=&\frac{\partial\tilde w_i^g(s,\tau,\alpha^{(i)})}{\partial s}+
\sum_{j=1}^{n-1}\frac{\partial \tilde w_i^g(s,\tau,\alpha^{(i)})}{\partial\hat y_j}\alpha_{js}^{(i)}(s,\tau,y').
\end{align}
Returning back  in  (\ref{eq:5.28})  to  $y'=\beta^{(1)}(s,\tau,\hat y')$  coordinates  we get
\begin{align}                   																\label{eq:5.29}
&\frac{\partial\tilde w_i^g(s,\tau,\hat y')}{\partial s}+
\sum_{j=1}^{n-1}\frac{\partial \tilde w_i^g(s,\tau,\hat y')}{\partial\hat y_j}\alpha_{js}^{(i)}(s,\tau,\beta^{(i)}(s,\tau,\hat y'))
\\
&=
\nonumber
e^{-b^{(i)}(s,\tau,\hat y')}v_{is}^g(s,\tau,\beta^{(i)}(s,\tau,\hat y'))
-\frac{d}{ds}b^{(i)}(s,\tau,\alpha(s,\tau,y'))\big|_{y'=\beta^{(i)}}\tilde w_i^g(s,\tau,\hat y').
\end{align}
Subtracting (\ref{eq:5.29}) 
for $i=1$  from  (\ref{eq:5.29})   for $i=2$ and taking  into account  (\ref{eq:5.24})
and (\ref{eq:5.22})  we get
\begin{align} 																	\label{eq:5.30}
&\sum_{j=1}^{n-1}\Big(\alpha_{js}^{(1)}(s,\tau,\beta^{(1)}(s,\tau,\hat y'))-
\alpha_{js}^{(2)}(s,\tau,\beta^{(2)}(s,\tau,\hat y'))\Big)\frac{\partial \tilde w_1^g(s,\tau,\hat y')}{\partial  \hat y_i}
\\
\nonumber
&+\Big(\frac{d}{ds}b^{(1)}(s,\tau,\alpha^{(1)}(s,\tau,y')\Big|_{y'=\beta^{(1)}}-
\frac{d}{ds}b^{(2)}(s,\tau,\alpha^{(2)}(s,\tau,y')\Big|_{y'=\beta^{(2)}}
\Big) \tilde w_1^g(s,\tau,\hat y'))=0
\end{align}
for all  $\tilde w_1^g(s,\tau,\hat y')$  where  $(s,\tau,\hat y')\in \Sigma\times \Gamma_2$.

Fix  $\tau=\tau', 0\leq \tau'<T_0-T_1$.  By the Density Lemma  \ref{lma:3.1}    $\{v_i^g(s,\tau',y')\}$
are dense  in $H_0^1(R_{2T_1}^{(i)}(\tau'))$,  
where  $g\in H_0^1(\Gamma_2\times\{T_1\leq  y_0\leq   T_2-\tau'\})$.

Let  $\tilde R_{2T_1}^{(i)}(\tau')$  be  the image  of  $R_{2T_1}^{(i)}(\tau')\cap\beta^{(i)}(\Sigma\times\overline\Gamma_2)$
under the map  (\ref{eq:5.9}).   Since  $\tilde w_i^g=e^{-b^{(i)}}v_i^g(s,\tau,\beta^{(i)}(s,\tau,\hat y'))$  we have
that $\tilde w_i^g(s,\tau',\hat y')$  are dense in $H_0^1(\tilde R_{2T}^{(i)}(\tau'))$.

The following lemma is similar to  arguments in  [E3],  pp 1749-1750.
\begin{lemma}																\label{lma:5.2}
Since  $\{w_1^g(s,\tau',\hat y'), g\in  H_0^1(\overline\Gamma_2\times\{T_1\leq y_0\leq T_2-\tau'\})\}$  are dense
in  $\tilde R_{2T_1}^{(1)}(\tau')$  we have
\begin{equation}																\label{eq:5.31}
\alpha_{js}^{(1)}(s,\tau',\beta^{(1)}(s,\tau',\hat y'))
=
\alpha_{js}^{(2)}(s,\tau',\beta^{(2)}(s,\tau',\hat y'))
 \ \ \mbox{on}\ \ \tilde R_{2T_1}^{(1)}(\tau'),
 \end{equation}
\begin{equation}																\label{eq:5.32}
\frac{d}{ds}b^{(1)}(s,\tau',\alpha^{(1)}(s,\tau', y'))\Big|_{y'=\beta^{(1)}}
=
\frac{d}{ds}b^{(2)}(s,\tau',\alpha^{(2)}(s,\tau', y'))\Big|_{y'=\beta^{(2)}}
 \ \ \mbox{on}\ \ \tilde R_{2T_1}^{(1)}(\tau').
 \end{equation}
\end{lemma}
{\bf Proof:}  
Let  $\gamma(s,\tau',\hat y')\in  C_0^\infty(\tilde R_{2T_1}^{(1)}(\tau'))$.  There exists  a sequence   
$\tilde w_1^{g_k}(s,\tau',\hat y')$  convergent  to  $\gamma(s,\tau',\hat y')$  in $H_0^1(\tilde R_{2T_1}^{(1)}(\tau'))$.  
Therefore  $\tilde w_1^{g_k}$  converges  weakly  to $\gamma(s,\tau',\tilde y')$.   Passing  in (\ref{eq:5.30})  to the limit
when $k\rw\infty$  we get
\begin{align}																		\label{eq:5.33}
 &\sum_{j=1}^{n-1}\big(\alpha_{js}^{(1)}(s,\tau',\beta^{(1)}(s,\tau',\hat y'))
-
\alpha_{js}^{(2)}(s,\tau',\beta^{(2)}(s,\tau',\hat y'))\big)\frac{\partial\gamma}{\partial\hat y_j}
\\
\nonumber
&+\Big(\frac{d}{ds}b^{(1)}(s,\tau',\alpha^{(1)}(s,\tau', y'))\Big|_{y'=\beta^{(1)}}
\\
\nonumber
&
-
\frac{d}{ds}b^{(2)}(s,\tau',\alpha^{(2)}(s,\tau', y'))\Big|_{y'=\beta^{(2)}}\Big)
\gamma(s,\tau',\hat y')=0.
\end{align}
For any  point  $(s,\hat y')\in  \tilde R_{2T_1}^{(1)}(\tau')$  
we can  find  $n $  
$C_0^\infty(\tilde R_{2T_1}^{(1)})$  functions  \linebreak $\gamma_1(s,\hat y'),...,\gamma_n(s,\hat y')$
such that the determinant of $n\times n$  matrix
$$
\begin{bmatrix}
&\frac{\partial\gamma_1}{\partial \hat y_1}&\dots &\frac{\partial \gamma_1}{\partial \hat y_{n-1}} &\gamma_1
\\
&\dots\                                                            &\dots                                      &\dots                      &\dots
\\
&\frac{\partial\gamma_n}{\partial  \hat y_1}&\dots &\frac{\partial \gamma_n}{\partial \hat y_{n-1}}  &\gamma_n
\end{bmatrix}
$$
is not equal  to zero  at the point  $(s,\hat y')$.  Therefore  (\ref{eq:5.31}),  (\ref{eq:5.32})  hold.
\qed

Repeating  the same arguments for any  $0\leq\tau'\leq T_2-T_1$  we get   that  
(\ref{eq:5.31}),  (\ref{eq:5.32})  hold  for any  $\tau',$  i.e.  it hold on  $\tilde X_{\Gamma_2}^{(1)}
= \bigcup_{0\leq \tau'\leq T_2-T_1}    \tilde R_{2T_1}^{(1)}(\tau')$, 
 since  $\bigcup_{0\leq \tau'\leq T_2-T_1}R_{2T_1}^{(1)}(\tau')=\Sigma\times\Gamma_2$  and  $\tilde X_{\Gamma_2}$ 
 is the image  of  $(\Sigma\times\overline\Gamma_2)\cap\beta^{(1)}(\Sigma\times\overline\Gamma_2)$  
 under the map (\ref{eq:5.9}).
This   proves (\ref{eq:5.25}).   To prove   (\ref{eq:5.26}) 
we note  that 
$$
\frac{d}{ds}b^{(i)}(s,\tau,\alpha^{(i)}(s,\tau,y'))\Big|_{y'=\beta^{(i)}}
=\frac{\partial b^{(i)}}{ds}+\sum_{j=1}^{n-1}\frac{\partial b^{(i)}}{\partial \hat y_j}\alpha_{js}^{(i)}(s,\tau,\beta^{(i)}(s,\tau,\hat y').
$$
Since (\ref{eq:5.31}),  (\ref{eq:5.32}) hold,  we have 
\begin{equation}																				\label{eq:5.34}
\frac{\partial}{\partial s}(b^{(1)}-b^{(2)})
+\sum_{j=1}^{n-1}\frac{\partial }{\partial \hat y_j}(b^{(1)}-b^{(2)})
\alpha_{js}^{(1)}(s,\tau,\beta^{(1)}(s,\tau,\hat y'))=0.
\end{equation}

Equation  (\ref{eq:5.34})  is  a linear  homogeneous equation  for  $b^{(1)}(s,\tau,\hat y')-
b^{(2)}(s,\tau,\hat y')$
on $\hat X_{\Gamma_2}^{(1)}$.  Since  $b^{(1)}=b^{(2)}=0$
when $y_n=0$,  we get 
\begin{equation}																				\label{eq:5.35}
b^{(1)}(s,\tau,\hat y')=b^{(2)}(s,\tau,\hat y')\ \ \ \mbox{on}\ \ \tilde X_{\Gamma_2}^{(1)}.
\end{equation}
It follows from  (\ref{eq:5.35})  and  (\ref{eq:5.22})  that
\begin{equation}																				\label{eq:5.36}
w_1^g(s,\tau,\hat y')=w_2^g(s,\tau,\hat y')  \ \ \mbox{on}\ \ \tilde X_{\Gamma_2}^{(1)},
\end{equation}
where $w_i^g=v_i^g(s,\tau,\beta^{(i)}(s,\tau,\hat y'))$.

\section{The conclusion of the local step}
\label{section 6}
\init

We shall prove the following theorem:
\begin{theorem}																					\label{theo:6.1}
Let  $L_1^{(i)}v_i^g=0,i=1,2.$  Make change of variables
\begin{equation}																				\label{eq:6.1}
\hat s=s,\ \ \hat\tau=\tau,\ \ \hat y'=\alpha^{(i)}(s,\tau,y'),\ i=1,2.
\end{equation}
Let  $\tilde L_1^{(i)}w^g=0$  be the operator  $L_1^{(i)}$  in  the new coordinates.  Then
the coefficients of  $\tilde L_1^{(1)}$   and $\tilde L_1^{(2)}$  are equal  on  $\tilde X_{\Gamma_2}^{(1)}$.
\end{theorem}
{\bf Proof:}
Equations $L_1^{(i)}v_i^g=0$  have the following form  in $(s,\tau,y')$ coordinates  (cf. (\ref{eq:2.27})):
\begin{align}																					\label{eq:6.2}
L_1^{(i)}v_i^g&=-2\frac{\partial}{\partial s}\Big(\frac{\partial}{\partial\tau}+iA_-^{(i)}\Big)v_i^g
-2\Big(\frac{\partial}{\partial\tau}+iA_-^{(i)}\Big)\frac{\partial}{\partial s}v_i^g
\\
\nonumber
&+\sum_{j=1}^{n-1}2\Big(\frac{\partial}{\partial y_j}-iA_j^{(i)}\Big)g_{i0}^{+,j}\frac{\partial}{\partial s}v_i^g
+\sum_{j=1}^{n-1}2\frac{\partial}{\partial s}\Big(g_{i0}^{+,j}\Big(\frac{\partial}{\partial y_j}-iA_j^{(i)}\Big)\Big)v_i^g
\\
\nonumber
&+\sum_{j,k=1}^{n-1}\Big(\frac{\partial}{\partial y_j}-iA_j^{(i)}\Big)g_{i0}^{jk}
\Big(\frac{\partial}{\partial y_k}-iA_k^{(i)}\Big)v_i^g+V_1^{(i)}v_i^g=0,
\end{align}
where  $i=1,2,\  g_{i0}^{+,j}=g_{i0}^{0j},\ V_1^{(i)}$ is the same  as in  (\ref{eq:2.27}).

Making the change of variables (\ref{eq:6.1})  in (\ref{eq:6.2}) we get:
\begin{align}																	\label{eq:6.3}
&\tilde L_1^{(i)}w_i^g(s,\tau,\hat y')
=-2J_i^{-1}(s,\tau,\hat y')\Big(\frac{\partial}{\partial s}+i\tilde A_+^{(i)}\Big)J_i\Big(\frac{\partial}{\partial\tau}+
i\tilde A_-^{(i)}\Big)w_i^g
\\
\nonumber
&-2J_i^{-1}\Big(\frac{\partial}{\partial\tau}+i\tilde A_-^{(i)}\Big)J_i\Big(\frac{\partial}{\partial s}+i\tilde A_+^{(i)}\Big)w_i^g
\\
\nonumber
&-\sum_{j=1}^{n-1}2J_i^{-1}\Big(\frac{\partial}{\partial \tau}+i\tilde A_-^{(i)}\Big)J_i\alpha_{js}^{(i)}(s,\tau,\beta^{(i)})
\Big(\frac{\partial}{\partial y_j}-i\tilde A_j^{(i)}\Big)w_i^g
\\
\nonumber
&-\sum_{j=1}^{n-1}2J_i^{-1}\Big(\frac{\partial}{\partial y_j}-i\tilde A_j^{(i)}\Big)J_i\alpha_{js}^{(i)}(s,\tau,\beta^{(i)})
\Big(\frac{\partial}{\partial \tau}+i\tilde A_-^{(i)}\Big)w_i^g
\\
\nonumber
&+\sum_{j,k=1}^{n-1}J_i^{-1}\Big(\frac{\partial}{\partial y_j}-i\tilde A_j^{(i)}\Big)J_i
\tilde g_{i0}^{jk}\Big(\frac{\partial}{\partial y_k}-i\tilde A_k^{(i)}\Big)w_i^g
\\
\nonumber
&+V_1^{(i)}(s,\tau,\beta^{(i)}(s,\tau,\hat y'))w_i^g(s,\tau,\hat y')=0,
\end{align}
where  $w_i^g(s,\tau,\hat y')=v_i^g(s,\tau,\beta^{(i)}(s,\tau,\hat y'))$,
\begin{multline}																		\label{eq:6.4}
\tilde g_{i0}^{jk}(s,\tau,\hat y') =\sum_{p,r=1}^{n-1}g_{i0}^{pr}(s,\tau,\beta^{(i)})\alpha_{jy_p}^{(i)}(s,\tau,\beta^{(i)})
\alpha_{ky_r}^{(i)}(s,\tau,\beta^{(i)})
\\
-2\alpha_{js}^{(i)}\alpha_{k\tau}^{(i)}-2\alpha_{j\tau}^{(i)}\alpha_{ks}^{(i)}
+2\sum_{p=1}^{n-1}g_{i0}^{+,p}(\alpha_{js}\alpha_{ky_p}+\alpha_{jy_p}\alpha_{ks}).
\end{multline}
We used in (\ref{eq:6.3}) that (see (\ref{eq:5.12}))
\begin{equation}                  																\label{eq:6.5}
\tilde g_{i0}^{+,j}(s,\tau,\hat y')=
\sum_{k=1}^{n-1}g_{i0}^{+,k}(s,\tau,\beta^{(i)}(s,\tau^{(i)},\hat y'))
\alpha_{jy_k}^{(i)}(s,\tau,\beta^{(i)})-\alpha_{j\tau}^{(i)}(s,\tau,\beta^{(i)})=0.
\end{equation}
Also we have 
\begin{multline}                  																\label{eq:6.6}
\tilde g_{i0}^{-,j}(s,\tau,\hat y')=
\sum_{k=1}^{n-1}g_{i0}^{-,k}(s,\tau,\beta^{(i)}(s,\tau,\hat y'))\alpha_{jy_k}^{(i)}(s,\tau,\beta^{(i)})
-\alpha_{j s}^{(i)}(s,\tau,\beta^{(i)})
\\
=
-\alpha_{js}^{(i)}(s,\tau,\beta^{(i)}),
\end{multline}
since $g_{i0}^{-,k}=0$  (cf.  (\ref{eq:6.2})).

Note that $A_+^{(i)}, A_-^{(i)}, A_j^{(i)}$  and  $\tilde A_+^{(i)},\tilde A_-^{(i)},\tilde A_j^{(i)}$
are related by the equality
\begin{equation}																		\label{eq:6.7}
 A_+^{(i)}ds+A_-^{(i)}d\tau -\sum_{j=1}^{n-1}A_j^{(i)}dy_j = 
\tilde A_+^{(i)}d\hat s+\tilde A_-^{(i)}d\hat\tau-\sum_{j=1}^{n-1}\tilde A_j^{(i)}d\hat y_j,
\end{equation}
where
$A_+^{(i)}=0,i=1,2,\ s=\hat s,\ \tau=\hat \tau,\ y_j=\beta_j(s,\tau,\hat y')$.

Note  that  (\ref{eq:5.21}),   (\ref{eq:5.26})  imply
\begin{equation}																	\label{eq:6.8}
J_1(s,\tau,\hat y')=J_2(s,\tau,\hat y')\ \ \ \mbox{in}\ \ \tilde X_{\Gamma_2}^{(1)}.
\end{equation}
The first order term containing $\frac{\partial}{\partial\tau}$  in (\ref{eq:6.3})  is equal to 
\begin{align}																	\label{eq:6.9}
&-2i\tilde A_+^{(i)}\Big(\frac{\partial}{\partial \tau}\Big) w_i^g-2iJ_i^{-1}\Big(\frac{\partial}{\partial\tau}\Big) J_i\tilde A_+^{(i)}w_i^g
+2i\sum_{j=1}^{n-1}\tilde A_j^{(i)}\alpha_{js}^{(i)}\frac{\partial}{\partial\tau}w_i^g
\\
\nonumber
&+i\sum_{j=1}^{n-1}2J_i^{-1}\Big(\frac{\partial}{\partial\tau}\Big)J_i\alpha_{js}^{(i)}\tilde A_j^{(i)}w_i^g.
\end{align}
It follows from (\ref{eq:6.7})  that
\begin{equation}															      \label{eq:6.10}
 A_+^{(i)}=\tilde A_+^{(i)}-\sum_{j=1}^{n-1}\tilde A_j^{(i)}\alpha_{js}(s,\tau,\beta^{(i)}(s,\tau,\hat y')).
\end{equation}
Since $ A_+^{(i)}=0$  we have  that   (\ref{eq:6.10})  implies that (\ref{eq:6.9})  is equal  to zero.

Taking into account  that  $\alpha_{js}^{(1)}(s,\tau,\beta^{(1)})=\alpha_{js}^{(2)}(s,\tau,\beta^{(2)}),\ 1\leq j\leq n-1$,
$J_1=J_2$
and 
$w_1^g(s,\tau,\hat y')=w_2^g(s,\tau,\hat y')$  we get
that  $\tilde L_1^{(1)}-\tilde L_1^{(2)}$  is a differential  operator  in $\frac{\partial}{\partial s},
\frac{\partial}{\partial y_1},...,\frac{\partial}{\partial y_n}$.  We have
\begin{equation}																\label{eq:6.11}
(\tilde L_1^{(1)}-\tilde L_1^{(2)})w_1^g=0.
\end{equation}
Since  $\{w_1^g,g\in C_0^\infty(\Gamma^{(1)}\times[T_1,T_2-\tau']\}$  are  dense
in  $H_0^1(\tilde R_{1T_1}(\tau'))$ we get as in  Lemma \ref{lma:5.2}  (cf.  [E3])  that 
all coefficients of  $\tilde L_1^{(1)}$  and  $\tilde L_1^{(2)}$  are equal  in $\tilde R_{1T_1}^{(1)}(\tau')$.
Since  $\tau'\in[0,T_2-T_1]$  is arbitrary,  we get  that  on  $\hat X_{\Gamma_2}^{(1)}$:
\begin{align}																		\label{eq:6.12}
&\tilde g_{10}^{jk}(s,\tau,\hat y')=\tilde g_{20}^{jk}(s,\tau,\hat y'),\ 1\leq j,k\leq n-1,
\\
\label{eq:6.13}
&\tilde A_-^{(1)}=\tilde A_-^{(2)}, \tilde A_j^{(1)}=\tilde A_j^{(2)}(s,\tau,\hat y), 1\leq j\leq n-1,
\\
\label{eq:6.14}
&V_1^{(1)}(s,\tau,\beta^{(1)}(s,\tau,\hat y'))=V_1^{(2)}(s,\tau,\beta^{(2)}(s,\tau,\hat y')).
\end{align}
Therefore  $\tilde L_1^{(1)}=\tilde L_1^{(2)}$  in  $\tilde X_{\Gamma_2}^{(1)}$.
\qed

This completes the proof  of Theorem \ref{theo:6.1}. 
Let $L_i'v_i^g=0$  be the equation  of the form  (\ref{eq:2.24}).  Making  the change  of variables  (\ref{eq:5.10})  we get
the equation 
$\tilde L_i'w_i^g=0,i =1,2,$  on $\tilde X_{\Gamma_2}^{(1)}$.

Note  that $w_1^g=w_2^g$  on  $\tilde X_{\Gamma_2}^{(1)}$.  We shall prove  that  $\tilde L_1'=\tilde L_2'$  on 
 $\tilde X_{\Gamma_2}^{(1)}$.

Let  $\hat g_i^{+,-},\hat g_i^{+,j},\hat g_j^{jk},1\leq j\leq n-1,1\leq k\leq n-1$,  be  the 
inverse metric  tensor  of  $L_i'$.
Note that for $L_1^{(i)}$
  we have
(cf.  (\ref{eq:2.27}))
$$
g_{i0}^{+,j}=\frac{\hat g_i^{+,j}}{\hat g_i^{+,-}},\ \ \ g_{i0}^{jk}=\frac{\hat g_i^{jk}}{\hat g_i^{+,-}},\ \ i=1,2.
$$
Therefore the equation  $\tilde L_i'w_i^g=0$  has  the inverse metric  tensor with elements  (cf. (\ref{eq:6.4}),  (\ref{eq:6.5}),
(\ref{eq:6.6}))
\begin{align}																\label{eq:6.15}
&\tilde g_i^{jk}=\hat g_i^{+,-}g_{i0}^{jk},\ \ \  1\leq j,k \leq n-1,
\\
\nonumber
&\tilde g_i^{-,k}=-\hat g_i^{+,-}\alpha_{ks}^{(i)}(s,\tau,\beta^{(i)}),\ \ \tilde g_i^{+,k}=0,\ \  1\leq k\leq n-1,\ \ i=1,2.
\end{align}
Since $\alpha_{ks}^{(1)}=\alpha_{ks}^{(2)}$  and  $\tilde g_{10}^{jk}=\tilde g_{20}^{jk}$
(see  (\ref{eq:6.12})),
we get that the metric  tensors of $\tilde L_1'$  and  $\tilde L_2'$
are equal  if we  can  prove  that
\begin{equation}																\label{eq:6.16}
\tilde g_1^{+,-}(s,\tau,\beta^{(1)}(s,\tau,\hat y'))=
\tilde g_2^{+,-}(s,\tau,\beta^{(2)}(s,\tau,\hat y')).
\end{equation}
We shall prove  first  that
\begin{equation}																\label{eq:6.17}
g_1^{(1)}(s,\tau,\beta^{(1)})=g_1^{(2)}(s,\tau,\beta^{(2)}),
\end{equation}
where  $ g_1^{(i)}=\big|\det[\hat g_i^{jk}]_{j,k=1}^{n-1}\big|^{-1}$  (see  (\ref{eq:2.22})).

Note that $V_1^{(i)}(s,\tau,y')$  has  the form  (\ref{eq:2.25})  for  $i=1,2$,  where  $A^{(i)}=\ln(g_1^{(i)})^{\frac{1}{4}}$.  Making the  change of variables (\ref{eq:6.1}) we get (cf. (\ref{eq:6.14})) 
\begin{equation}																\label{eq:6.18}
V_1^{(1)}(s,\tau,\beta^{(1)}(s,\tau,\hat y'))-V_2^{(2)}(s,\tau,\beta^{(2)}(s,\tau,\hat y'))=0.
\end{equation}
Note that  the metric  tensors for $\tilde L_1'$  and  $\tilde L_2'$  are  equal  on $\tilde X_{\Gamma_2}$.
Let $\tilde A^{(i)}(s,\tau,\hat y')=A^{(i)}(s,\tau,\beta^{(i)}(s,\tau,\hat y'))$.

Using the equality 
$$
\tilde A_{y_j}^{(1)}\tilde A_{y_k}^{(1)}-\tilde A_{y_j}^{(2)}\tilde A_{y_k}^{(2)}
=(\tilde A_{y_j}^{(1)}-\tilde A_{y_j}^{(2)})A_{y_k}^{(1)}+
(\tilde A_{y_k}^{(1)}-\tilde A_{y_k}^{(2)})\tilde A_{y_j}^{(2)}
$$
and similar  equality  involving derivatives  in $s$  and $\tau$ we can represent (\ref{eq:6.18}) 
as  a homogeneous second order   hyperbolic  equation  in $\tilde A^{(1)}-\tilde A^{(2)}$
with the coefficients depending  on  $\tilde A^{(1)}$  and $\tilde A^{(2)}$.
Since the Cauchy  data  for  $\tilde A^{(1)}-\tilde A^{(2)}=0$
at  $y_n=0$  (cf.  Lemma 2.1)  we get, by the uniqueness of  the Cauchy  problem  (cf.  [T],  [RZ]),  that
$\tilde A^{(1)}=\tilde A^{(2)}$
in $\tilde X_{\Gamma_2}^{(1)}$.  Therefore  (\ref{eq:6.17}) holds.

Note that
$\tilde g_i^{jk}=\tilde g_i^{+,-}\tilde g_0^{jk}$.  Therefore 
$$
g_1^{(i)}=\det[\tilde g_i^{jk}]_{j,k=1}^{n-1}=(\tilde g_i^{+,-})^{n-1}\det[g_{0i}^{jk}]_{j,k=1}^{n-1}.
$$
Since
$\tilde g_{10}^{jk}=\tilde g_{20}^{jk}$  and  (\ref{eq:6.17})  holds,  we
get
\begin{equation}																		\label{eq:6.19}
(\tilde g_1^{+,-}(s,\tau,\beta^{(1)}))^{n-1}=(\tilde g_2^{+,-}(s,\tau,\beta^{(1)}))^{n-1},
\end{equation}
and this proves  (\ref{eq:6.16}),  since we assumed that $n>1$.
   Therefore  metric  tensors of  $\tilde L_1'$  and  $\tilde L_2'$
are equal.  Combining  this with  (\ref{eq:6.13}),  (\ref{eq:6.14})  we get $\tilde L_1'=\tilde L_2'$  on  
$\tilde X_{\Gamma_2}^{(1)}\supset \Sigma\times\overline\Gamma_1 $. 

{\bf Remark 6.1}. Change $\Gamma_2$  to $\Gamma_1$.  We have 
 $\beta^{(i)}(\Sigma\times \overline \Gamma_1)\subset \beta^{(i)}(\Sigma\times \overline \Gamma_2)$.
 Since  $\beta^{(i)}(\Sigma\times \overline \Gamma_1)\subset X_{1T_1}^{(i)}$ and     
 $X_{1T_1}^{(i)}\subset(\Sigma\times \overline \Gamma_2)$, we get
$\beta^{(i)}(\Sigma\times \overline \Gamma_1)\subset (\Sigma\times\overline \Gamma_2)$.
Therefore,  $\beta^{(i)}(\Sigma\times \overline \Gamma_1)\subset 
 (\Sigma\times \overline \Gamma_2)\cap \beta^{(i)}(\Sigma\times \overline \Gamma_2)=
X_{\Gamma_2}^{(i)}$.
Applying  the map  (\ref{eq:5.9})  to  $\beta^{(i)}(\Sigma\times \overline \Gamma_1)\subset
X_{\Gamma_2}^{(i)}$
we get
 $\Sigma\times \overline \Gamma_1\subset  \tilde X_{\Gamma_2}^{(i)}$.
Therefore,   $\tilde L_1'=\tilde L_2'$  on 
$\Sigma\times\overline\Gamma_1$.
\\
\

We shall summarize the results of \S\S 2-6.
\begin{theorem}[Local step]														\label{theo:6.2}
Consider
two initial boundary value problems
\begin{align}																			\label{eq:6.20}
&L^{(i)}u_i=0\ \ \ \mbox{in}\ \ D_0^{(i)}\times\R,
\\
\nonumber
&u_i=0\ \ \ \mbox{for}\ \ x_0\ll 0,
\\
\nonumber
&u_i\big|_{\partial D_0^{(i)}}=f,\ \ i=1,2,
\end{align}
where
$L^{(i)}$  have the form  (\ref{eq:1.1}).  Suppose  $\Gamma_0\subset \partial D_0^{(1)}\cap\partial D_0^{(2)}$
and 
suppose the BLR  condition holds  for  $L^{(1)}$  on $[t_0,T_{t_0}]$.  Suppose
 the corresponding  DN operators  $\Lambda^{(i)}$  are equal  on $\Gamma_0\times(t_0,T_2),\ T_2\geq T_{t_0}$,
i.e.  $\Lambda^{(1)}f=\Lambda^{(2)}f$    on  $\Gamma_0\times(t_0,T_{2})$ for all  $f$ with
support  in $\overline{\Gamma_0}\times[t_0,T_{2}]$.   Let  $T_2-T_1$  be small. 
Suppose  coefficients  of $L^{(1)}$  and  $L^{(2)}$  are  analytic  in $x_0$.

Let  $\varphi^{(i)}$  be the changes  of variables  (\ref{eq:2.14}) for  $i=1,2$  and  let  $\beta^{(i)},i=1,2,$  be 
the changes of variables  (\ref{eq:5.10}).  Let  $c_i$  be the  gauge transformation (\ref{eq:2.20}), (\ref{eq:2.21})  for
$i=1,2.$  Then 
\begin{equation}																			\label{eq:6.21}
\beta^{(1)}\circ c_1\circ\varphi^{(1)}\circ L^{(1)}=\beta^{(2)}\circ c_{2}\circ \varphi^{(2)}\circ L^{(2)}\ \ \mbox{on}\ \ 
\Sigma\times \overline\Gamma_1,
\end{equation}
where 
\\
$\Sigma=\{(s,\tau),s\geq 0,\tau\geq 0, s+\tau\leq T_2-T_1\}=
\{(y_0,y_n):0\leq y_n\leq \frac{T_2-T_1}{2},\linebreak  T_1+y_n<y_0<T_2-y_n\}$.
\end{theorem}


\section{The global step}
\label{section 7}
\init
 Let $L^{i}u_i=0$ in $D_i=D_0^{(i)}\times \R,i=1,2,\ u_i=0$ for $x_0\ll 0,\ \partial D_0^{(1)}\cap\partial D_0^{(2)}\supset \Gamma_0$ and  
 $u_i\big|_{\partial D_0^{(i)}\times\R}=f,i=1,2,\ f$  has a compact  support  in $\overline{\Gamma_0}\times\R$.

First  we extend  the Theorem \ref{theo:6.2} for a larger  time interval.
 
 Let $[t_1,t_2]$  be an arbitrary  time interval.  Let $[t_0,T_{t_0}]$  be such that $T_{t_0}\leq t_1$  and
 the BLR  condition  holds  on  $[t_0,T_{t_0}]$.  
 Thus  the  BLR condition  is satisfied  on $[t_0,t]$  for any  $t\in  [t_1,t_2]$.  Let
 $\Gamma_1$  be  arbitrary  connected  part  of  $\Gamma_0,\ \overline{\Gamma_1}\subset \Gamma_0$.  Note  that  
 we  do not require    $\overline{\Gamma_1}$  to be small.
 
 Let $\psi_{0i}^\pm(x_0,x',x_n),i=1,2,$  be the solution  of the form (\ref{eq:2.4}) in  
 $[t_0-1,t_2+1]\times\overline{\Gamma'}\times[0,\e_\pm]$  
 where  $\overline{\Gamma_1}\subset \Gamma'\subset\Gamma_0$.
 
 We  impose  the  following  initial  conditions  on  $\psi_{0i}^\pm, i=1,2$,
 \begin{equation}																				\label{eq:7.1}
 \psi_{0i}^+\big|_{x_n=0}=x_0,\ \ \psi_{0i}^-\big|_{x_n=0}=-x_0.
 \end{equation}
 Such solutions exist  in 
  $[t_0-1,t_2+1]\times\overline{\Gamma'}\times[0,\e_0]\subset D_0^{(i)} \times\R$,  when  $\e_0$  is small.  
  We choose  $\psi_{0i}^\pm$  such   that  (\ref{eq:2.6})  is satisfied and we choose  $\e_1>0$  such  that
  $\e_1\leq \e_0$
 and  $\{0<x_n<\e_1,x'\in\overline{\Gamma'},x_0\in [t_0-1,t_2+1]\}$  do not intersect
 $\partial D_0^{(i)}\times\R$.
 
 Let  $\varphi_{ji}(x_0',x',x_n),1\leq j\leq n-1,$  be the solutions of the linear equations  (cf.  (\ref{eq:2.7}))
 \begin{equation}																				\label{eq:7.2}
 \sum_{p,k=0}^n g_i^{pk}(x)\psi_{0ix_p}^-\varphi_{jix_k}=0\ \ \mbox{in}\ \ 
 [t_0-1,t_2+1]\times\overline{\Gamma'}\times[0,\e_1]
 \end{equation}
 with initial conditions
 \begin{equation}                    				 													\label{eq:7.3}
 \varphi_{ji}\big|_{x_n=0}=x_j,\ \ \ 1\leq j\leq n-1.
 \end{equation}
 Similarly to  (\ref{eq:2.14})  consider  the map $(y_0^{(i)}(x),y_i'(x),y_n^{(i)}(x))=(\varphi_0^{(i)},\varphi_i',\varphi_n^{(1)}),\ 
 x\in[t_0-1,t_2+1]\times\overline{\Gamma'}\times[0,\e_1]$,
 where 
 \begin{align}																					\label{eq:7.4}
 &y_0^{(i)}(x)=\frac{\psi_{0i}^+-\psi_{0i}^-}{2},
 \\
 \nonumber
 &y_j^{(i)}(x)=\varphi_{ji}(x),
 \\
 \nonumber
 &y_n^{(i)}(x)=-\frac{\psi_{0i}^++\psi_{0i}^-}{2}.
 \end{align}
 As in  (\ref{eq:2.15})   we have  that  the map  $(x_0,x',x_n)\rw(y_0,y',y_n)$  is the identity  when  $x_n=0$:
 \begin{equation}																				\label{eq:7.5}
 y_0^{(i)}\big|_{x_n=0}=x_0,\ \ \ y_j^{(i)}\big|_{x_n=0}=x_j,\ 1\leq j\leq n-1,\ \ \ y_n^{(i)}\big|_{x_n=0}=0.
 \end{equation}
 Let $u_s=\frac{1}{2}(u_{y_0}-u_{y_n}),u_\tau=-\frac{1}{2}(u_{y_0}+u_{y_n})$.
 Making the change of variables  (\ref{eq:7.4})   in $L^{i}u_i=0$,  the  gauge  transformation (\ref{eq:2.18}),  (\ref{eq:2.21})
 and the change   of unknown  function  (\ref{eq:2.26}),  we get  in $t_0\leq y_0\leq t_2,\ 0\leq y_n\leq T_0,\ 
 y'\in \overline{\Gamma'},\ T_0$ is small,  the equation of the form
 $$
 L_1^{(i)}u_1^{(i)}=0,\ \ y\in\hat\Omega_0,
 $$
 where $L_1^{(i)}$  has the form  (\ref{eq:2.28}).
 Here 
 \begin{equation}																			\label{eq:7.6}
 \overline{\Gamma_1}\subset\Gamma'\subset \Gamma_0,\ 
 \hat \Omega_0\stackrel{def}{=}[t_0,t_2]\times\overline{\Gamma'}\times[0,T_0].
 \end{equation}
 We assume that $u_1^{(i)}$  satisfy  the zero  initial  conditions
 $$
 u_1^{(i)}=\frac{\partial u_1^{(i)}}{\partial y_0}
 =0\ \ \mbox{when}\ \ y_0=t_{0}
 $$
 and  
 $$
 u_1^{(i)}\big|_{y_n=0}=f,\ \ i=1,2.
 $$
 We also assume that  DN  operators  for  $L^{i}$  and subsequently  for  $L_1^{(i)}$  are  equal on  $[t_0,t_2)\times
 \overline{\Gamma'}$.
 
 Note that the change of variables
 \begin{equation}																					\label{eq:7.7}
 \hat y_n=y_n,\ \ \hat y_0=y_0,\ \ \hat y_j'=\alpha_j^{(i)}(y_0,y',y_n),\ 1\leq j\leq n-1,
 \end{equation}
 where $ \alpha^{(i)}$  are the same  as in (\ref{eq:5.9}),  (\ref{eq:5.11}),
 are also defined globally on $\hat\Omega_0$.

 Let  $[T_1,T_2]\subset [t_1,t_2]$  be arbitrary  such that $T_2-T_1=2T_0.$
 
 Applying Theorem \ref{theo:6.1}  to  the interval  $[T_1,T_2]$  we get  that 
 the coefficients of  $\tilde L_1^{(1)}$  and  $\tilde L_1^{(2)}$  and  the coefficients of  $L_1'$  and  $L_2'$  are  equal 
 on $\Sigma_{T_1T_2}\times\overline{\Gamma_1}$  where $\Sigma_{T_1T_2}
 =\{0\leq y_n\leq T_0,T_1+y_n\leq y_0\leq T_2-y_n\}$. 
  We assume  that $\Gamma'\supset\overline{\Gamma_1}$  is such that
 $\overline{\Gamma_2}\subset \overline{\Gamma_3}\subset\Gamma_0'$  for 
 all  $[T_1,T_2]\subset[t_1,t_2]$.  Here  
 $\Gamma_2,\Gamma_3$  are  defined as in \S3.  Note that  $\Gamma_2,\Gamma_3$  may depend  on  $[T_1,T_2]$.
 
 If two intervals $[T_1,T_2]$  and  $[T_1',T_2']$  intersect,  then the coefficients of  $\tilde L_1^{(1)}$  and
 $\tilde L_1^{(2)}$  coincide  in $(\Sigma_{T_1T_2}\bigcup\Sigma_{T1'T_1'})\times \Gamma_1$.
 
 Therefore  coefficients of $\tilde L_1^{(1)}$  and  $\tilde L_1^{(2)}$  and consequently the coefficients of $L_1'$  and
 $L_2'$ (cf.  (\ref{eq:6.2}),  (\ref{eq:6.3})) coincide for $0\leq y_n\leq T_0,\ y'\in \overline{\Gamma_1},\ 
 t_1+T_0<y_0<t_2-T_0.$
 
 Therefore we proved 
 \begin{lemma}																			\label{lma:7.1} 
 Suppose $[t_1,t_2]$  is arbitrary large, $T_0>0$ is small,  $t_0$  is such  that the BLR condition is satisfied on  $[t_0,t_1]$.
 Let $\Omega_0=\{y_0\in[t_0+T_0,t_2-T_0],\ y'\in\overline{\Gamma_1}, 
 y_n\in[0,\frac{T_0}{2}]\}$.
 Assume that the coefficients of  $L^{(i)}$  are analytic  in $x_0,i=1,2$. 
 Then
 $$
 \beta^{(1)}\circ c_1\circ \varphi^{(1)}\circ L^{(1)}=
   \beta^{(2)}\circ c_2\circ \varphi^{(2)}\circ L^{(2)}  \ \ \ \mbox{on} \ \ \Omega_0.\quad\quad\qquad\qquad
   \qed
  $$
   \end{lemma}
 Let  $\Omega_i=(\beta^{(i)}\varphi^{(i)})^{-1}\Omega_0,\ i=1,2$.  Note that $\Omega_i\subset D_0^{(i)}\times
 [t_0-1,t_2+1]$  since  $T_0$  is small.  We have that $\Phi_2=(\beta^{(1)}\varphi^{(1)})^{-1}\beta_2\varphi^{(2)}$
 maps  $\Omega_2$  onto $\Omega_1$.  Note  that  $\partial\Omega_1\cap\partial\Omega_2\supset \Gamma_1\times[t_0,t_2]$  and
  $\Phi_2=I$  on $\Gamma_1\times[t_0+T_0,t_2-T_0]$.  Note also  that $\beta^{(i)}\circ c_i$  can be represented   as 
 $c_i'\circ\beta^{(i)}$  where  $c_i'$  is the gauge  transformation in $(y_0,\hat y',y_n)$  coordinates.  Analogously,
 $(\beta^{(1)}\circ c_1\circ\varphi^{(1)})^{-1}\beta^{(2)}\circ c_2\circ\varphi=c_3\circ\Phi_2$,  where $c_3$  is the gauge 
 transformation.  Therefore    
 \begin{equation}																	\label{eq:7.8}
 c_3\circ \Phi_2\circ L^{(2)}=L^{(1)}\ \ \ \mbox{in}\ \ \Omega_1.
 \end{equation}
 Let $B$ be a smooth domain  in $D_0^{(1)}$  such  that 
 $\partial B\cap\partial D_0^{(1)}=\gamma_1\subset\Gamma_0$.  Suppose $B$  is 
 small and such  that $B\times[t_1+1,t_2-1]\subset
 \Omega_1$.
 
 Let $S_2=\Phi_2^{-1}(B\times[t_1+1,t_2-1]\subset D_0^{(2)}\times\R$ and let 
 $S_2^+=\Phi_2^{-1}(B\times\{x_0=t_2-1\}), \ S_2^-=\Phi_2^{-1}(B\times\{x_0=t_1+1\})$.
 Let $\tilde S_2^\pm$  be space-like surfaces in $D_0^{(2)}\times \R$  such  that
 $\tilde S_2^+$  is the extension of $S_2^+$ and $\tilde S_2^-$ is the extension of  $S_2^-$.
 
 We assume  that  the projections of  $\tilde S_2^\pm$  on  $D_0^{(2)}$  is 
 $D_0^{(2)}$.  Let  $D_1^{(2)}$  be the domain  in  $D_0^{(2)}\times\R$  bounded  by  $\tilde S_2^+$
 and  $\tilde S_2^-$ (cf.  Fig. 7.1).
 
 It follows  from [Hi],  Chapter  8,  that  there exists an extension $\tilde\Phi_2$  of  $\Phi_2$  from  $S_2\subset D_1^{(2)}$ 
  to  
 $D_1^{(2)}$  such that  $\tilde\Phi_2\big|_{\Gamma_0\times [t_1+1,t_2-1]}=I$.
 
 Define  $\overline D_1^{(3)}=\tilde\Phi_2(\overline D_1^{(2)})$.
 There exists  also an extension  $\tilde c_3$  of  the gauge $c_3$  from $S_2$  to  $D_1^{(2)}$  such  that  
 $\tilde c_3=1$  on
 $\Gamma_0\times[t_1+1,t_2-1]$.  Let $L^{(3)}=\tilde c_3\circ\tilde\Phi_2\circ L^{(2)},\ L^{(3)}$  is defined  on 
  $D_1^{(3)}$.  Thus 
 $L^{(3)}=L^{(1)}$  on  $B\times[t_1+1,t_2-1]$.  Note that
 $D_1^{(3)}\cap(D_0^{(1)}\times [t_1+1,t_2-1])\supset
 B\times[t_1+1,t_2-1],
 \ 
  \partial' D_1^{(3)}\cap (\partial D_0^{(1)}\times[t_1+1,t_2-1])\supset \Gamma_0\times
 [t_1+1,t_2-1].
 $
 We denote  by  $\partial'D_1^{(3)}$  the lateral (time-like)  part of $\partial D_1^{(3)}$  and by $\partial_\pm D_1^{(3)}$   the top
 and the bottom  space-like  parts  of  $\partial D_1^{(3)}$,  i.e.  $\partial D_1^{(3)}=\partial' D_1^{(3)}\cup
 \partial_+D_1^{(3)}\cup\partial_- D_1^{(3)}$.
 \\
\begin{tikzpicture}[scale=0.55]

\draw (-2,0) arc (180: 360 :2 and 1);
\draw[dashed] (2,0) arc(0:180:2 and 1);
\draw (-2,-1) -- (-2,-4);

\draw (2,0) -- (2,-4);
\draw (-2,-4) arc (180: 360 :2 and 1);
\draw[dashed] (2,-4) arc(0:180:2 and 1);
\draw(2,1) -- (2,-5);
\draw(-2,1)-- (-2,-5);

\draw (2,0) arc (0: 180 :2 and 1.5);
\draw(-2,-4) arc(180:360:2 and 1.5);

\draw (-1,-0.9).. controls (-1.1,0.2) and (-1,1) .. (0,1.5);
\draw [dashed](0,1.5).. controls (0.5,1.4).. (1,0.9); 

\draw(-1,-4.9) .. controls (-0.6,-5.3) .. (0,-5.5);
\draw[dashed](0,-5.5).. controls   (1,-5.1)  and (1.2,-4.1) ..  (1,-3.1);

\draw(-1,-4.9)--(-1,-0.9);
\draw[dashed](1,-3.1)--(1,0.9);

\draw (-1.7,1.5)  node {$\tilde S_2^+$};
\draw (1.8,-5.45)  node {$\tilde S_2^-$};
\draw (3.4, -2.5) node {$D_0^{(2)}\times \R$};

\end{tikzpicture}
\\
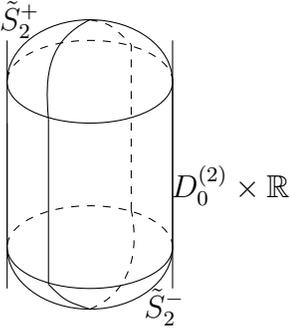
{\bf Fig. 7.1.}  The almost   cylindrical domain  $D_1^{(2)}$ is the part  of $D_0^{(2)}\times\R$  bounded  from above and  from  below
by space-like surfaces $\tilde S_2^+$ and  $\tilde S_2^-$.
\\
\

 The following lemma  is the key  lemma of  this section.  
 It allows to reduce the solution  of the inverse problem  to  an inverse  problem  in a smaller domain.
 \begin{lemma}																	\label{lma:7.2}
 Consider  two  initial-boundary value  problem  
 $L^{(1)}u_1=0$  in  $D_0^{(1)}\times[t_1,t_2]$  and  $L^{(3)}u_3=0$  in    $D_1^{(3)}$,
 $$
 u_1\big|_{x_0=t_1}=\frac{\partial  u_1}{\partial x_0}\Big|_{x_0=t_1}=0,\ \ \ x\in  D_0^{(1)},
 $$
 $$
 u_3\big|_{\partial_-D_1^{(3)}}=\frac{\partial u_3}{\partial x_0}\Big|_{\partial_- D_1^{(3)}}=0,
 \ \ 
 u_1\big|_{\partial D_0^{(1)}\times[t_1,t_2)}=f_1,\ u_3\big|_{\partial' D_1^{(3)}}=f_3.
 $$
 We assume that  $(\partial D_0^{(1)}\times [t_1,t_2]) \cap\partial D_1^{(3)}
 \supset\Gamma_0\times[t_1,t_2]$.
 Assume that  $L^{(1)}=L^{(3)}$  in a smooth domain  $B\times[t_1,t_2]$  where $B\times[t_1,t_2]\subset
 (D_0^{(1)}\times[t_1,t_2])\cap D_1^{(3)},\ \gamma_1=\partial D_0^{(1)}\setminus\Gamma_0,\ 
 \tilde\Gamma_3=\partial' D_1^{(3)}\setminus(\Gamma_0\times(t_1,t_2)),\ \partial B=\gamma_0\cup\gamma_0'$,
 where  $\gamma_0,\gamma_0'$  are smooth,
 $\gamma_0\subset\Gamma_0$,  (cf.  Fig. 7.2).
 
 Suppose $\Lambda_1=\Lambda_3$  on  $\Gamma_0\times [t_1,t_2]$.
 
 Consider
 $L^{(1)}u_1=0$  and  $L^{(3)}u_3=0$  in  smaller  domains 
 $(D_0^{(1)}\setminus B)\times(t_1+\delta,t_2-\delta)$
 and  $(D_1^{(3)}\cap(t_1+\delta,t_2-\delta))\setminus  (B\times(t_1+\delta,t_2-\delta))$.  Note that
 $\partial(D_0^{(1)}\setminus B)\supset(\Gamma_0\setminus\gamma_0)
 \cup\gamma_0')$.  Then
 $\Lambda_1'=\Lambda_3'$ are equal  on 
  $((\Gamma_0\setminus\gamma_0)\cup\gamma_0')\times(t_1+\delta,t_2-\delta)$ for some
 $\delta>0$.  Here  $\Lambda_1',\Lambda_3'$   are  DN  operators  for  the 
 initial-boundary value problem
 \begin{align}
 \nonumber
& L^{(1)}u_i'=0\ \ \mbox{in}\ \ (D_0^{(1)}\setminus B)\times(t_1+\delta,t_2-\delta),
 \\
 \nonumber
 & L^{(3)}u_3'=0\ \ \mbox{in}\ \ (D_1^{(3)}\cap (t_1+\delta,t_2-\delta))\setminus (B\times(t_1+\delta,t_2-\delta)),
 \\
 \nonumber
 &u_1'\big|_{x_0=t_1+\delta}=\frac{\partial u_1'}{\partial x_0}\Big|_{x_0=t_2+\delta}=0,
 \\
 \nonumber 
 &u_3'\big|_{\partial_-(D_1^{(3)}\cap (t_1+\delta,,t_2-\delta))}
  =\frac{\partial u_3'}{\partial x_0}\Big|_{\partial_- (D_1^{(3)}
 \cap(t_1+\delta,t_2-\delta))}=0,
 \\
 \nonumber
 &u_1'\big|_{((\Gamma_0\setminus\gamma_0)\cup\gamma_0')\times(t_1+\delta,t_2-\delta)}=f,\ \ \
 u_1'\big|_{(\partial  D_0^{(1)}\setminus\Gamma_0)\times(t_1+\delta,t_2-\delta)}=0,
 \\
  \nonumber
  &u_3'\big|_{((\Gamma_0\setminus\gamma_0)\cup\gamma_0')\times(t_1+\delta,t_2-\delta)}=f,\ \ \
 u_3'\big|_{((\partial'  D_1^{(3)}\cap (t_1+\delta,t_2-\delta))\setminus(\Gamma_0\times(t_1+\delta,t_2-\delta))}=0,
 \\
 \nonumber
 &\mbox{supp}\,f\subset(((\Gamma_0\setminus\gamma_0)\cup\gamma_0')\times(t_1+\delta,t_2-\delta)).
 \end{align}
 \end{lemma}




\

\begin{tikzpicture}[scale=1]
\draw(0,0) circle [x radius = 3, y radius = 2.5]; 
\draw[line width=2pt] (-1.5,-2.15) .. controls (-1.2,-2.3) and (-1,-1.4).. (0,-1.5);
\draw[line width = 2pt](0,-1.5) .. controls (1,-1.4)  and (1.2,-2.3) .. (1.5,-2.15);  
\draw[dashed](-2.7,-1.1) ..controls (-1,-4)  and (1,-4).. (2.7,-1.1); 
\draw[line width=2pt] (-2.7,-1.1) .. controls   (-2.4,-1.6)  and (-2.2,-1.8)     .. (-1.5,-2.15);
\draw[line width=2pt] (2.7,-1.1).. controls (2.4,-1.6) and (2.2,-1.8)..     (1.5,-2.15);

\draw(0.2,-1.1) node {$\gamma_0'$};
\draw(0.4,-2.7) node {$\gamma_0$};
\draw(0.1,-2) node {$B$};
\draw(0.5,-3.5) node {$\gamma_\varepsilon$};
\draw(2,2.3) node {$\gamma_1$};
\draw(3,-1.8) node {$\Gamma_0\setminus\gamma_0$};
\draw(-2.9,-1.8)node {$\Gamma_0\setminus\gamma_0$};

\filldraw[black] (-2.7,-1.1) circle (2pt)
                      (2.7,-1.1) circle (2pt)
                       (-1.5,-2.15) circle (2pt)
                       (1.5,-2.15) circle (2pt);

\end{tikzpicture}
\\
{\bf Fig 7.2.} 
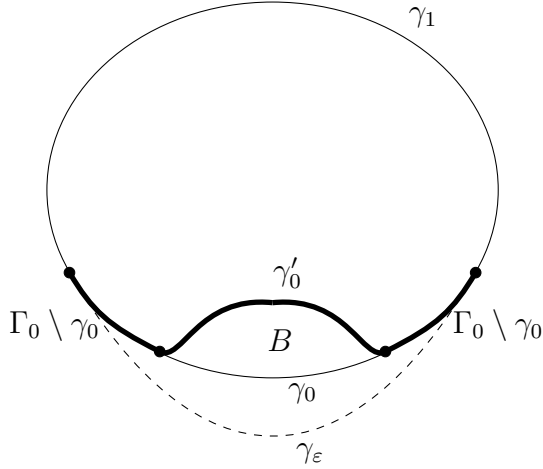
The boundary of $B$  is $\gamma_0\cup\gamma_0'$.  
\\
The boundary  of $D_\e$  is 
$\gamma_\e\cup\Gamma_0,\ \partial D_0^{(1)}=\Gamma_0\cup\gamma_1$.
\\

  To prove  Lemma \ref{lma:7.2}  we will need  the following version  of the Runge  theorem 
 about  the approximation  of solutions of the equation  in a smaller domain  by  solutions of the same equation  in a 
 larger domain.
 \begin{lemma}																				\label{lma:7.3}
  Denote by $D_{\e}$  the domain  bounded by $\Gamma_0$  and  $\gamma_\e$  such that 
$\gamma_\e\cup\gamma_1$  is smooth.
 Let $W_0$  be  the space of $v\in H_s((D_0^{(1)}\setminus B)\times(t_1,t_2)),\ s\geq 1$,  such that 
 \begin{align}																				\label{eq:7.9}
 &v\big|_{\gamma_1}=0,\ \ v\big|_{x_o=t_1}=\frac{\partial v}{\partial x_0}\Big|_{x_0=t_1}=0,
 \ \ x\in (D_0^{(1)}\setminus B),
 \\
 \nonumber
 &L^{(1)}v=0\ \ \mbox{in}\ \ (D_0^{(1)}\setminus B)\times(t_1,t_2),
 \end{align}
 where  $\gamma_1=\partial D_0^{(1)}\setminus\Gamma_0$.
 
 Denote by  $K$  the closure  of $W_0$  in $L_2((D_0^{(1)}\setminus B)\times (t_1,t_2))$.
 Consider  the space  $W$  of  $u(x)\in  H_s((D_0^{(1)}\cup D_\e)\times(t_1,t_2)),s\geq 1$  such that 
 \begin{align}																				\label{eq:7.10}
 &L^{(1)}u=0\ \ \mbox{in}\ \ (D_0^{(1)}\cup D_{\e})\times(t_1,t_2),\ u\big|_{(\gamma_1\cup\gamma_\e)\times(t_1,t_2)}=0,
 \\
 \nonumber
 &u\big|_{x_0=t_1}=\frac{\partial u}{\partial x_0}\Big|_{x_0=t_1} =0,\ \ x\in  D_0^{(1)}cup D_\e.
 \end{align}
 Then the closure of the restrictions  of $W$   to $L_2((D_0^{(1)}\setminus B)\times(t_1,t_2))$  is also equal  to $K$.
 Thus any function  $v\in W_0$  in  ($D_0^{(1)}\setminus B)\times (t_1,t_2)$ 
 can be  approximated in  $L_2((D_0^{(1)}\setminus B)\times (t_1,t_2))$  norm  by the functions in $W$.
 \end{lemma}
 {\bf Proof:}
 Let $K^\perp$  be the orthogonal complement
 of $K$  in  $L_2((D_0^{(1)}\setminus B)\times(t_1,t_2))$.  
 Take  any $g\in K^\perp$  and denote  by $g_0$  the extension  of $g$  by 
 zero outside $(D_0^{(1)}\setminus B)\times(t_1,t_2)$.  Let  $w$  be the solution  of the initial-boundary value
 problem
 \begin{align}																			\label{eq:7.11}
 &L_1^*w=g_0,\ x\in (D_0^{(1)}\cup D_{\e })\times(t_1,t_2),
 \\
 \nonumber
 &w\big|_{x_0=t_2}=\frac{\partial w}{\partial x_0}\Big|_{x_0=t_2}=0,\ \ x\in D_0^{(1)}\cup D_\e,
 \\
 \nonumber
 &w\big|_{\partial(D_{0}^{(1)}\cup D_{\e })\times(t_1,t_2)}=0,
 \end{align}
 where $L_1^*$  is the formally   adjoint  to $L^{(1)}$.
 Note that  $\partial(D_0^{(1)}\cup D_{\e 0})=\gamma_1\cup \gamma_\e$  (see Fig.7.2).
 
 By [H] and [E6] (see also  Lemma \ref{lma:3.3})  such $w(x)$  exists  and belongs to  
 $H_1((D_0^{(1)}\cup D_\e))\times(t_1,t_2))$.  We shall show  that  $w=0$ in  $(B\cup D_\e)\times(t_1,t_2)$.  Let
 $\varphi\in C_0^\infty((B\cup D_{\e})\times(t_1,t_2))$  and  let  $u(x)$  be the solution  of
 \begin{align}																			\label{eq:7.12}
 &L^{(1)}u=\varphi,\ \ x\in  (D_0^{(1)}\cup D_{\e })\times(t_1,t_2)
 \\
 \nonumber
 &u\big|_{x_0=t_1}=\frac{\partial u}{\partial x_0}\Big|_{x_0=t_1}=0,
 \ \ \ u\big|_{\partial(D_0^{(1)}\cup D_{\e })\times(t_1,t_2)}
 =0,
 \end{align}
 (cf.  [H], [E6] and Lemma \ref{lma:3.3}),  i.e.  $u\in  W_0$  since  $\varphi=0$  in  $(D_0^{(1)}\setminus B)\times(t_1,t_2)$.
 
 Consider the $L_2$  inner  product  $(\varphi,w)$  in  $(D_0^{1)}\cup D_{\e })\times (t_1,t_2)$.
 Since  $\varphi=L^{(1)}u$  we get  $(\varphi,w)=(L^{(1)}u,w)$.  Integrating  by parts  we have  $(L_1u,w)=
 (u,L_1^*w)=(u,g_0)=0$  since  $u\in W_0, g_0\in K^\perp$.  Therefore  $(\varphi,w)=0,\ \forall \varphi$.  Thus
 $w=0$  in $(B\cup D_{\e })\times(t_1,t_2)$.
 
 Let  now  $\tilde w$ be any function  in $W$.  We  have  $(\tilde w,g_0)_0=(\tilde w,L_1^*w)_0$,
 where  $(\ ,\,)_0$ means that  we integrate  over $(D_0^{(1)}\setminus B)\times(t_1,t_2)$.  Since
 $w=0$  in  $(B\cup D_{\e})\times(t_1,t_2)$,  we have that
 \begin{equation}         															\label{eq:7.13}
 w\big|_{(\Gamma_0\setminus\gamma_0)\cup \gamma_0')\times(t_1,t_2)}=
 \frac{\partial w}{\partial \nu}\Big|_{(\Gamma_0\setminus\gamma_0)\cup \gamma_0')\times(t_1,t_2)}=0,
 \end{equation}
 where $\frac{\partial}{\partial \nu}$  is the normal derivative.
 
 Note that $(\Gamma_0\setminus \gamma_0)\cup \gamma_0'=\partial(D_{\e }\cup B)\setminus\gamma_\e$.
 Since $w$  satisfies   the homogenous equation  $L_1^*w=0$  in $D_\e\cup B$  the  restrictions  of  
 $w$  and all  derivatives  on 
 $\partial(D_{\e}\cup B)$  exists  by the partial  hypoellipticity (see,  for example,  [E5]).
 Note that  $\tilde w$  and  $w$  have zero values  on  $\gamma_1$.  Therefore,  integrating  by parts,  we have
 $$
 (\tilde w,L_1^* w)_0=(L^{(1)}\tilde w, w)_0=0,
 $$
 since  $L^{(1)}\tilde w=0$  in  $(D_0^{(1)}\setminus B)\times(t_1,t_2)$.  Therefore $(\tilde w,g_0)_0=0,\ \ \forall g_0\in K^\perp,$
  \  i.e. $\tilde w\in \overline K$.
 
 To make the integration  by  parts rigorous  we  approximate  
 $\gamma_0'\cup(\Gamma_0\setminus\gamma_0)$  by  $\gamma_{\e_1}'$,  similar to  $\gamma_\e,\ \gamma_{\e_1}'
 \subset D_\e\cup B$.  Note that  $w=0$ in $D_\e\cup B$.  Therefore integrating
   by parts over domain bounded  by $\gamma_1\cup \gamma_{\e_1}'$,
   and  taking  the limit 
    when $\gamma_{e_1}'\rw\gamma_0'\cup(\Gamma_0\setminus\gamma_0)$ we get $(\tilde w,g_0)=0$.
 \qed
 
  Now we  shall proof  Lemma \ref{lma:7.2}.
 
 Let $\mbox{supp}\,f\subset \Gamma_0'\times(t_1,t_2),\ \Gamma_0'=(\Gamma_0\setminus\gamma_0)\cup\gamma_0'$.
 Let  $v_1$  be  the solutions of  
 \begin{align}																\label{eq:7.14}
 &L^{(1)}v_1=0, \ x\in (D_0^{(1)}\setminus B)\times(t_1,t_2),
 \\
 \nonumber
 &v_1\big|_{x_0=t_1}=\frac{\partial v_1}{\partial x_0}\Big|_{x_0=t_1}=0,
 \\
 \nonumber 
 &v_1\big|_{\partial (D_0^{(1)}\setminus B)\times(t_1,t_2)}=f_1,
 \end{align}
 where $\partial(D_0^{(1)}\setminus B)=\Gamma_0'\cup\gamma_1,\ f_1=0$  on $\gamma_1\times(t_1,t_2), \ f_1=f$
 on $\Gamma_0'\times(t_1,t_2)$.
 
 Let  $v_3$  be  solution  of  $L^{(3)}v_3=0$  in  $D_1^{(3)}\setminus(B\times(t_1,t_2))$
 \begin{align}												 						\label{eq:7.15}
 &v_3\big|_{\partial_-D_1^{(3)}}=\frac{\partial v_3}{\partial x_0}\Big|_{\partial_-D_1^{(3)}}=0,
 \ \ \ v_3\big|_{(\partial'D_1^{(3)}\setminus(\Gamma_0\times(t_1,t_2))}=0,
 \\
 \nonumber
&v_3\big|_{\Gamma_0'\times(t_1,t_2)}=f.
\end{align}
 
 Let  $\Lambda_1'$  be  the DN operator  for  (\ref{eq:7.14})  and $\Lambda_3'$  be the 
 DN  operator  for  (\ref{eq:7.15}).  Assuming  that
  $\Lambda_1=\Lambda_2$  on  $\Gamma_0\times(t_1,t_2)$  we shall prove  that 
 $$
 \Lambda_1'f\Big|_{\Gamma_0'\times(t_1+\delta,t_2-\delta)}=\Lambda_2' f\Big|_{\Gamma_0'\times(t_1+\delta,t_2-\delta)}
 $$
 for all $f$  with supports in $\Gamma_0'\times(t_1+\delta,t_2-\delta)$.  By Lemma \ref{lma:7.3}  there exists 
  a sequence  of smooth solutions  $w_{n1}\in W_0$  such that 
 \begin{equation}															
\nonumber
 \|v_1-w_{n1}\|_0\rw 0, \ \ n\rw \infty,
 \end{equation}
 where $\|v_1\|_0$  is  the norm in  $L_2((D_0^{(1)}\setminus B)\times(t_1,t_2))$.    
 Note that 
 $
 L^{(1)}w_{n1}=0$  in  $D_0^{(1)}\times(t_1,t_2),\ w_{n1}\big|_{\gamma_1\times(t_1,t_2)}=0,\ 
 w_{n1}\big|_{x_0=t_1}=\frac{\partial w_{n1}}{\partial x_0}\big|_{x_0=t_1}=0$,
 where $\gamma_1=\partial D_0^{(1)}\setminus \Gamma_0$.  Let  $f_n=w_{n1}\big|_{\Gamma_0\times(t_1,t_2)}$.
 Denote by $w_{n3}$  the solution  of
 \begin{align}	 														\label{eq:7.16}
 L^{(3)}w_{n3}=0\ \ \mbox{in}\ \ D_1^{(3)},\ \ &w_{n3}\big|_{\partial'D_1^{(3)}\setminus(\Gamma_0\times(t_1,t_2))}=0,\ \ 
 w_{n3}\big|_{\Gamma_0\times(t_1,t_2)}=f_n,
 \\
 \nonumber
 &w_{n3}\big|_{\partial_-D_1^{(3)}}=\frac{\partial w_{n3}}{\partial x_0}\Big|_{\partial_-D_1^{(3)}}=0.
 \end{align} 
 Since $\Lambda_1=\Lambda_2$  on  $\Gamma_0\times (t_1,t_2)$,  we have
 \begin{equation}															\label{eq:7.17}
 \frac{\partial w_{n1}}{\partial \nu}\Big|_{\Gamma_0\times(t_1,t_2)}=
 \frac{\partial w_{n3}}{\partial \nu}\Big|_{\Gamma_0\times(t_1,t_2)}
 \end{equation}
  Since $\gamma_o\subset \Gamma_0$,  the equality
 (\ref{eq:7.17})  implies
 $$
 w_{n1}\big|_{\gamma_0\times(t_1,t_2)}=w_{n3}\big|_{\gamma_0\times(t_1,t_2)},\ \ \ 
 \frac{\partial w_{n1}}{\partial \nu}\Big|_{\gamma_0\times(t_1,t_2)}=
 \frac{\partial w_{n3}}{\partial \nu}\Big|_{\gamma_0\times(t_1,t_2)}.
 $$
 We have  $L^{(1)}=L^{(3)}$  on  $B\times(t_1,t_2)$.    Using the uniqueness theorem  of  [RZ]  and  [T],  we get
 \begin{equation}															\label{eq:7.18}
 w_{n1}=w_{n3}  \ \ \mbox{in}\ \ B\times(t_1+\delta,t_2-\delta),
 \end{equation}
 where $\delta >0$  is determined  by  the metric and by the domaine $B$  (cf. Fig.7.2).   In particular,
 \begin{align}															\label{eq:7.19}
 &w_{n1}\big|_{\gamma_0'\times(t_1+\delta,t_2-\delta)}=w_{n3}\big|_{\gamma_0'\times(t_1+\delta,t_2-\delta)},
 \\
 \nonumber
 &\frac{\partial w_{n1}}{\partial \nu}\Big|_{\gamma_0'\times(t_1+\delta,t_2-\delta)}=
 \frac{\partial w_{n3}}{\partial \nu}\Big|_{\gamma_0'\times(t_1+\delta,t_2-\delta)}.
 \end{align}
 Therefore
 \begin{align}															\label{eq:7.20}
& w_{n1}\big|_{\Gamma_0'\times(t_1+\delta,t_2-\delta)}=w_{n3}\big|_{\Gamma_0'\times(t_1+\delta,t_2-\delta)},
 \\ 
 \nonumber
 &\frac{\partial w_{n1}}{\partial \nu}\Big|_{\Gamma_0'\times(t_1+\delta,t_2-\delta)}=
 \frac{\partial w_{n3}}{\partial \nu}\Big|_{\Gamma_0'\times(t_1+\delta,t_2-\delta)},
 \end{align}
 where $\Gamma_0'=(\Gamma_0\setminus\gamma_0)\cup\gamma_0'$,  i.e.  $\Lambda_1'f_n'=\Lambda_2'f_n'$
 on  $\Gamma_0'\times(t_1+\delta,t_2-\delta)$,  where  $f_n'=w_{n1}\big|_{\Gamma_0'\times(t_1+\delta,t_2-\delta)}
 =w_{n3}\big|_{\Gamma_0'\times(t_1+\delta,t_2-\delta)}$.  We have
 \begin{multline}															\label{eq:7.21}
\|f- f_n'\|_{-\frac{1}{2},\Gamma_0'\times(t_1+\delta,t_2-\delta)}
 =\|f-f_n'\|_{-\frac{1}{2},\partial( D_0^{(1)}\setminus B)\times(t_1+\delta,t_2-\delta)}
 \\
 \leq
 \|v_1-w_{n1}\|_{0,(D_0^{(1)}\setminus B)\times (t_1+\delta_1,t_2-\delta)},
 \end{multline}
 since
 $\partial(D_0^{(1)}\setminus B)=\Gamma_0'\cup\gamma_1$  and  $f=f_n'=0$  on  $\gamma_1\times(t_1+\delta,t_2-\delta)$.
 
 In  (\ref{eq:7.21})  we again use the partial hypoellipticity 
 property that  restrictions of  solutions of  $L^{(1)}u=0$  to  
 the noncharacteristic  boundary  exists  for any  Sobolev's space
 $H_s$  (cf. [E5]).  The  same is true  for all  normal derivatives  of  $u_1$,  and  the same  estimates hold
  as in the case of
 positive  $s>0$  (cf  [E5] and [E1]).
 
 By Lemma \ref{lma:3.3}  (see [H], [E6])  we have 
 \begin{multline}															\label{eq:7.22}
 \Big\|\frac{\partial v_1}{\partial\nu}
 -\frac{\partial w_{n1}}{\partial\nu}\Big\|_{-\frac{3}{2},\partial(D_0^{(1)}\setminus B)\times(t_1+\delta,t_2)-\delta}
 \leq
  \big\|f-f_n'\big\|_{-\frac{1}{2},\partial(D_0^{(1)}\setminus B)\times(t_1+\delta,t_2-\delta)}
 \end{multline}
 Analogously  we have
 \begin{multline}															\label{eq:7.23}
 \Big\|\frac{\partial v_3}{\partial\nu}
 -\frac{\partial w_{n3}}{\partial\nu}\Big\|_{-\frac{3}{2},\partial'(D_1^{(3)}\cap(t_1+\delta,t_2-\delta)\setminus (B\times(t_1+\delta,t_2-\delta))}
 \\
 \leq
  \big\|f-f_n'\big\|_{-\frac{1}{2},\partial'(D_1^{(3)}\cap(t_1+\delta,t_2-\delta)\setminus (B\times(t_1+\delta,t_2-\delta))}
 \end{multline}
 Note that
 
  $ \big\|f-f_n'\big\|_{-\frac{1}{2},\partial(D_0^{(1)}\setminus B)\times(t_1+\delta,t_2-\delta)}=
  \big\|f-f_n'\big\|_{-\frac{1}{2},\partial'(D_1^{(3)}\cap(t_1+\delta,t_2-\delta)\setminus (B\times(t_1+\delta,t_2-\delta))}$.

 Therefore,  taking  the  limit  as  $n\rw\infty$,  we get,  using (\ref{eq:7.20}),  that
 \begin{equation}																	\label{eq:7.24}
 \frac{\partial v_1}{\partial \nu}\Big|_{\Gamma_0'\times(t_1+\delta,t_2-\delta)}=
 \frac{\partial v_3}{\partial \nu}\Big|_{\Gamma_0'\times(t_1+\delta,t_2-\delta)}.
 \end{equation}
 Thus we proved that
 $$
 \Lambda_1'f=\Lambda_2'f\ \ \mbox{on}\ \ \Gamma_0'\times(t_1+\delta,t_2-\delta)
 $$
 for any  $f$  with $\mbox{supp}\,f\subset\Gamma_0'\times(t_1+\delta,t_2-\delta)$.
 \qed
 
Using  Lemma \ref{lma:7.2}  we reduce  the inverse problem  in $D_0^{(1)}\times(t_1,t_2)$  to the inverse  problem  in
smaller domains $(D_0^{(1)}\setminus B)\times(t_1+\delta,t_2-\delta)$  and  we can continue  this process
starting  from  
 $(D_0^{(1)}\setminus B)\times(t_1+\delta,t_2-\delta)$  instead of $D_0^{(1)}\times(t_1,t_2)$.  

In all lemmas below we assume that DN  operators  for  $L^{(1)}$  and  $L^{(2)}$  are equal  on
$\Gamma_0\times[t_1,t_2]$ and that the time  interval  $[t_1,t_2]$  is large enough.
We shall continue to call  the coordinates $(y_0,\hat y',y_n)$,  given  by the map $\beta^{(i)}\varphi^{(i)}$, the Goursat
coordinates  for  $\tilde L_1^{(i)},i=1,2$.
\begin{lemma}																	\label{lma:7.4}
Let  $\Gamma_1\subset \Gamma_0$  and  let  $\tilde \Gamma_1\subset \Gamma_0$  be such  that  $\overline \Gamma_1
\subset\tilde\Gamma_1$.  We assume  that  $\tilde\Gamma_1\subset\tilde\Gamma_2\subset\tilde\Gamma_3
\subset \Gamma_0$  where  $\tilde\Gamma_j,1\leq j\leq 3,$  are the same as  $\Gamma_j,1\leq j\leq 3,$  in \S 3. 
Suppose the Goursat coordinates  for $L^{(1)}$  exists  in $\Omega_1=(t_1,t_2)\times\tilde\Gamma_3\times[0,\e_0]$,  i.e.
$L^{(1)}$  has the form $\tilde L_1^{(1)}$ in these  coordinates (we include the gauge  transformation (\ref{eq:2.21})  in 
$\tilde L_1^{(i)})$.   Suppose the Goursat  coordinates  for  $L^{(2)}$  exist in 
$(t_1,t_2)\times(\overline {\tilde\Gamma}_3\setminus \Gamma_1)\times [0,\e_0]$.  Let  $\tilde\Omega_2=
(t_1,t_2)\times(\tilde\Gamma_1\setminus\Gamma_1)\times[0,\e_0]$
and  suppose  $\tilde L_1^{(2)}=\tilde L_1^{(1)}$  in  $\tilde\Omega_2$.  Then  $L^{(2)}$  has also Goursat  coordinates
in  $\Omega_3=(t_1,t_2)\times\Gamma_1\times[0,\e_0]$,   and $\tilde L_1^{(1)}=\tilde L_1^{(2)} $  in  $\Omega_3$.
\end{lemma}
 
 {\bf Proof:}
 Let  $y^{(i)}=\psi^{(i)}(x)$  be the transformation  to the Goursat  coordinates,  and
 let  $\frac{D\psi^{(i)}(x)}{Dx}$ be  the Jacobi matrix  of this transformation. We have
 $$
 [\hat g_i^{jk}(y)]=\frac{D\psi^{(i)}(x)}{D(x)}[g_i^{jk}(x)]\Big(\frac{D\psi^{(i)}}{Dx}\Big)^{T},\ i=1,2,
 $$
 where $[\hat g_i^{jk}]^{-1}$  is the metric  tensor  in the  Goursat coordinates.  The Goursat coordinates degenerate 
 at point  $y_i^{(0)}=\psi^{(i)}(x^{(0)})$,  when  $\det\frac{D\psi^{(i)}(x)}{Dx}\rw \infty$ for
 $y\rw y^{(0)}$  (or $x\rw x^{(0)}$)  \ (cf. (\ref{eq:2.11})).  We call  such point  a focal point.
 We shall prove that  there is no focal points for  $L^{(2)}$  in $\Omega_3$.
 
 We have
  $$
\det [\hat g_i^{jk}(y)]=\det [g_i^{jk}(x)]\Big(\det\frac{D\psi^{(i)}}{Dx}\Big)^{2}.
 $$
 Suppose  there exists  the focal  point  $y^{(0)}=(y_0^{(0)},y_0',y_n^{(0)}),
 y_n^{(0)}<\e_0,y_0'\in\Gamma_1$  such  that  there  is no focal  points  for $L^{(2)}$  when  $y_n<y_n^{(0)}$
 for all  $y_0\in [t_1,t_2], y'\in\overline\Gamma_1.$
 
 Since  $L^{(1)}$  and $L^{(2)}$  have  Goursat  coordinates  for  $y_n<y_n^{(0)}$  we  get,  by Lemma \ref{lma:7.1},
 that $\tilde L_1^{(1)}=\tilde L_1^{(2)}$ 
in $(t_1,t_2)\times\overline\Gamma_1\times[0,y_n^{(0)}-\e],  \forall\e >0$,  
 and hence  $[\hat g_1^{jk}]=[\hat g_2^{jk}]$
 for $\e>0$.
 Since  $\det[\hat g_2^{jk}]=\det[\hat g_1^{jk}]$  for $y_n<y_n^{(0)}-\e$,  we have that 
 $\big(\det \frac{D\psi^{(2)}}{Dx}\big)^2=\frac{\det[\hat g_1^{jk}]}{\det[ g_2^{jk}]}$
 is bounded when  $\e\rw 0$.  Therefore  $y^{(0)}=(y_0^{(0)},y_0',y_n^{(0)})$  is  not a focal  point for  $L^{(2)}$.
 Thus  $L^{(2)}$  has   no focal points  in $\Omega_3$  and  then,  by Lemma \ref{lma:7.1},  we have  
 $\tilde L_1^{(1)}=\tilde L_1^{(2)}$ in  $\Omega_3$  (cf. [E2]).
 \begin{lemma}																		\label{lma:7.5}
 Assume  that DN operators  for  $L^{(1)}$  and  $L^{(2)}$ are equal  on  $\Gamma_0\times[t_1,t_2]$.
 Let  $\overline\Gamma_1\subset\Gamma_0$.
 Assume  that  the Goursat coordinates for  $L^{(1)}$  exists  in  
 $(t_1,t_2)\times\overline \Gamma_1\times[0,\frac{T_0}{2}]$.   Then  the Goursat coordinates  for $L^{(2)}$  also
 exists in $\Omega_1=(t_1+\delta,t_2-\delta)\times\overline\Gamma_1\times[0,\frac{T_0}{2}]$  for some  $\delta>0$  
 and $\tilde L_1^{(1)}=\tilde L_1^{(2)}$  in  $\overline\Omega_1$,  where  $\tilde L_1^{(i)}$  are  the operators  $L^{(i)}$
 in  the Goursat coordinates.
 \end{lemma}
 {\bf Proof:}  
 Let  $\overline\Gamma_1\subset\tilde\Gamma_1,\  \overline{\tilde\Gamma}_1\subset \Gamma_0$.
 If  $0\leq y_n\leq \e$,  where $\e>0$  is small  enough,  then  
 $\tilde\Gamma_1\subset\tilde\Gamma_2\subset \tilde\Gamma_3\subset\Gamma_0$,  where
 $\tilde\Gamma_j,j=1,2,3,$  are  the same as in Lemma \ref{lma:7.4}.   Applying Lemma \ref{lma:7.1}  
 we get  that the Goursat  coordinates  for 
 $\tilde L_1^{(1)}$ and $\tilde L_1^{(2)}$
 exist  in $\Omega_{1\e}=[t_1,t_2]\times\tilde\Gamma_1\times[0,\e]$  and 
 $$
 \tilde L_1^{(1)}=\tilde L_1^{(2)} \ \ \mbox{in}\ \ \overline\Omega_{1\e}.
 $$
 Let  $\Sigma_1$  be the surface  in $(y',y_n)$ space
 such that 
 $y_n=0$  on $\tilde \Gamma_1\setminus\Gamma_1,0\leq y_n\leq \e$  on  $\partial\Gamma_1,y_n=\e$  on  $\Gamma_1$.
 Note that $\Sigma_1$   is not  smooth  since  it has edges
 when  $y_n=0,y'\in\partial\Gamma_1$  and when  $y_n=\e,y'\in\partial\Gamma_1$.  We shall smooth  $\Sigma_1$
 by replacing it by  smooth surface  $\tilde\Sigma_1$,  where $\tilde\Sigma_1$  differs  from  $\Sigma_1$  in 
 a neighborhood  of edges  having  the size $O(\e)$.
 Let $\Sigma_2$ be the surface,  where $y_n=\e$  when $y'\in \tilde\Gamma_1\setminus
 \Gamma_1(\e),\ \Gamma_1(\e)$  is the  $\e$-neighborhood  of  $\Gamma_1,\e\leq y_n\leq 2\e, $
 when $y'\in\partial\Gamma_1(\e),\ y_n=2\e$   when  $y'\in\Gamma_1(\e)$  (cf.  Fig. 7.3).
\\
\\

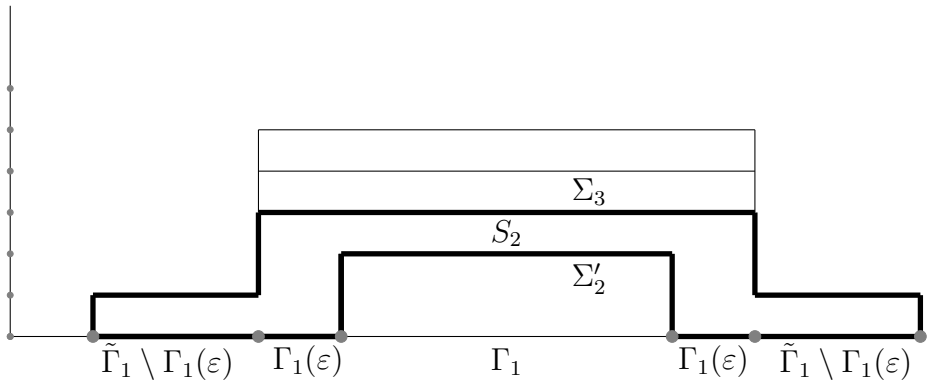
\begin{tikzpicture}[scale = 1.1]
\draw(0,0)--(11,0);
\draw(0,0)--(0,4);
\draw[line width = 2pt](1,0)--(1,0.5);
\draw[line width = 2pt](3,0.5)--(3,1.5);
\draw[line width = 2pt] (1,0.5)--(3,0.5);
\draw[line width = 2pt](1,0)--(4,0);
\draw[line width = 2pt] (4,0)--(4,1);
\draw[line width = 2pt](4,1)--(8,1);
\draw[line width = 2pt] (8,0)--(11,0);

\draw[line width = 2pt](8,0)--(8,1);
\draw[line width = 2pt](3,1.5)--(9,1.5);
\draw[line width = 2pt] (9,0.5)--(9,1.5);

\draw[line width = 2pt](1,0)--(4,0);
\draw[line width = 2pt] (4,0)--(4,1);

\draw[line width = 2pt](9,0.5)--(11,0.5);
\draw[line width = 2pt] (11,0.5)--(11,0);

\draw(3,1.5) -- (3,2.5);
\draw(3,2)--(9,2);
\draw(3,2.5) -- (9,2.5);
\draw(9,1.5) -- (9,2.5);

\filldraw[gray] (0,0.5) circle (1pt)
                      (0,1) circle (1pt)
                       (0,1.5) circle (1pt)
                       (0,2) circle (1pt);
                       (0,2.5) circle (1pt);
   \filldraw[gray] (0,3) circle (1pt);
                       (0,3.5) circle (1pt);
    \filldraw[gray]  (0,2.5) circle (1pt);
 \filldraw[gray]  (0,0) circle (1pt);           
 
  \filldraw[gray]  (0,0) circle (1pt); 
     \filldraw[gray]  (3,0) circle (2pt);   
      \filldraw[gray]  (1,0) circle (2pt);   
       \filldraw[gray]  (4,0) circle (2pt);   
        \filldraw[gray]  (9,0) circle (2pt); 
            \filldraw[gray]  (8,0) circle (2pt);   
             \filldraw[gray]  (11,0) circle (2pt);        
    
    \draw(6,1.25) node {$S_2$};         
      \draw(7,0.7) node {$\Sigma_2'$};       
      \draw(7,1.75) node {$\Sigma_3$};  
\draw(1.9,-0.3) node {$\tilde\Gamma_1\setminus \Gamma_1(\e)$};
\draw(3.6,-0.3)  node {$\Gamma_1(\e)$};

\draw(10.1,-0.3) node {$\tilde\Gamma_1\setminus \Gamma_1(\e)$};
\draw(8.5,-0.3)  node {$\Gamma_1(\e)$};
\draw(6,-0.3) node {$\Gamma_1$};

\end{tikzpicture}
  \\
 \\
{\bf Fig. 7.3.} 
\\
$\Sigma_2'$  is the surface  $\{y_n=0\ \mbox{when}\ y'\in \tilde\Gamma'\setminus\Gamma_1,\ 
 0\leq y_n\leq 2\e
\ \mbox{when}\ y'\in\partial \Gamma_1,\linebreak
 y_n=2\e \ \mbox{when}\ y'\in \Gamma_1\}$,  
\\
$\Sigma_3$  is the surface
 $\{y_n=\e \ \mbox{when}\  y'\in \tilde\Gamma_1\setminus\Gamma_1(\e), \ 
  \e\leq y_n\leq 3\e
\ \mbox{when}\ y'\in\partial \Gamma_1(\e),\linebreak y_n=3\e \ \mbox{when}\ y'\in \Gamma_1(\e)\}$,  
\\
 $S_2$  is the region between $\Sigma_3$  and $\Sigma_2'$.
\\
\\
\

 Let  $\tilde\Sigma_2$  be the  smoothing  of $\Sigma_2$.  Denote  by  $\tilde S_1$  the  domain  between  $\tilde\Sigma_1$  and $\tilde\Sigma_2$  when  $y'\in\tilde\Gamma_1$.
  Since  $\tilde L_1^{(1)}=\tilde L_1^{(2)}$ for  $0\leq y_n\leq \e$,  we have,  by Lemma \ref{lma:7.2},  that  DN operators  for  
  $L^{(1)}$  and $L^{(2)}$  are equal   on $\tilde\Sigma_1\times(t_1+\delta_1,t_2-\delta_1)$  for some $\delta_1>0$.
 
 Suppose $\e>0$  is small and such that 
 we can introduce the
 Goursat coordinates for $L^{(1)}$   in 
 $\tilde S_1\times[t_1+\delta_1,t_2-\delta_1]$.
   
Note that $\e$  and $\delta_1$ are determined by $L^{(1)}$  only  and are independent  of $L^{(2)}$.
 It follows from Lemma \ref{lma:7.4}  that Goursat coordinates  for  $L^{(2)}$  also hold  on 
 $\tilde S_1\times(t_1+\delta_1,t_2-\delta_1)$
 and
 $\tilde L_1^{(1)}=\tilde L_1^{(2)}$  in  $\tilde S_1\times(t_1+\delta_1,t_2-\delta_1)$.
 
 By  Lemma \ref{lma:7.2} DN  operators  for  $L^{(1)}$  and  $L^{(2)}$  are  equal  on  
 $\tilde\Sigma_2\times(t_2+\delta_1,t_2-\delta_1)$.
 
 Let  $\Sigma_2'$  be  the surface in $(y',y_n)$  space 
   such  that  $y_n=0$  on  $\tilde\Gamma_1\setminus\Gamma_1,\ 
 0\leq y_n\leq 2\e$  on  $\partial\Gamma_1,\ y_n=2\e$  on $\Gamma_1$   and let  $\Sigma_3$  be the surface
  where  $y_n=\e$  on  $\tilde\Gamma_1\setminus\Gamma_1(\e),\ \e\leq y_n\leq 3\e$  on 
 $\partial\Gamma_1(\e),y_n=3\e$  on  $\Gamma_1(\e)$.  Let  $\tilde\Sigma_2'$  and $\tilde\Sigma_3$  be the 
 smoothing  of  $\Sigma_2',\Sigma_3$  and let $\tilde S_2$  be the domain  between 
 $\tilde\Sigma_2'$  and $\tilde\Sigma_3$ when  $y'\in \tilde\Gamma_1$.  Since  DN  operators  for  $L^{(1)}$  and
 $L^{(2)}$ are equal  on  $\tilde\Sigma_2'\times(t_1+\delta_1,t_2-\delta_1)$  and since the Goursat coordinates for 
 $L^{(1)}$  hold  on  $\tilde S_2\times(t_1+\delta_1,t_2-\delta_1)$,   Lemma \ref{lma:7.4}
 implies  that the Goursat coordinates  hold
 for  $L^{(2)}$  in   $S_2\times(t_1+\delta_1,t_2-\delta_1)$,
 $L_1^{(1)}=L_1^{(2)}$  in  $S_2\times(t_1+\delta_1+\delta_2,t_2-\delta_1-\delta_2)$  
 for some $\delta_2>0$,
 and DN  operators  for $L^{(1)}$ 
 and $L^{(2)}$  are equal on $\tilde\Sigma_3\times(t_1+\delta_1+\delta_2,t_2-\delta_1-\delta_2)$.
 
 Analogously,  for $k>2$  denote by  $\Sigma_k'$  the surface 
   such  that  $y_n=0$  on  $\tilde\Gamma_1\setminus\Gamma_1, 0\leq y_n\leq k\e$  on 
    $\partial\Gamma_1$  and $y_n=k\e$  on  $\Gamma_1$.  
 Let $\Sigma_{k+1}$  be
 the surface,  where $y_n=\e$  for  $\tilde\Gamma_1\setminus \Gamma_1(\e),\ \e\leq y_n\leq 
 (k+1)\e$ on $\partial\Gamma_1(\e)$ and  $y_n=(k+1)\e$  on  $\Gamma_1(\e)$.
 
 Denote by $\tilde\Sigma_k'$  and  $\tilde\Sigma_{k+1}$  the smoothing  of  $\Sigma_k',\Sigma_{k+1}$.
 Let  $\tilde S_k$  be the domain between $\tilde\Sigma_k'$ and  $\tilde \Sigma_{k+1}$.  Applying successively  
 the same arguments for
 $k=3,...,m$,  we prove  as above that $\tilde L_1^{(2)}=\tilde L_1^{(1)}$  in 
 $\tilde S_k\times[t_1+\sum_{j=1}^k\delta_j,t_2-\sum_{j=1}^k\delta_j],k=3,...,m$.
 
 Let $m$  be such  that $(m+1)\e\geq\frac{T_0}{2}$.  Then we get that $\tilde L_1^{(2)}=\tilde L_1^{(2)}$ 
 in $(t_1-\delta,t_2+\delta)\times\overline\Gamma_1\times[0,\frac{T_0}{2}]$,  where $\delta=\sum_{j=1}^m\delta_j$.
 Note that we assume  that  $[t_1,t_2]$  is large.  Thus  $t_2-t_1\gg \delta$.
 \qed

  Suppose  that after several  applications of Lemma \ref{lma:7.2}  we have
 $$
 L^{(1)}u_1=0\ \ \mbox{in }\ \ D_0^{(1)}\times[t_1,t_2],
 $$
 $$
 L^{(m)}u_2=0\ \ \mbox{in }\ \ D_1^{(m)},
 $$
 where we are  considering  the interval  $[t_1,t_2]$ instead  of $[t_1+\delta,t_2-\delta]$  for the  simplicity  of notations.
 We assume that  
 $$
 u_1\big|_{x_0=t_1}=
 \frac{\partial u_1}{\partial x_0}\Big|_{x_0=t_1}=0,
\ \ \
  u_2\big|_{\partial_-D_1^{(m)}}=
 \frac{\partial u_2}{\partial x_0}\Big|_{\partial_- D_1^{(m)}}=0.
$$
 We also assume  that  $\Gamma_0\times(t_1,t_2)\subset\partial' D_1^{(m)}\cap(\partial D_0^{(1)}\times(t_1,t_2))$
 and $\Omega_1\times(t_1,t_2)\subset(D_0^{(1)}\times(t_1,t_2))\cap D_1^{(m)}$  (cf.  Fig.  7.4).
 
 We assume  that $L^{(1)}=L^{(m)}$  in  $\Omega_1\times(t_1,t_2)$  and that DN  operators $\Lambda_1$
 and  $\Lambda_2$  are 
 equal  on  $\partial\Omega_1\times(t_1,t_2)$  and  $\Gamma_0\times (t_1,t_2)$.  Here,  as above,  $\partial'D_1^{(m)}$
 means  the time-like  part  of the boundary of  $D_1^{(m)}$.
 It follows from Lemmas \ref{lma:7.4}, \ref{lma:7.5}  that  the enlargement  of the domain  $\Omega_1$  depends   only on  $L^{(1)}$
 and does not depend  on  $L^{(2)}$.   Therefore,  as in  [E2],  we arrive  to the situation  when  $\Omega_1$  and  $D_1^{(0)}$  are
 close.   To apply   Lemma  \ref{lma:7.2}  to the domain   $(D_1^{(0)}\setminus\Omega_1)\times\R$  we need  new tools. 
 
 When $\partial\Omega_1$  and  $\partial D_0^{(1)}$  are close,   there is a narrow  domain 
 $\sigma_1\subset D_0^{(1)}\setminus\Omega_1$ 
  such  that $\gamma_1\subset\partial D_0^{(1)},\gamma_0\subset\partial\Omega_1$  and  the 
 distance between  $\gamma_0$  and $\gamma_1$  is small (cf.  Fig. 7.4)
 \\
\\
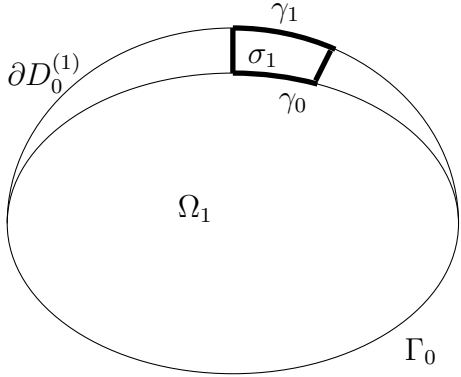
\begin{tikzpicture}[scale=1]
\draw (0,0) circle [x radius =3, y radius = 2];
\draw (3,0) arc (0:180: 3 and 2.6);
\draw [line width =2pt](0,2) -- (0,2.6);
\draw [line width = 2pt] (1.1,1.9) -- (1.3, 2.3);
\draw [line width = 2 pt] (0,2) arc (90:68:3 and 2);
\draw [line width = 2pt] (0,2.6) arc (90:63: 3 and 2.6);

\draw(-0.5,0.2) node {$ \Large\Omega_1$};
\draw(-2.5, 2) node {$\partial D_0^{(1)}$};
\draw(0.8,1.6) node {$\gamma_0$};
\draw(0.4,2.2) node {$\sigma_1$};
\draw(0.7,2.8) node {$\gamma_1$};
\draw(2.5,-1.7) node {$\Gamma_0$};
\end{tikzpicture}
 \\
 \\
{\bf Fig. 7.4.}  $\gamma_1\subset \partial D_0^{(1)},\gamma_0\subset\Omega_1,\ \gamma_1$  and  $\gamma_0$
 are close. 
 \\
 \
 
 Introduce Goursat coordinates for $L^{(1)}$  and $L^{(m)}$ near $\gamma_0\times(t_1,t_2)$.  We assume that
 operators $L^{(1)}$ and $L^{(m)}$  are defined  in domains slightly larger than  $D_0^{(1)}\times(t_1,t_2)$  and 
 $D_1^{(m)}$.  Let $L_1^{(1)}$  and $L_1^{(m)}$  be the operators $L^{(1)}$  and $L^{(m)}$  in corresponding Goursat 
 coordinates.  Let $y=\varphi_1(x)$  be the transformation to the Goursat coordnates  for $L^{(1)}$.
 Let $\sigma_0=(t_1,t_2)\times\gamma_0\times[0,\e_0]$  be  the domain where the Goursat coordinates for $L^{(1)}$  hold.
 We assume that $\sigma_1$  is so small  that  $\varphi_1(\sigma_1\times(t_1,t_2))\subset\sigma_0$.  Let
 $\tau_0=\varphi_1(\gamma_1\times(t_1,t_2))$,  i.e.  $\tau_0$  is the image  of part of the boundary  
 $\partial D_0^{(1)}\times(t_1,t_2)$  in  
 Goursat  coordinates.  Denote  by  $\sigma_0^+$  the part of $\sigma_0$  between 
 $\gamma_0\times(t_1,t_2)$   and $\tau_0$,  i.e.  $\sigma_0^+$  is  the image of  $\sigma_1\times(t_1,t_2)$  
 in Goursat  coordinates.
 
 We assume  that the Goursat coordinates for $L^{(m)}$  also hold  in $(t_1,t_2)\times\gamma_0\times[0,\e_0]$.  
 Moreover,  applying
 Lemmas \ref{lma:7.4},  \ref{lma:7.5}
 repeatedly we get that 
  $\tilde L_1^{(m)}=\tilde L_1^{(1)}$
 in $\hat\sigma_0^+$,  where $\hat\sigma_0^+=\sigma_0^+\cap(t_1+\delta,t_2-\delta)$.
 Here  $\tilde L_1^{(1)},\tilde L_1^{(m)}$ are operators  $L^{(1)},L^{(m)}$  in Goursat coordinates.
 
 Consider the initial-boundary value problem  for $\tilde L_1^{(1)}$  in Goursat coordinates
 \begin{align}																		\label{eq:7.25}
 &\tilde L_1^{(1)}\tilde u_1=0\ \ \mbox{in}\ \ \hat\sigma_0^+,
 \\
 \nonumber
 &\tilde u_1\big|_{x_0=t_1+\delta}=\frac{\partial\tilde u_1}{\partial x_0}\Big|_{x_0=t_1+\delta}=0,
 \\
 \nonumber
 &\tilde u_1\big|_{(t_1+\delta,t_2-\delta)\times\gamma_0}=f,\ \ \ \tilde u_1\big|_{\partial \hat\sigma_0^+}=0,
 \end{align}
 where $\mbox{supp}\,f\subset(t_1+\delta,t_2-\delta)\times\gamma_0$.
 
 Consider  now  the equation $\tilde L_1^{(m)}\tilde u_2=0$  in $\hat\sigma_0^+.$
 \begin{align}	
 \nonumber																	
 &\tilde L_1^{(m)}\tilde u_2=0\ \ \mbox{in}\ \ \hat\sigma_0^+,
 \\
 \nonumber
 &\tilde u_2\big|_{x_0=t_1+\delta}=\frac{\partial\tilde u_2}{\partial x_0}\Big|_{x_0=t_1+\delta}=0,
 \\
 \nonumber
 &\tilde u_2\big|_{(t_1+\delta,t_2-\delta)\times\gamma_0}=f,
 \end{align}
 where $f$ is the same as in (\ref{eq:7.25}).
 
 Since $\tilde L_1^{(1)}=\tilde L_1^{(m)}$  in $\hat\sigma_0^+$ and  since  
 $\Lambda_1^{(1)}f =\Lambda_1^{(2)}f$  on  $(t_1+\delta,t_2-\delta)\time\gamma_0$  we  have,  by  the unique  continuation  theorem  
 of [RZ]  and  [T]  
  that 
 $\tilde u_1=\tilde u_2$  in  $\hat\sigma_1^+\cap (t_1+2\delta,t_2-2\delta)$.
 Therefore  by the continuity
 $\tilde u_2\big|_{\partial \hat \sigma_1^+\cap (t_1+2\delta,t_2-2\delta)}=0$.
 
 Let $\sigma_2^+=\varphi_2^{-1}(\hat\sigma_0^+),\ \tau_2=\varphi_2^{-1}(\partial\hat\sigma_0^+)$,
 where
 $y=\varphi_2(x)$  is the transformation to the Goursat coordinates for  $L^{(m)}$.
 
 We shall show  that $\tau_2$  is a part  of  the   
 boundary of  $D_1^{(m)}$.
 
 Construct  the geometric optic solution  $v_1(y)$  for $\tilde L_1^{(1)}u_1=0$
 in Goursat coordinates  as in  (\ref{eq:5.1}).  Since $\tau_0$  is the boundary of the domain  $\sigma_0^+$  and 
 since the zero Dirichlet boundary condition holds on $\hat\tau_0=\tau_0\cap(t_1+\delta,t_2-\delta)$  
 this solution must reflect  at  $\hat\tau_0$  (cf.  [E7]).
 
 Consider now the geometric optics solution  $v_2(y)$  for $\tilde L_1^{(m)}u_2=0$ with
 the same initial condition.  Since $\tilde L_1^{(1)}v_1=\tilde L_1^{(m)}v_2$
 in $\hat\sigma_0^+$ we have that  $v_1(y)=v_2(y)$  before the reflection at $\hat\tau_0$.   If
 $\tau_2=\varphi_2^{-1}(\hat\tau_0)$  is not a part  of the boundary  of $\partial D_1^{(m)}$  there will be no
 reflection  for $v_2(y)$  at  $\hat\tau_0$.   Thus,  the solutions $v_1(y)$  and $v_2(y)$  will be different  in
 $\hat\sigma_0^+$  near $\hat\tau_0$.  This  contradicts   the fact that $v_1(y)=v_2(y)$  in  $\hat\sigma_0^+$.
 
 Therefore 
  $\varphi_m=\varphi_2^{-1}\varphi_1$  maps the boundary  
 $\gamma_1\times(t_1,t_2)$  of  $\partial D_0^{(1)}\times(t_1,t_2)$
 on the part of boundary of $\partial D_1^{(m)}$.
 Let  $\sigma_2\subset D_1^{(m)}$  be the image  
 of $\varphi_m=\varphi_2^{-1}\varphi_1(\sigma_1\times(t_1+2\delta,t_2-2\delta))$.  
 Let $\partial_\pm((\overline\Omega_1\times[t_1+2\delta,t_2-2\delta])\cup\overline\sigma_2)$  
 be the space-like parts  of the boundary  of 
  $(\overline\Omega_1\times[t_1+2\delta,t_2-2\delta])\cup\overline\sigma_2$. 
   Extend $\partial_+((\overline\Omega_1\times[t_1+2\delta,t_2-2\delta])\cup\overline\sigma_2)$ and
  $\partial_-((\overline\Omega_1\times[t_1+2\delta,t_2-2\delta])\cup\overline\sigma_2)$ 
  as space-like surfaces  $S_m^+$  and  $S_m^-$  to the whole domain
  $D_1^{(m)}$.  Let  $\tilde D_1^{(m)}$  be  part  of  $D_1^{(m)}$  bounded  from below and above  by $S_m^-$  
  and  $S_m^+$,
  respectively. 
  Note
 that $\partial' D_1^{(m)}\supset \Gamma_0\times(t_1+2\delta,t_2-2\delta)$.
 Define  $\tilde\varphi_m=I$  on  $\Omega_1\times (t_1+2\delta,t_2-2\delta),
 \ \tilde\varphi_m=I$  on $\Gamma_0\times(t_1+2\delta,t_2-2\delta),\ 
 \tilde\varphi_m=\varphi_m$  on  $\sigma_1\times(t_1+2\delta,t_2-2\delta)$.  Let $\Phi_{m+1}$  be the  extension  of 
 $\tilde\varphi_m$  (cf.  [Hi])  to  the whole  domain  $\tilde D_1^{(m)}$.  Denote  by  $D_1^{(m+1)}$  the image  of  
 $\tilde D_1^{(m+1)}$ under  the map  
 $\Phi_{m+1}$.  If $c_m$  is a gauge transformation  on $\overline\sigma_1\times[t_1+2\delta,t_2-2\delta]$  we denote by
  $\tilde c_{m+1}$  
 the extension  of  $c_m$  to  $D_1^{(m+1)}$  such  that  $\tilde c_{m+1}=1$  on  
 $\Omega_1\times[t_1+2\delta,t_2-2\delta],\ \tilde c_{m+1}=1$ 
  on 
 $\Gamma_0\times(t_1+2\delta,t_2-2\delta)$.
 
 We just proved the following lemma:
 \begin{lemma}           																	\label{lma:7.6}
 Let  $L^{(1)}$  and  $L^{(m+1)}=\tilde c_{m+1}\circ \Phi_{m+1}\circ L^{(m)}$  be operators  in  
 $D_0^{(1)}\times[t_1+2\delta,t_2-2\delta]$  and  
 $D_1^{(m+1)}$,  respectively.
 Then  $\Omega_2\times(t_1+2\delta,t_2-2\delta)\subset D_1^{(m+1)}\cap (D_0^{(1)}\times(t_1+2\delta,t_2-2\delta))$  where
 $\overline\Omega_2=\overline\Omega_1\cup\overline\sigma_1$ 
  and $L^{(1)}=L^{(m+1)}$  in  $\Omega_2\times(t_1+2\delta,t_2-2\delta)$.
 \end{lemma}
 We shall proceed 
 with  the enlargement  of the domain 
  $\Omega_2$ using  Lemmas \ref{lma:7.2} and \ref{lma:7.6}.    Therefore after
  finite number of steps  (cf. [E2]) 
  we get a domain  $D_1^{(N)}$,  operator $L^{(N)}$  on  $D_1^{(N)}$  and  the  map  $\Phi_N$  of  $D_1^{(N)}$
  onto  $D_0^{(1)}\times(t_1+\delta_N,t_2-\delta_N)$  such  that
  $c_N\circ\Phi_N\circ L^{(N)}=L^{(1)}$  in  $D_0^{(1)}\times(t_1+\delta_N,t_2-\delta_N)$  for some  $\delta_N>0$.
  Here  $c_N$  is  the gauge  transformation.  Remind  that $\tilde \Phi_2$  is the diffeomorphism  of  $ D_1^{(2)}$
  onto  $D_1^{(3)},\ \Phi_3$  is the diffeomorphism  of $\tilde D_1^{(3)}\subset D_1^{(3)}$  onto 
  $D_1^{(4)}$,  etc. ... $\tilde\Phi_{N-1}$ 
   is the map   of  $\tilde D_1^{(N-1})\subset D_1^{(N-1)}$  onto  $D_1^{(N)}$  and  $\Phi_N$  is  the map of
  $D_1^{(N)}$ onto  $D_0^{(1)}\times(t_1+\delta_N,t_2-\delta_N)$.
  
  Therefore,  the diffeomorphism  $\Phi^{-1}=\Phi_1^{-1}\Phi_3^{-1}...\Phi_N^{-1}$  maps  
  $D_0^{(1)}\times[t_0+\delta_N,t_2-\delta_N]$  onto  
  $D_1^{(2)}$.  Thus  $\Phi$  maps  $D_1^{(2)}$  onto
    $D_0^{(1)}\times[t_1+\delta_N,t_2-\delta_N]$.
  
  Note  that  $D_1^{(2)}$  is an almost  cylindrical  domain  in  $D_0^{(2)}\times\R$,  i.e.  $D_1^{(2)}=
  D_0^{(2)}\times\{S^-(x_1,...,x_n)\leq x_0\leq S^+(x_1,...,x_n)\}$,  where  $x_0=S^\pm(x_1,...,x_n)$   are
  space-like surfaces,  $(x_1,...,x_n)\in D_0^{(2)}$.
  
     Note that  $[t_1,t_2]$  is arbitrary  large and  therefore  
 $[t_1',t_2']=[t_1+\delta,t_2-\delta]$  is also   arbitrary  large.  
 Therefore we obtained  the following theorem:
 \begin{theorem}																	\label{theo:7.7}
 Let  $L^{(1)}$  and  $L^{(2)}$  be two  operators  in  $D_0^{(1)}\times\R$  and  $D_0^{(2)}\times\R$,  respectively.
 Suppose  $\Gamma_0\subset\partial D_0^{(1)}\cap\partial D_0^{(2)}$  and the DN  operators,
 corresponding to $L^{(i)}$,
  are equal  on
 $\Gamma_0\times\R$  for all $f$  that have  a compact support  in $\overline\Gamma_0\times\R$.
 Suppose  that the conditions  (\ref{eq:1.2}), (\ref{eq:1.6}) hold for $L^{(i)}, i=1,2$,  and the  coefficients  of $L^{(1)}$  and
 $L^{(2)}$  are analytic  in $x_0$  in  $D_0^{(i)}\times\R,\ i=1,2$. 
 Suppose  for each  $t_0\in R$  there exists  $T_{t_0}$
 such that the BLR  condition is satisfied  for $L^{(1)}$  on  $[t_0,T_{t_0}]$.
 Let  $[t_1',t_2']$  be  an arbitrary sufficiently large time  interval .  Then there exists  a diffeomorphism  $\Phi^{-1}$  of
 $\overline D_0^{(1)}\times[t_1',t_2']$  on an almost  cylindrical  domain
  $\overline D_1^{(2)}\subset\overline D_0^{(2)}\times \R,\ \Phi=I$  on $\Gamma_0\times[t_1',t_2']$
  and  there exists a gauge  transformation  $c(y)$  on  $D_1^{(2)}, \ |c(y)|=1$ 
 on  $ D_1^{(2)},\ c(y)=1$  on  $\Gamma_0\times[t_1',t_2']$  such that 
 $$
 c\circ \Phi^{-1}\circ L^{(2)}=L^{(1)}\ \ \mbox{on}\ \ D_0^{(1)}\times[t_1',t_2'].
 $$
 \end{theorem}
 Now we shall use Theorem  \ref{theo:7.7}  to prove  Theorem \ref{theo:1.2}.
 
 {\bf Proof of Theorem \ref{theo:1.2}}
 Let $L^{(i)}$  be two operators  in $D_0^{(i)}\times\R,i=1,2,\ \Gamma_0\subset \partial D_0^{(1)}\cap\partial D_0^{(2)}$
 and all  conditions  of Theorem \ref{theo:7.7}  are satisfied.
 
 Let  $(t_{j1},t_{j2})$  be  an interval  as in Theorem \ref{theo:7.7}  
 and  $\bigcup_{j=-\infty}^\infty(t_{j1},t_{j2})=\R$.
 We have $\overline D_0^{(1)}\times \R\subset\bigcup_{j=-\infty}^\infty\overline D_0^{(1)}\times [t_{j1},t_{j2}]$.
 It follows from Theorem \ref{theo:7.7}  that  for each  $j\in \Z$  there exists a diffeomorphism
 $\Psi_j$  on  $D_j^{(1)}\times[t_{j1},t_{j2}]$  and a gauge transformation  $c_j$  such  that  $\Psi_j=I$  and
 $c_j=1$  on  $\overline \Gamma_0\times[t_{j1},t_{j2}]$,
 and
 \begin{equation}																	\label{eq:7.26}
 c_j\circ \Psi_j^{-1}\circ L^{(2)}=L^{(1)}\ \ \mbox{in}\  \ \overline D_0^{(1)}\times[t_{j1},t_{j2}].
 \end{equation}
 In  (\ref{eq:7.26})  $\Psi_j$  is a diffeomorphism of  $\overline D_0^{(1)}\times[t_{j1},t_{j2}]$  onto an almost cylindrical 
domain $\overline D_0^{(2)}\times\{S_j^-(x_1,...,x_n)\leq x_0\leq S_j^+(x_1,...,x_n)\}$,
 where $x_0= S_j^\pm(x_1,...,x_n)$
 are space-like  surfaces,  $\Psi_j=I$  on  $\overline\Gamma_0\times[t_{j1},t_{j2}],\ |c_j(x)|=1$  for all  
 $x\in  \overline D_0^{(1)}\times[t_{j1},t_{j2}],\ c_j=1$
on  $\Gamma_0\times[t_{j1},t_{j2}]$.

We shall show  that
\begin{equation}																	\label{eq:7.27}
\Psi_j=\Psi_{j+1},\ \ c_j=c_{j+1}
\end{equation}
on  $\overline D_0^{(1)}\times[t_{j+1,1},t_{j2}]$ where  $[t_{j+1,1},t_{j2}]$  is the intersection  of  $[t_{j1},t_{j2}]$  and
$[t_{j+1,1},t_{j+1,2}]$.

Let  $[t_{j1},t_{j+1,2}]=[t_{j1},t_{j2}]\cup[t_{j+1,1},t_{j+1,2}]$.
Note  that the proof  of Theorem \ref{theo:7.7}   consists  of several steps,  each step  is dealing  with  the 
arbitrary  large  time interval.  Consider,   for  example,  Lemma  \ref{lma:7.1}.

The changes  of variables  (\ref{eq:7.4}),  (\ref{eq:7.5})  and  (\ref{eq:7.7})  are defined  on arbitrary  time interval  
$[t',t'']$.  Thus  they are defined  on  $[t_{j1},t_{j+1,2}]$.
Therefore  the maps  (\ref{eq:7.4}),  (\ref{eq:7.5}),   (\ref{eq:7.7})
are  defined  on  $[t_{j1},t_{j2}]$  and  $[t_{j+1,1},t_{j+1,2}]$  and they coincide  on 
$[t_{j+1,1},t_{j2}]=[t_{j1},t_{j2}]\cap[t_{j+1,1},t_{j+1,2}]$.

Also  using  the  extension theorems  of  [Hi],  \S 8   we can find  the extension $\tilde\Phi_2\cap[t_{j1},t_{j+1,2}]$ 
of $\Phi_2\cap[t_{j1},t_{j+1,2}]$ 
 (cf.  (\ref{eq:7.8}))  in two steps:  First  extend  $\tilde\Phi_2\cap[t_{j1},t_{j+1,1}]$
 and then extend  continuously
  $\Phi_2\cap[t_{j+1,1},t_{j+1,2}]$ 
to get  a continuous  extension  of $\Phi_2\cap[t_{j1},t_{j+1,2}]$.
This way  we show  that  (\ref{eq:7.27})   holds  for some  $j\in \Z$.
Therefore  starting   with  $\Psi_0$  on  $[t_{01},t_{02}]$   we  can  construct  $\Psi_1$  on  $[t_{11},t_{12}]$  such
that   $\Psi_0=\Psi_1$  on  $[t_{01},t_{02}]\cap [t_{11},t_{12}]$.  Continuing  this construction  we get
(\ref{eq:7.27})  for any  $j\in \Z$.

Let  $\Psi=\Psi_j, c=c_j$  on  $[t_{j1},t_{j2}],\forall j\in \Z$.   Then  $\Psi$  is a proper  diffeomorphism  
of  $\overline D_0^{(1)}\times \R$  onto  $\overline D_0^{(2)}\times \R,  \Psi=I$  on  $\Gamma\times\R$   and
$$
c\circ \Psi^{-1}\circ L^{(2)}=L^{(1)} \ \ \mbox{on}\ \ \overline D_0^{(1)}\times \R.
$$

\end{document}